\documentclass[11pt]{amsart}
\usepackage{txfonts}
\usepackage{amsmath,yhmath}
\usepackage{amsxtra}
\usepackage{amscd}
\usepackage{amsthm}
\usepackage{amsfonts}
\usepackage{amssymb}
\usepackage{eucal}
\usepackage{pmgraph}
\usepackage{stmaryrd}
\usepackage[usenames,dvipsnames]{color}
\usepackage{graphics}
\newcommand{\newcom}{\newcommand}

\newcom{\al}{\alpha}
\newcom{\be}{\beta}
\newcom{\eps}{\epsilon}
\newcom{\del}{\delta}
\newcom{\D}{\Delta}
\newcom{\g}{\gamma}
\newcom{\G}{\Gamma}
\newcom{\ka}{\kappa}
\newcom{\Lam}{\Lambda}
\newcom{\lam}{\lambda}
\newcom{\Om}{\Omega}
\newcom{\om}{\omega}
\newcom{\Si}{\Sigma}
\newcom{\si}{\sigma}
\newcom{\tht}{\theta}
\newcom{\dtri}{\nabla}
\newcom{\td}{\tilde}
\newcom{\tri}{\triangle}
\newcom{\oo}{\infty}
\newcom{\vphi}{\varphi}
\newcom{\cA}{{\mathcal A}}
\newcom{\cB}{{\mathcal B}}
\newcom{\cC}{{\mathcal C}}
\newcom{\cD}{{\mathcal D}}
\newcom{\cF}{{\mathcal F}}

\newcom{\cL}{{\mathcal L}}

\newcom{\cK}{{\mathcal K}}
\newcom{\cP}{{\mathcal P}}
\newcom{\cR}{{\mathcal R}}
\newcom{\cS}{{\mathcal S}}
\newcom{\cT}{{\mathcal T}}
\newcom{\cU}{{\mathcal U}}
\newcom{\cX}{{\mathcal X}}
\newcom{\cY}{{\mathcal Y}}
\newcom{\cN}{{\mathcal N}}
\newcom{\cH}{{\mathcal H}}

\newcom{\mN}{{\mathfrak N}}

\newcom{\R}{\Bbb R}
\newcom{\N}{\Bbb N}
\newcom{\Z}{\Bbb Z}
\newcom{\C}{\Bbb C}
\newcom{\E}{\Bbb E}

\newcom{\bx}{\bar x}
\newcom{\bz}{\bar z}
\newcom{\tx}{\tilde x}
\newcom{\tz}{\tilde z}
\newcom{\f}{\frac}
\newcom{\di}{\displaystyle\int}
\newcom{\ds}{\displaystyle\sum}
\newcom{\dl}{\displaystyle\lim}
\newcom{\ov}{\overline}
\newcom{\sset}{\subset}
\newcom{\wt}{\widetilde}
\newcom{\pa}{\partial}
\newcom{\na}{\nabla}
\newcom{\co}{\cdot}
\newcom{\suml}{\sum\limits}
\newcom{\supl}{\sup\limits}
\newcom{\intl}{\int\limits}
\newcom{\infl}{\inf\limits}
\newcom{\disp}{\displaystyle}
\newcom{\non}{\nonumber}
\newcom{\no}{\noindent}
\newcom{\QED}{$\square$}

\def\ef{\hphantom{MM}\hfill\llap{$\square$}\goodbreak}

\def\dive{\mathop{\rm div}\nolimits}
\def\curl{\mathop{\rm curl}\nolimits}

\newtheorem{athm}{\bf \t}[section]
\newenvironment{thm} [1] {\def\t{#1}\begin{athm} \bf \rm} {\end {athm}}
\newcom{\bthm}{\begin{thm}}\newcom{\ethm}{\end{thm}}

\newcom{\beq}{\begin{equation}}
\newcom{\eeq}{\end{equation}}
\newcom{\ben}{\begin{eqnarray}}
\newcom{\een}{\end{eqnarray}}
\newcom{\beno}{\begin{eqnarray*}}
\newcom{\eeno}{\end{eqnarray*}}

\numberwithin{equation}{section}

\begin{document}

\title{Water waves problem with surface tension in a corner domain II: the local well-posedness}

\author{Mei Ming}
\address{School of Mathematics and Statistics,Yunnan University, Kunming 615000, P. R. China}
\address{School of Mathematics, Sun Yat-sen University, Guangzhou 510275, P. R. China}
\email{mingmei@ynu.edu.cn}

\author{Chao Wang}
\address{School of  Mathematical Science, Peking University, Beijing 100871, P. R. China}
\email{wangchao@math.pku.edu.cn}

\date{Dec. 24th, 2018 }

\maketitle

 \begin{abstract}

Based on the a priori estimates in our previous work \cite{MW}, we continue to investigate the water-waves problem in a bounded two-dimensional corner domain  in this paper.  We prove the local well-posedness of the solution to the water-waves system when the contact angles are less than $\f{\pi}{16}$.

 \end{abstract}

\tableofcontents

\section{Introduction}
We consider the irrotational incompressible water-waves problem in a two-dimensional bounded corner domain $\Om_t$ with an  upper free surface $\Gamma_t$ and a fixed  bottom $\Gamma_b$.  This domain contains two contact points $p_{l}, p_{r}$ with contact angles $\om_l, \om_r$, which are the intersections of $\G_t$ and $\G_b$.

The water-waves problem on $\Om_t$ can be expressed as the  system of  velocity $v$ and pressure $P$:
 \[
\mbox{(WW)}\qquad \left\{
\begin{array}{l}
\pa_t v+v\cdot \na v=-\na P+{\bf g}, \\
\dive v=0,\quad \curl v=0\qquad \hbox{on}\quad \Om_t\\
P|_{\Gamma_t}=\sigma \kappa,\\
\pa_t+v\cdot \na \quad \textrm{is tangent to}\quad \bigcup_t\Gamma_t,\\
v\cdot N_b |_{\Gamma_b}=0,\\
\beta_c v_i=\sigma(\cos{\om_s}-\cos{\om_i})\qquad\hbox{at}\quad p_i\ (i=l,r),
\end{array}
\right.
\]
where  $\kappa$  is the mean curvature of the free surface, $\sigma$ is the  surface tension coefficient, $\om_i$ ($i=l, r$) are the contact angles between $\Gamma_t,  \Gamma_b$, and $ g$ is the gravity coefficient with ${\bf g}=-g {\bf e_z}$ the gravity vector.  Moreover, we denote by $v_i$ the upward tangential component of the velocity at the corner points along $\G_b$:
\[
v_l=-v\cdot\tau_b\quad\hbox{at}\  p_{l},\quad\hbox{and}\quad v_r=v\cdot\tau_b\quad\hbox{at}\  p_{r}.
\]

The last condition in \mbox{(WW)} describes the motion at the contact points, which was studied in \cite{RE} and has been used in our previous work \cite{MW}.  Here the stationary contact angle $\om_s$ is a physical constant depending  on the materials of the bottom and the fluid, and $\beta_c$ denotes the  effective friction coefficient. In fact, this condition tells us that slip velocity is dominated by the unbalanced Young stress,   which  is an effective variation of Young's law (1805) for  stationary contact angles \cite{Young}.
This kind of conditions are very common and widely discussed, see \cite{BEIMR, CDA, SA, GL}. One can see in our previous work \cite{MW} that, there is some dissipation corresponding to this condition at the contact points naturally.  Mathematically,  this condition turns out to be some kind of Lopatinsky condition,  and it is necessary for solving the linear system for the iteration as well as for proving the energy estimates. 
Moreover, a similar condition was used in \cite{GT} for the Stokes flow.

\includegraphics*[-50,0][300,120]{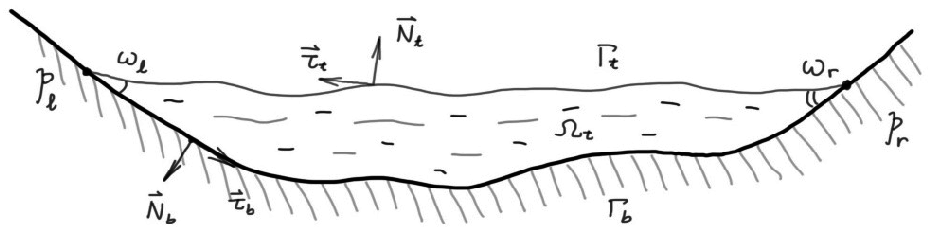}

 \subsection{Some known results}
 Let us  recall some previous works on the well-posedness for the water-waves problems. When  we say  ``classical " water-waves problems, we refer to  the water-waves problems  with  a smooth free surface, and the fluid boundaries satisfy $\Gamma_t\cap \Gamma_b= \varnothing$. There is a rich literature on the classical water-waves problems.

Firstly, we recall the results on the local well-posedness for the classical water-waves problems. To begin with, we have a quick review on some works about the irrotational case. Some early works such as Nalimov \cite{Na}, Yosihara \cite{Yo1, Yo2} and Craig \cite{Craig} established  the local well-posedness with small data  in two-dimensional case. A breakthrough is done by  Wu \cite{Wu1, Wu2} which removed the smallness condition and proved that the Taylor sign condition
 \beno
 -\na_{N_t}P|_{\Gamma_t} \geq c_0>0
 \eeno
always holds as long as $\Gamma_t$ is not self-intersection. Later on,   some more different methods were applied to prove the local well-posedness. 
Iguchi \cite{Iguchi} and Ambrose \cite{AM} studied the local well-posedness in two-dimensional case respectively.
In \cite{Lannes}, Lannes   proved the finite-depth case under Eulerian coordinates.  Later, Ming and Zhang \cite{MZ} generalized Lannes's paper to the case with surface tension.  Alazard, Burq and Zuily in \cite{ABZ1, ABZ2, ABZ3} used the tools of paradifferential operators to prove the local well-posdenss in a low-regularity case. Moreover, Alvarez-Samaniego and Lannes \cite{AL} considered the large-time existence for the problem under the shallow-water regime.

Concerning the rotational case, there are also  many works  on the local well-posedness using various methods.
Christodoulou and Lindblad \cite{CL} were the first to prove a priori estimates based on the geometry of the moving domain,  and later  Lindblad \cite{Lin} proved the existence of solutions using Nash-Moser iteration. Coutand and Shkoller \cite{CS} proved the local well-posedness under Lagrangian coordinates. Zhang and Zhang \cite{ZZ} used the Clifford analysis introduced by Wu \cite{Wu2} to solve the problem. Shatah and Zeng \cite{SZ, SZ2} treated the problem in a geometric way,  where they used the equation of the mean curvature. Meanwhile, a similar geometric approach had also been used by Beyer and G\"unther \cite{BG1, BG2} to study the irrotational problem for some star-shaped domains. Recently, Wang, Zhang, Zhao and Zheng \cite{WZZZ} proved the local well-posedness in low-regularity case. For more results on the local well-posedness, the readers can check the book by Lannes \cite{LannesBook}, and Iguchi, Tanaka and Tani \cite{Ig-Ta}, Ogawa and Tani \cite{OT1, OT2}, Schweizer \cite{Sch}, Ambrose and Masmoudi \cite{AM1, AM2}  etc..

\medskip

For the global well-posedness, the first result was given by Wu \cite{Wu4} which proved  the almost-global existence for the gravity problem in two dimensions. Later, Wu \cite{Wu5} and  Germain, Masmoudi and Shatah \cite{GMS} showed the global existence of gravity waves in three dimensions respectively by different methods.  Moreover, Alazard and Delort \cite{AD} and Ionescu and Pusateri \cite{IP} studied the global regularity for  gravitational water-waves systems in two dimensions independently. Recently, Hunter, Ifrim and Tataru \cite{HIT1,HIT2, HIT3} used the conformal-mapping method to give another proof of the global existence for the gravitational problem in two dimensions. For more results on the global well-posedness, readers can check \cite{Deng, Wang1} and their references.

\medskip

Compared to the classical water-waves case, when we say ``non-smooth" water-waves problems, we mean that there are contact points on the fluid boundaries i.e. $\Gamma_t\cap \Gamma_b\neq \varnothing$,  or the free surface is not smooth.   In fact, theoretical research on this field only started several years ago and there remains a lot of open problems.  Alazard, Burq and Zuily \cite{ABZ3} proved the local well-posedness for the special case when the contact angle is equal to $\pi/2$. In this case, they used symmetrizing and periodizing methods to turn this problem into a classical case. Later, Kinsey and Wu proved the local well-posedness for the two-dimensional water waves with angled crests when the wall is vertical, see \cite{WuK,Wu3}. Recently, de Poyferr\'e \cite{Poyferre} gave a priori estimates for the water-waves problem in a bounded corner domain without  surface tension under the assumption of  small contact angles. Meanwhile, under the assumption of small contact angles, the authors  proved a priori estimates for the water-waves problem in a corner domain with surface tension, see \cite{MW}. For both  the two results \cite{Poyferre, MW}, one important observation is that  small contact angles can prevent the appearance of  singularities from the corners.

On the other hand,  Lannes and M\'etivier \cite{LM} solved the local well-posedness for the Green-Naghdi equations in a beach-type domain, which is a shallow-water model of the water-waves problem. Lannes \cite{Lannes1} addressed the floating-body problem and proposed a new formulation of the water-waves problem that can be easily  generalized in order to take into account the presence of a floating body. Very recently, Lannes and Iguchi \cite{LI} proved some sharp results for initial boundary value problem with a free boundary arising in wave-structure interaction, and it contains the floating problem in the shallow-water regime. Besides, Guo and Tice considered  a priori estimates for the contact line problem in case of the stokes equations, see \cite{GT}. Later, Tice and Zheng proved the local well-posedness of the contact line problem in 2D Stokes flow, see \cite{TZ}.

In the end, we also mention some results concerning  geometric singularities on the free surfaces for the water-waves problems. In \cite{CCFG}, the authors showed the existence of a wave which is  given initially as the graph of a function and then can overturn at a later time. Later on, the authors in \cite{CCFGG} proved the existence of some  ``splash" singularities. Moreover, this result was extended in \cite{CS2} to three-dimensional case and some other models.

 \medskip

 \subsection{Main results and ideas}
 In this paper, we  prove the local well-posedness of system \mbox{(WW)} on a bounded two-dimensional  corner domain, which is based on our previous work \cite{MW}. The following theorem states our main result:
 \bthm{Theorem}\label{main theorem}
{\it 
Assume that the initial data $(\G_0,v_0) \in H^{8.5}\times H^{7.5}(\Om_0)$ and the initial contact angles $\om_{i0}\in(0, \pi/{16})$ for $i=l, r$.  
When the compatibility conditions at $t=0$, namely
\ben\label{eq:com cond}
\beta_c\pa_t^k v_i(0)=\sigma\pa_t^k\big(\cos{\om_s}-\cos{\om_i(0)}\big)\qquad\hbox{at}\quad p_i \,(i=l,r),\  k=0,1,2,3,4
\een
are satisfied,
there exists a small constant $T>0$ depending on the initial data such that system \mbox{(WW)} has a unique solution $(\G_t, v) \in C([0, T];H^{8.5})\times C([0, T];H^{7.5}(\Om_t))$. Moreover, the solution $(\G_t, v)$ is locally well-posed.
}
\ethm

\bthm{Remark}
{\it According to our previous paper \cite{MW} and compactness arguments, we can get  $(\G_t, v) \in C([0, T];H^{s})\times C([0, T];H^{s-1}(\Om_t))$ with $4\leq s\leq 8.5$.
}
\ethm

Similarly as in \cite{SZ, SZ2},  the pressure $P$ in \mbox{(WW)} is regained by the velocity $v$ and the mean curvature $\ka$. In fact, we  decompose the pressure into two parts:
\[
P=\si \ka_\cH+P_{v, v},
\]
where the first part  is the harmonic extension of $\ka$ and the second part $P_{v,v}$ is decided by $v$ and ${\bf g}$:
\[\left\{\begin{array}{ll}
\D P_{v,v}=-tr (\na v)^2\qquad\hbox{on}\quad\Om_t,\\
P_{v,v}|_{\G_t}=0,\qquad \na_{N_b} P_{v,v}|_{\G_b}=\na_v N_b\cdot v+N_b\cdot{\bf g}.
\end{array}\right.
\]
Therefore, as long as we have $(\G_t, v)$, the whole water-waves system \mbox{(WW)} is recovered immediately.

  \bthm{Remark}
 {\it  In this paper, we need smaller contact angles compared to \cite{MW}. The reason is that when we construct the Cauchy sequence, we need a higher regularity. To avoid the singularity from the corners, the range of the angles is decided by $P_{v,v}$ using Remark 5.20 \cite{MW1} for the related mixed-boundary elliptic problem.  The local well-posedness of system \mbox{(WW)} with general contact angles still remains an open problem.
  }\ethm

  \medskip

Now, we explain the main ideas of this paper.  Firstly, we need to choose a good formulation to construct approximate solutions. Inspired by \cite{SZ2}, we introduce a universal coordinate map $\Phi_{S_t}: \Gamma_*\to \Gamma_t$ which can reduce our system into a system defined on a fixed domain.  $d_{\Gamma_t}$ is used as the ``distance" between $\Gamma_*$ and $\Gamma_t$, where $\G_*$ is some reference upper surface.

Secondly, based on the mean curvature $\ka$ of $\G_t$, we introduce a new quantity on $\G_*$:
\[
\mN_a=\cN(\ka)\circ\Phi_{S_t}+a^3 d_{\G_t},
\]
 which is different from the modified mean curvature $\ka_a$ in \cite{SZ2}. Here $\cN$ stands for the Dirichlet-Neumann operator. We derive the evolution equation and the boundary conditions for $\mN_a$ from $\mbox{(WW)}$:
\beno
\left\{\begin{array}{ll}
D_{t*}^2 \mN_a+\sigma \cA(d_{\G_t})\mN_a=R_0\qquad\hbox{on}\quad \G_*,\\
D_{t*}\cA(d_{\G_t}) \mN_a+\f{\sigma^2}{\beta_c}( \sin \om_{i})^2 \circ\Phi_{S_t} (\na_{\tau_t}(\cA(d_{\G_t})\mN_a\circ \Phi_{S_t}^{-1}))\circ\Phi_{S_t}= R_{c,i}\quad \textrm{at} \quad p_{i*}\ (i=l,r),
\end{array}
\right.
\eeno
 where $R_0, R_{c,i}$ are remainder terms defined by $d_{\Gamma_t}, \pa_t d_{\Gamma_t}$ and $v$, and the third-order operator $\cA(d_{\G_t})$ is defined as
\[
  \cA (d_{\G_t}) f=-\big( \cN \Delta_{\Gamma_t} (f\circ\Phi^{-1}_{S_t})\big)\circ \Phi_{S_t},\qquad\hbox{for some $f$ on $\G_*$}.
\]
As mentioned before, the boundary conditions at $p_i$ play a key role in the energy estimates and in the iteration, which have been used in a different version for the equation of $J=\na\ka_\cH$ in \cite{MW}. The part involving the boundary conditions in our papers is completely new compared to the classical water waves or the other works on non-smooth water waves.

Moreover,  the velocity $v$  is  recovered by $\pa_t d_{\G_t}$ in the iteration, so the  energy estimates and the iteration depend on the free surface $d_{\G_t}$ and its boundary conditions. The boundary conditions for $d_{\G_t}$ (i.e. information of the contact points)  take the form of ODE:
 \[
 d''_i(t)=\mathfrak B_i,\quad i=l, r,
 \]
where $d_i=d_{\G_t}\big|_{p_i}$,  and  $\mathfrak B_i$ depends on $d_{\G_t}, \pa_t d_{\G_t}, v$ as well.  These conditions come naturally from the definition of the velocity on $\G_t$ and Euler equation, which are very important when one wants to retrieve $d_{\G_t}$ from $\mN_a$.

As a result, the above system of $(\mN_a, d_l, d_r)$ provides a closed system to construct  approximate solutions. Due to the presence of  the contact points, we do not expect the system is smooth even if the initial data is smooth enough. In  \cite{MW}, we need small contact angels to avoid singularities from the corresponding elliptic systems.  In this paper, we assume that the contact angels are smaller than ${\pi}/{16}$,  which ensures that a related  mixed-boundary
elliptic system (see Lemma \ref{existence of D system}) has a solution in $H^9(\Om_t)$. Compared to the regularity considered in \cite{ABZ1, ABZ2, ABZ3}, $H^9(\Om_t)$ is a much higher regularity, and a lower regularity is still desirable.

Even with small contact angles, the choice of proper energy functionals is still made very carefully.  
To prove the energy estimates of the linear system for the iteration, we need to use the material derivative $D_t$  more often than $\na_{\tau_t}$. We choose the following energy and dissipation functionals
\[
\begin{split}
E_h(t, \bar f, \pa_t \bar f)&= \big\|\na_{\tau_t}\cN\D_{\G_t}D_t\bar f\big\|^2_{L^2(\G_t)}+\si\big\|\na \cH(\D_{\G_t}\cN\D_{\G_t}\bar f) \big\|^2_{L^2(\Om_t)}+\|D_t\bar f\|^2_{L^2(\Gamma_t)} +\|\bar f\|^2_{L^2(\Gamma_t)},\\
F_{h, i}(t, \bar f, \pa_t \bar f)&= ( \sin \om_i)^2\big|\na_{\tau_t}\cN\D_{\G_t}D_t\bar f\big|^2\qquad \textrm{at}\quad p_i\ (i=l, r).
\end{split}
\]
Notice that the dissipation only takes place at the contact points, and one can check our previous work \cite{MW} for more details.

\medskip

In the end, we emphasize the differences between our paper with \cite{SZ2}. We use the geometric approach introduced by \cite{SZ2}. Compared to \cite{SZ, SZ2},  some new difficulties appear due to the presence of the corners. Firstly, during the construction of approximate solutions, we choose to use the equation for a new quantity $\mN_a$,  while the modified mean curvature $\kappa_a=\ka\circ\Phi_{S_t}+a^2 d_{\G_t}$ is used in \cite{SZ2}. The reason of using  $\mN_a$ is that it's more convenient to derive the boundary conditions for  $\mN_a$ at the contact points,  and we do not have the information for $\kappa_a$ at the same time. Besides, if we choose $\cN(\kappa_a)$ instead of $\mN_a$, we need to maintain 
% $\dive \na \kappa_{\cH}=0$ and 
$\int_{\Gamma_t}\cN(\kappa_a)\circ\Phi^{-1}_{S_t}ds=0$ in the iteration, which  makes the iteration much  more complicated. But choosing $\mN_a$ we  needs no restriction. 
%Moreover, the equation of $\mN_a$ can recover the mean curvature by an elliptic equation. 
Secondly, when we recover the free surface from $\mN_a$,  the boundary conditions $d_i=d_{\G_t}|_{p_i} (i=l, r)$ are needed essentially  to solve the related elliptic equation, see Proposition \ref{d and  N ka}.  The system of $(\mN_a, d_l, d_r)$ together makes sure that our iteration sequence  converges and goes back to the solution to system \mbox{(WW)}. Thirdly, the definition of the energy functionals as well as the dissipations are totally different in our paper, and the details involving the contact points in the energy estimates are completely new.

\subsection{Organization of the paper}
In Section 2, we give some useful lemmas.
In Section 3 the free surfaces and the domains are defined.  In Section 4, we  recover the velocity from the free surface. Meanwhile, we also give the equivalent formulation of the problem. Section 5 deals with the existence of the solution to the linear problem and proves  higher-order energy estimates. In Section 6, we use an iteration scheme to finish the proof for the local well-posedness.

\subsection{Notations}

\noindent - $\Om_*$ is the reference domain with the boundary $S_*=\G_*\cup \G_{b*}$. Here $\G_*$ is the upper boundary and $\G_{b*}$ is the fixed bottom. $\tau_*$ is the unit tangential vector of $\G_*$.\\
\noindent - $X$ and $p$ are both used to denote a point in $\Om_t$ or sometimes $\Om_*$.\\
\noindent - The entire fixed bottom is denoted as $\G_{fix}$.\\
\noindent - $\Pi$: the second fundamental form where $\Pi(w)=\na_w N_t\in T_X\G_t$ for $w\in T_X\G_t$.\\
\noindent- $\Pi(v,w)$ denotes $\Pi(v)\cdot w$. \\
\noindent - $\kappa=tr \Pi=\na_{\tau_t}N_t\cdot \tau_t$ is the mean curvature.\\
\noindent - $\cH(f)$ or $f_\cH$ is the harmonic extension for some function $f$ on $\G_t$, which is defined by the elliptic system
\[
\left\{\begin{array}{ll}
\D \cH(f)=0\qquad\hbox{on}\quad \Om_t,\\
\cH(f)|_{\G_t}=f,\quad \na_{N_b}\cH(f)|_{\G_b}=0.
\end{array}\right.
\]
\noindent - $\top$ denotes the tangential component of a vector.\\
\noindent - $\D^{-1}(h,g)$ is defined as  the solution $u$ to the system
\[
\left\{\begin{array}{ll}
\D u=h\qquad \hbox{on}\quad \Om_t,\\
u|_{\G_t}=0,\qquad \na_{N_b}u|_{\G_b}=g.
\end{array}
\right.
\]
\noindent - $\bar\cD$ is the covariant derivative on $T\G_t$ under Lagrangian coordinates, and $\cD$ is the corresponding derivative under Eulerian coordinates.\\ 
%\noindent - $\cD^2f(\tau_1, \tau_2)=D^2f(\tau_1, \tau_2)-(\Pi(\tau_1)\cdot\tau_2)\na_{N_t} f$  for any two %vectors $\tau_1,\tau_2\in T_X \G_t$.\\
\noindent - $M^*$ denotes the transport of a matrix $M$.\\
\noindent - $w^\perp$ on $\G_t$:   the normal component $(w\cdot N_t)\,N_t$.  \\
\noindent - $w^\top$ on $\G_t$: the tangential component $(w\cdot \tau_t)\,\tau_t$. Sometimes we also use $w^\top$ on $\G_b$ with a similar definition.\\
\noindent - $d_i=d_i(t)$ stands for the value of $d_{\G_t}$ at the contact points $p_{i}$ $(i=l ,r)$.\\
\noindent -  Sometimes we need to identify the signs between the left and the right contact point in boundary conditions; We always take $``+"$ for   $i=l$,  and $``-"$ for $i=r$.
\\
\noindent - $F=F(u_1, u_2,\dots, u_m)$ denotes that the higher-order terms in   function $F$ are $u_1, u_2,\dots, u_m$.\\
\noindent - $C=C(\|u_1\|, \|u_2\|,\dots, \|u_m\|)$ denotes a constant $C$ in the form of a polynomial  of some norms for  $u_1, u_2,\dots, u_m$.\\
\noindent - $D_{t*}=\pa_t+v^*\cdot \na$ is defined by $v^*=D\Phi^{-1}_{S_t}(v\circ\Phi_{S_t}-\pa_td_{\G_t}\mu)$ on $\G_*$.\\
\noindent - $D_\tau=\pa_\tau+v_\tau\cdot \na$ is defined on $\G_t$ and $\Om_t$ with some parameter $\tau$, and $v_\tau$ is induced by $\Phi_{S_t}$ with $v_\tau=(\pa_\tau d_{\G_t}\mu)\circ \Phi^{-1}_{S_t}$ on $\G_t$.\\
\noindent - $\lceil s\rceil=s$ when $s>0$ is an integer, and $\lceil s\rceil=m+1$ when $s=m+\epsilon$ for some  $\epsilon\in (0,1)$ and $m\in \N$.

\section{Some preliminaries}

Firstly, we present some useful lemmas on elliptic systems adjusted from \cite{MW1, MW}.

\bthm{Lemma}\label{existence of D system}
{\it Assume that $2\leq s\leq 9$. Let $f\in H^{s-1/2}(\G_{*})$ and $g\in H^{s-1/2}(\G_{b*})$  satisfying
\[
f|_{p_{i*}}=g|_{p_{i*}},\ i=l, r.
\]
If the contact angles of $\Om_*$ are less than $\f\pi{2(\lceil s\rceil-1)}$, then the system
\[
\left\{\begin{array}{ll}
\D u=0\qquad\hbox{on}\quad \Om_*,\\
u\big|_{\G_{*}}=f,
\quad u\big|_{\G_{b*}}=g,
\end{array}\right.
\]
admits a unique solution $u\in H^s(\Om_*)$ satisfying the estimate
\[
\|u\|_{H^s(\Om_*)}\le C\big(\|f\|_{H^{s-1/2}(\G_{*})}+\|g\|_{H^{s-1/2}(\G_{b*})}\big),
\]
with the constant $C$ depending on $\Om_*$ and $s$.
}
\ethm
{\bf Proof}. Checking Theorem 4.7,  Proposition 5.19 and Remark 5.20  in \cite{MW1} for the Dirichlet problem, one can obtain the desired results. 
\ef

\medskip

 \bthm{Lemma}\label{embedding}(Lemma 5.8 \cite{MW})
{\it We have the following embeddings:
\[
\|u\|_{L^4(\Om_t)}\le C(\|\Gamma_t\|_{H^{5/2} })  \|u\|_{H^{1/2}(\Om_t)},
\]
for any $u\in H^{1/2}(\Om_t)$,  and
\[
\|u\|_{L^\infty(\Om_t)}\le C(\|\Gamma_t\|_{H^{5/2} }) \|u\|_{H^{s_1}(\Om_t)},
\]
for  any $u\in H^{s_1}(\Om_t)$ with $s_1=1+\epsilon$ ($\epsilon>0$ is a small constant). Moreover, for any $ f\in H^{s_2}(\G_t)$ with $s_2=1/2+\epsilon$, the embedding holds:
\[
\| f \|_{L^\infty(\G_t)}\le C(\|\Gamma_t\|_{H^{5/2} }) \|f\|_{H^{s_2}(\G_t)}.
\]
Similar result holds for the case of $\G_b$.
}
\ethm

\bthm{Lemma}\label{Harmonic extension H1 estimate}(Lemma 5.9 \cite{MW})
{\it Let $\cH(f)$ be the harmonic extension of a function $f\in H^{1/2}(\G_t)$.
 Then  one has $\cH(f)\in H^1(\Om_t)$ satisfying the following estimate
\[
\|\cH(f)\|_{H^1(\Om_t)}\le C(\|\Gamma_t\|_{H^{5/2} }) \|f\|_{H^{1/2}(\G_t)}.
\]
}
\ethm
Besides, we need to consider the following mixed boundary problem sometimes:
\beq\label{eq:elliptic}
\left\{
\begin{array}{l}
\Delta u=h\qquad \textrm{on}\quad \Om_t, \\
u|_{\Gamma_t}=f,\quad\quad \na_{N_b} u|_{\Gamma_b}=g.
\end{array}
\right.
\eeq
\bthm{Lemma}\label{est:elliptic} (Elliptic estimates)
{\it Let the contact angles $\om_i\in(0,\f\pi{2(\lceil s\rceil-1)})$ and $s\geq 2$. Then the following estimate holds for system \eqref{eq:elliptic}:
\beno
\|u\|_{H^{s}(\Om_t)}\leq C(\| \Gamma_t\|_{ H^{s-1/2}  })\big(\|h\|_{H^{s-2}(\Om_t)}+ \|f\|_{H^{s-1/2}(\Gamma_t)}+\|g \|_{H^{s-3/2}(\Gamma_b)}\big).
\eeno
}
\ethm
\begin{proof}
When $s\leq 4$, the case has been proved in Theorem 5.1 in \cite{MW}. For the higher-order estimates, we apply Proposition 5.19 in \cite{MW1}. When contact angles are less than $\f\pi{2(\lceil s\rceil-1)}$, there is no singular part in our elliptic estimate thanks to Remark 5.20 \cite{MW1}.
\end{proof}

Moreover, the trace theorem on $\G_t$, $\G_b$  is quoted directly from Theorem 5.3 \cite{MW}.
\bthm{Lemma}\label{trace thm PG} (Traces on $\G_t$ or $\G_b$)
{\it Let the integer $ l\in[0,s-1/2)$ with $s-l>1/2$, we define the map
\beno
u\to \{u, \na_{N_j} u,\dots \na_{N_j}^l u\}|_{\Gamma_j},
\eeno
 for $u\in \cD(\bar \Om_t)$ where $N_j$ is the unit outward normal vector on $\Gamma_j$ with $\Gamma_j$ taking $\Gamma_b$ or $\Gamma_t$. Then, the map has a unique continuous extension as an operator from
\beno
H^{s}(\Om_t)\quad \textrm{onto} \quad  \Pi^l_{k=0}H^{s-k-1/2}(\Gamma_j).
\eeno
Moreover, one has the estimate for $0\le k\le l$:
\[
\| \na^k_{N_b} u\|_{H^{s-k-1/2 }(\Gamma_b )}+\|  \na^k_{N_t}u\|_{H^{s-k-1/2 }(\Gamma_t )}\le    C(\|\Gamma_t\|_{H^{s-1/2} })\|u\|_{H^{s}(\Om_t)}.
\]
}
\ethm

\medskip

The Dirichlet-Neumann (D-N) operator $\cN$ is defined by 
\[
\cN(f)=\na_{N_t}\cH(f)\qquad\hbox{on}\quad \G_t
\]
for a function $f$ defined on $\G_t$, which is an important operator  in water waves. 
We would like to recall some useful properties of the Dirichlet-Neumann operator here.
\bthm{Lemma}\label{DN-1 operator}{\it Let  $s\in (1, 6.5]$ and $\om_i \in (0, \f{\pi}{2(\lceil s\rceil-1)})$ for $i=l, r$ .\\
\noindent (1) The D-N operator $\cN$ is an order-1 operator on $\G_t$:
\[
\|\cN(f)\|_{H^{s-1}(\G_t)}\le C\big(\|\G_t\|_{H^s}\big)\|f\|_{H^s(\G_t)};
\]
\noindent (2) The following estimate holds:
\[
\|f\|_{H^s(\G_t)}\le C\big(\|\G_t\|_{H^s}\big)\big(\big\|\cN(f)\big\|_{H^{s-1}(\G_t)}+\|f\|_{L^2(\G_t)}\big);
\]
Moreover,  when $\cN(f)\in \big(H^{1/2}(\G_t)\big)^*$ (the dual space of $H^{1/2}(\G_t)$),  one has
\[
\|f\|_{H^{1/2}(\G_t)}\le C\big(\|\G_t\|_{H^{2.5}}\big)\big(\big\|\cN(f)\big\|_{H^{1/2}(\G_t)^*}+\|f\|_{L^2(\G_t)}\big);
\]
\noindent (3) When $\int_{\G_t}f\,ds=0$ and $\int_{\G_t}g\,ds=0$, the inverse D-N operator $\cN^{-1}$ makes sense, i.e.
\[
\cN(f)=g\quad\hbox{implies}\ f=\cN^{-1} (g),
\]
and  the above two inequalities in (2) hold without the term $\|f\|_{L^2(\G_t)}$ on the right side.
}
\ethm
\begin{proof}
Although similar results as in this lemma have been proved in \cite{SZ, Poyferre},  we  still provide some details here.  The first result comes directly from Lemma \ref{trace thm PG}. For the second and third results, we need to consider the following system
\[
\left\{\begin{array}{ll}
\Delta u=0\qquad\hbox{on}\quad \Om_t,\\
\na_{N_t} u|_{\Gamma_t}= g, \quad \na_{N_b} u|_{\Gamma_b}= 0,
\end{array}\right.
\]
 where $f=u|_{\G_t}$ and  the compatibility condition holds:
\[
\int_{\G_t}g\,ds =0.
\]
 Since the elliptic estimates have been proved in Theorem 5.10 \cite{MW}, the proof lies in the existence of the variational solution $u\in \dot H^1(\Om_t)$.

\noindent Step 1: The existence of the solution $u$ and the $\big(H^{1/2}(\G_t)\big)^*$ case.   Defining the variation space
\[
\mathcal V=\Big\{v\in \dot H^1(\Om_t)\,\big|\,\int_{\Om_t}vdX<\infty \Big\},
\] 
one writes the variation equation for $u$ as
\[
\int_{\Om_t}\na u\cdot \na v dX=\int_{\G_t} g\,v\,ds,
\]
with $\forall v\in \mathcal V$.

Due to the compatibility condition, one has
\[
\int_{\G_t} g\,v\,ds=\int_{\G_t} g\,(v-\bar v)\,ds\le \|g\|_{H^{1/2}(\G_t)^*}\|v-\bar v\|_{H^{1/2}(\G_t)},
\]
where $\bar v=\big(\int_{\Om_t}dX\big)^{-1}\int_{\Om_t}vdX$.
Applying Lemma \ref{trace thm PG} and Poincar\'e's inequality, one obtains
\[
\int_{\G_t} g\,v\,ds\le C(\|\G_t\|_{H^{2.5}})\|g\|_{H^{1/2}(\G_t)^*}\|v\|_{\dot H^1(\Om_t)},
\]
as long as $g\in \big(H^{1/2}(\G_t)\big)^*$.

Consequently, applying Lax-Milgram's theorem, one concludes that there exists a variation solution $u\in \mathcal V$ to the Neumann problem satisfying the estimate
\[
\|\na u\|_{L^2(\Om_t)}\le C(\|\G_t\|_{H^{2.5}})\|g\|_{H^{1/2}(\G_t)^*}.
\]
 
On the other hand, since $u=\cH(f)$, one knows directly that
\[
\|\cH(f)\|_{L^2(\Om_t)}\le C \big(\|\na \cH(f)\|_{L^2(\Om_t)}+\|f\|_{L^2(\G_t)}\big),
\]
where the constant $C$ depends on the size of the domain $\Om_t$. Combining these two inequalities above with Lemma \ref{trace thm PG}, one derives the second estimate in (2).

In order to prove (3), when $\int_{\G_t} f\,ds=0$, one has
\[
\|\cH(f)\|_{L^2(\Om_t)}=\|\cH(f)-\bar f\|_{L^2(\Om_t)}\le C\|\na\cH(f)\|_{L^2(\Om_t)},
\]
where $\bar f=\big(\int_{\G_t}ds\big)^{-1}\int_{\G_t}f\,ds$. As a result,
one derives by Lemma \ref{trace thm PG} that
\[
\|f\|_{H^{1/2}(\G_t)}\le C(\|\G_t\|_{H^{2.5}})\|\cN(f)\|_{H^{1/2}(\G_t)^*},
\]
which implies that the D-N operator $\cN$ is invertible in this case and the lower-order estimate in (3) follows.

\medskip
\noindent Step 2: The higher-order estimate for (2) and (3). The estimates  follow from Step 1 and Theorem 5.10 \cite{MW}, as long as one checks carefully from Remark 5.20 \cite{MW1} or \cite{PG1} that the contact angle $\om\in (0, \f{\pi}{2(s-1)})$ for an integer $s$ under the required regularity of this lemma.

\end{proof}

\section{Definitions of surfaces and domains}

In this section, we firstly define a coordinates system based on a reference domain $\Om_*$, and then we construct surfaces according to the coordinates system. %Our formulation here follows from \cite{SZ2}.

To begin with, we fix a reference domain $\Om_*$, where $\Om_*$ can be taken as the initial domain $\Om_0$. The boundary of the reference domain $\Om_*$ is denoted by $S_*$, which contains two parts indeed: The upper surface $\G_{*}$ and the bottom $\G_{b*}$. The corresponding contact points are  noted as  $p_{i*}$ $(i=l, r)$ with contact angles 
\[
\om_{i*}\in (0, \pi/16).
\]
 The unit outward normal vectors and tangent vector are denoted by $N_{*},\,N_{b*}$ and $\tau_{*},\,\tau_{b*}$ accordingly.

\subsection{Definition for surfaces}

Since the domain will be fixed if the boundary is fixed, we consider about how to define the boundary, or the upper free surface. We use some oblique coordinates system on $\G_*$ to define surfaces $S$ near $S_*$. In fact, these surfaces will be set to be in a  neighborhood of $S_*$.

Firstly, we introduce a unit upward vector field $\mu\in H^s(\G_{*}, \cS^1)$ for some large $s$ such that
\[
\mu\cdot N_{*}\ge c_0\qquad\hbox{on}\quad \G_{*},\qquad\hbox{and}\quad\mu|_{p_{l*}}=-\tau_{b*}|_{p_{l*}},
\  \mu|_{p_{r*}}=\tau_{b*}|_{p_{r*}}
\]
with some fixed constant $c_0\in(0,1)$. Note that  this condition holds at $p_{l*}, p_{r*}$ since the contact angle stays in $(0,\pi/2)$ in this paper.

From the implicit function theorem, there exists a small constant $d_0>0$  such that the map
 \[
 \Phi: \G_*\times[-d_0,d_0]\rightarrow \R^2\qquad\hbox{where}\quad
\Phi(p, d) \triangleq p+d\,\mu(p)
 \]
is an $H^s$ diffeomorphism  from its domain to a neighborhood of $\G_*$. As a result, this coordinate system identifies each upper surface $\G$ close to $\G_*$ with a unique function $d_{\G}: \G_*\rightarrow  \R$.    Plugging $d_{\G}$ into $\Phi$, one writes
\[
\Phi_S: \ \G_*\rightarrow \G\subset \R^2,\qquad\hbox{where}\quad \Phi_S(p)=p+d_{\G}(p)\mu(p).
\]
Sometimes, the function $d_{\G}(p)$ is also used as the expression of the upper surface $\G$. Moreover, we extend $\Phi_S$ to be defined on the whole boundary $S_*$. In fact, let
\[
\Phi_S: \ \G_{b*}\rightarrow \G_{b}\subset \G_{fix}\qquad \hbox{satisfying}\quad \Phi_S\in H^s(\G_{b*},\G_b),\  \Phi_{S}(p_{i*})=p_{i*}+d_{\G}(p_{i*})\mu(p_{i*}) \ (i=l, r).
\]
Notice that $\G_{b*}$ and $\G_{b}=\Phi_S(\G_{b*})$ are both parts from the entire fixed bottom $\G_{fix}$.

Consequently, we denote  the surface of the domain by
\[
S=\G\cup \G_b,
\]
which is defined by $\Phi_S$. One can see that in our case, as long as we know the free upper surface $\G$, the whole boundary $S$ and  the domain $\Om$ are fixed. So all we need is to concentrate on the upper free surface $\G$.

\medskip
We are going to consider free surfaces varying near the reference surface $\G_*$  in the following set.
\bthm{Definition}\label{Lambda neighborhood}
{\it Let $\delta>0$ and $s\in(1.5, 8.5]$.  One defines the set of upper free surfaces
\[
\Lam(S_*,s,\del, {\pi}/{16})\triangleq \{\, \G\,\big|\,  \|d_{\G}\|_{H^s(\G_{*})} <\del,\ \|\Phi_S-Id_{S_*}\|_{H^s(\G_{b*})}\le \del, \, \om_i\in(0, \pi/16), i=l, r \}
\]
as a neighborhood of $\G_*$.
}
\ethm
When the constant $\del$ is taken small enough,  for any $\G\in \Lam(S_*,s,\del, {\pi}/{16})$,  $\Phi_S$ is a diffeomorphism both in $H^s(\G_{*},\G)$ and $H^s(\G_{b*},\G_b)$, and the contact angles $\om_i$ lie in $(0, \pi/16)$. Moreover, from the definition of $\Phi_S$, one can see that the norm $\|\Phi_S-Id_{S_*}\|_{H^s(\G_{b*})}$ on the bottom can be controlled by $d_{\G}$ as well.

\subsection{Harmonic coordinates}
Let
\[
\cT_S:\ \Om_*\rightarrow \Om\qquad\hbox{with}\quad \cT_S=\cH_*\big(\Phi_S-Id_{S_*}\big)+Id.
\]
Here $\cH_*\big(\Phi_S-Id_{S_*}\big)$ is the harmonic extension of $\Phi_S-Id_{S_*}$ satisfying the system with Dirichlet boundary conditions
\beq\label{H* system}
\left\{\begin{array}{ll}
\D \cH_*\big(\Phi_S-Id_{S_*}\big)=0\qquad\hbox{on}\quad \Om_*,\\
\cH_*\big(\Phi_S-Id_{S_*}\big)\big|_{\G_{*}}=d_{\G}\mu,
\quad \cH_*\big(\Phi_S-Id_{S_*}\big)\big|_{\G_{b*}}=\Phi_S|_{\G_{b*}}-Id_{\G_{b*}}.
\end{array}\right.
\eeq
One can see immediately from this definition that, the boundary of $\Om$ is the surface $S=\G\cup \G_b$.

Recalling the definition for $\Phi_S$, one knows that the compatibility condition for the Dirichlet boundary conditions is satisfied, i.e.
\[
d_{\G}\mu(p_{i*})=\big(\Phi_S|_{\G_{b*}}-Id_{\G_{b*}}\big)(p_{i*}) \ (i=l, r).
\]

%To consider about the regularity for $\cH_*(\Phi_S-Id_{S_*})$, one needs the following lemma.
%\bthm{Lemma}\label{existence of D system}
%{\it Assume that $2\leq s\leq 9$. Let $f\in H^{s-1/2}(\G_{t*})$ and $g\in H^{s-1/2}(\G_{b*})$  satisfying
%\[
%f|_{p_{i*}}=g|_{p_{i*}},\ i=l, r.
%\]
%Then the system
%\[
%\left\{\begin{array}{ll}
%\D u=0\qquad\hbox{on}\quad \Om_*,\\
%u\big|_{\G_{t*}}=f,
%\quad u\big|_{\G_{b*}}=g.
%\end{array}\right.
%\]
%admits a unique solution $u\in H^s(\Om_*)$ satisfying the estimate
%\[
%\|u\|_{H^s(\Om_*)}\le C\big(\|f\|_{H^{s-1/2}(\G_{t*})}+\|g\|_{H^{s-1/2}(\G_{b*})}\big),
%\]
%with the constant $C$ depending on $\Om_*$ and $s$.
%}
%\ethm
%{\bf Proof}. Checking The Proposition 5.19 and Remark 5.20  in \cite{MW1} for the Dirichlet problem, one can obtain the desired results.
%\ef

Applying Lemma \ref{existence of D system} on $\cH_*\big(\Phi_S-Id_{S_*}\big)$, one finds immediately that $\cH_*\big(\Phi_S-Id_{S_*}\big)\in H^{s+0.5}(\Om_*,\R^2)$ for $1.5\le s\leq 8.5$ with corresponding estimate
\[
\begin{split}
\|\cH_*(\Phi_S-Id_{S_*}\big)\|_{H^{s+0.5}(\Om_*)}&\le C \big(\|d_{\G}\mu\|_{H^s(\G_{*})}+\|\Phi_S-Id\|_{H^s(\G_{b*})}\big)\\
&\le C\big(\|d_{\G}\|_{H^s(\G_{*})}+\|\Phi_S-Id\|_{H^s(\G_{b*})}\big),
\end{split}
\]
where the constant $C=C(\Om_*,\mu)$ is uniform in $\Lam(S_*,s,\del, \pi/16)$. Moreover, one has
\[
\|\na \cT_S-I\|_{H^{s-0.5}(\Om_*)}=\|\na \cH_*(\Phi_S-Id_{S_*}\big)\|_{H^{s-1/2}(\Om_*)}
\le C\big(\|d_{\G}\|_{H^s(\G_{*})}+\|\Phi_S-Id\|_{H^s(\G_{b*})}\big),
\]
which implies that $\cT_S$ is a diffeomorphism from $\Om_*$ to $\Om=\cT_S(\Om_*)$.

\medskip
As a result, the map $\cT_S$ can be used as coordinates on $\Om_*$.

\bigskip

\section{Equivalent formulation of the problem}

We are going to introduce the new quantity $\mN_a$ and derive an equivalent system of \mbox{(WW)}.
From now on, we consider a family of upper free surfaces $\G_t$ with time variable $t$ in $\Lam(S_*,s,\del, \pi/16)$. The unit outward normal vector is denoted by $N_t$ and the unit tangent vector is $\tau_t$, while the domain is  $\Om_t$.

To begin with, denoting by $\Om_0$ and $\G_0$ the initial domain and upper surface respectively, the velocity filed $v$ induces a  flow map $U(t,\cdot):\ S_0\rightarrow S_t$ by
\[
U(0,\cdot)=Id_{S_0},\quad \pa_t U(t,\cdot)=v(t,U(t,\cdot)),
\]
and the material derivative is
\[
D_t\triangleq \pa_t+\na_v.
\]

\subsection{Some commutators} We recall some commutators involving $D_t$ from \cite{SZ, MW} on the surface $\G_t$  or in the domain $\Om_t$.  

To start with, we recall directly from \cite{SZ, MW} that
\[
D_t\tau_t=(\na_{\tau_t}v\cdot N_t)\,N_t,\quad D_tN_t=-\big((\na v)^*N_t\big)^\top\qquad\hbox{on}\quad \G_t.
\]

Moreover, we have
\beq\label{Dt k expression}
D_t\kappa=-\D_{\G_t}v\cdot N_t-2\Pi(\tau_t)\cdot \na_{\tau_t}v \qquad\hbox{on}\quad \G_t.
\eeq
\medskip

The commutators  are listed here:\\
\noindent 1. $[D_t,\cH]$. Recalling  directly from \cite{MW}, we obtain for a function $f$ on $\G_t$ that
\beq\label{commutator Dt H}
[D_t,\cH]f=\D^{-1}\big(2\na v\cdot\na^2f_{\cH}+\D v\cdot\na f_{\cH},\, (\na_{N_b}v-\na_v N_b)\cdot \na f_{\cH}\big) \qquad\hbox{on}\quad \Om_t,
\eeq
where $\D^{-1}(h,g)$ and $\cH$ are already defined in the Notation part.

\noindent 2. $[D_t, \cN]$.
Again, we quote directly that
\beq\label{commutator DN}
\begin{split}
[D_t,\,\cN]f=&\na_{N_t}\D^{-1}\Big(2\na v\cdot \na^2f_{\cH}+\D v\cdot\na f_{\cH},\,(\na_{N_b}v-\na_{v}N_b)\cdot\na f_{\cH}\Big)\\
&\,-\na_{N_t}v\cdot \na f_{\cH}-\na_{(\na f_{\cH})^\top} v\cdot N_t\qquad\qquad\hbox{on}\quad \G_t.
\end{split}
\eeq

\noindent  3. $[D_t,\D_{\G_t}]$.  For a function $f$ on $\G_t$, one has  
\beq\label{commutator surface delta}
[D_t,\,\D_{\G_t}]f=2\cD^2f\big(\tau_t,(\na_{\tau_t}v)^\top\big)-(\na f)^\top\cdot \D_{\G_t}v+\kappa \na_{(\na f)^\top} v\cdot N_t  \qquad\hbox{on}\quad \G_t.
\eeq

\noindent 4. $[D_t,\,\D^{-1}]$. We have
\beq\label{commutator Dt D-1}
\begin{split}
D_t\D^{-1}(h,g)=&\D^{-1}(D_th,\,D_tg)\\
&\,
+\D^{-1}\Big(2\na v\cdot\na^2\D^{-1}(h,g)+\D v\cdot \na \D^{-1}(h,g),\,
(\na_{N_b}v-\na_{v}N_b)\cdot \na\D^{-1}(h,g)\Big).
\end{split}
\eeq

\subsection{Recover of the velocity}
%The equation for the velocity $v$ of $\Om_t$ as well as the equation of $\ka$  will be derived.

Based on the construction of the free surface, one can  recover the velocity fields $v$. In fact,  we start with the evolution of the boundary $\G_t$  expressed by $d_{\G_t}$.

%Although some of the following computations are already done in \cite{SZ2}, we still provide some details in order to be self-content.

To begin with,  the normal component of $\pa_t \Phi_{S_t}$ represents  the normal component of velocity $v$ on $\G_t$, which means
\[
\pa_t \Phi_{S_t}\cdot (N_t\circ \Phi_{S_t})=(v\cdot N_t)\circ \Phi_{S_t},
\]
where
\[
\pa_t\Phi_{S_t}=\pa_t d_{\G_t}\mu\qquad\hbox{on}\quad \G_*.
\]
Therefore  we  have
\beq\label{normal derivative of v}
\pa_t d_{\G_t}=\f{(v\cdot N_t)\circ \Phi_{S_t}}{\mu\cdot (N_t\circ \Phi_{S_t})}\quad\hbox{i.e.} \quad
v\cdot N_t=(\pa_t d_{\G_t}\mu)\circ \Phi^{-1}_{S_t}\cdot N_t.
\eeq
Since the velocity $v$ in our paper is assumed to be irrotational, we define
\beq\label{v potential expression}
v=\na \phi
\eeq
with $\phi$ satisfying
\beq\label{phi system}
 \left\{\begin{array}{ll}
\Delta \phi=\xi \gamma\qquad\hbox{on}\quad \Om_t,\\
\na_{N_t} \phi|_{\Gamma_t}= (\pa_t d_{\G_t} \,\mu)\circ \Phi^{-1}_{S_t}\cdot N_t, \quad \na_{N_b} \phi|_{\Gamma_b}= 0,
\end{array}\right.
\eeq
where
\[
\gamma=\big(\int_{\Om_t}dX\big)^{-1},\quad\hbox{ and }\quad \xi=\int_{\Gamma_t} v\cdot N_t\,ds.
\]
The definition of $\phi$ is motivated by \cite{SZ2}. Here, we must point out that by the definition of $d_{\G_t}$, we can not get $\dive v=0$ right now. But we prove this in the last section. 
%Moreover, we will show in the end that $\dive v=0$ indeed.

\medskip
On the other hand, the flow map $U(t)$ together with the diffeomorphism $\Phi_{S_t}$ induces a conjugate flow map $U_*:\ S_*\rightarrow S_*$ by
\[
U_*(t,\cdot)\triangleq \Phi^{-1}_{S_t}\circ U(t,\cdot)\circ\Phi_{S_0},
\]
and  the corresponding velocity $v^*$ is given by
\[
\pa_t U_*(t,\cdot)=v^*(t, U_*(t,\cdot)).
\]

We express $v^*$ on the upper surface $\G_*$ in terms of $v$ and $\Phi_{S_t}$. To begin with, we know from the definition of $U_*(t)$ that
\[
\Phi_{S_t}\circ U_*(t)=U(t)\circ \Phi_{S_0}.
\]
Applying $\pa_t$ on both sides and constraining the computation on the surfaces $\G_t$ and $\G_*$ lead to
\[
\pa_t\Phi_{S_t}\circ U_*(t)+\big(D\Phi_{S_t}\circ U_*(t)\big)\,\pa_tU_*(t)=\pa_tU(t)\circ \Phi_{S_0}.
\]
Consequently, one obtains
\beq\label{v* expression}
v^*=D\Phi^{-1}_{S_t}(v\circ \Phi_{S_t}-\pa_t d_{\G_t}\mu)\qquad\hbox{on}\quad \G_*.
\eeq
Recalling from \eqref{normal derivative of v}, one  has
\[
v^{*}= D\Phi_{S_t}^{-1}(v^\top\circ \Phi_{S_t} -\pa_t d_{\G_t} \mu^{\top} ).
\]
Denoting the material derivative related to $U_*(t,\cdot)$ and $v^*$ by
\[
D_{t*}\triangleq \pa_t+\na_{v^*},
\]
a direct computation shows that
\ben\label{change of Dt}
(D_t f)\circ \Phi_{S_t}=D_{t^*}(f\circ \Phi_{S_t}) ,
\een
for any function $f$ on $\G_t$.

\bigskip

Next, we want to consider about the variation of $v^*$ with respect to a parameter $\tau$. Therefore, $\Phi_{S_t}$ and  $\G_t$ depends on $\tau$. Rewriting \eqref{v* expression} by
\[
(D\Phi_{S_t})v^*=v\circ \Phi_{S_t}-\pa_t d_{\G_t}\mu
\]
and taking $\pa_\tau$ on both sides lead to
\[
\pa_\tau(D\Phi_{S_t})v^*+(D\Phi_{S_t})\pa_\tau v^*=\pa_\tau (v\circ \Phi_{S_t})-\pa_\tau\pa_td_{\G_t}\mu,
\]
where
\[
\pa_\tau (v\circ \Phi_{S_t})=(D_\tau v)\circ \Phi_{S_t},\quad
\pa_\tau(D\Phi_{S_t})=D(\pa_\tau d_{\G_t}\mu).
\]
Consequently, we arrive at the expression
\beq\label{pa v*}
\pa_\tau v^*=(D\Phi_{S_t})^{-1}\Big((D_\tau v)\circ \Phi_{S_t}-\pa_\tau\pa_t d_{\G_t}\mu-D(\pa_\tau d_{\G_t}\mu)\,D\Phi^{-1}_{S_t}(v\circ \Phi_{S_t}-\pa_t d_{\G_t}\mu)\Big).
\eeq
For the moment, we still need to rewrite $D_\tau v$ in \eqref{pa v*}. In fact, one has immediately by recalling  from \eqref{v potential expression} that
\[
D_\tau v=D_\tau \na \phi=\na D_\tau \phi+[D_\tau,\,\na]\phi,
\]
with  $[D_\tau,\,\na]=-\na v_\tau\cdot \na \phi$.

Applying $D_\tau$ on the system \eqref{phi system} for $\phi$, one can write
\[
 D_\tau \phi\triangleq u_1+u_2,
 \]
 where $u_1$ satisfies
\[
 \left\{\begin{array}{ll}
\Delta u_1=0\qquad\hbox{on}\quad \Om_t,\\
\pa_{N_t} u_1|_{\Gamma_t}= (\pa_\tau\pa_t d_{\G_t} \,\mu)\circ \Phi^{-1}_{S_t}\cdot N_t, \quad \pa_{N_b} u_1|_{\Gamma_b}= 0,
\end{array}\right.
\]
and $u_2$ the remainder part satisfies
\[
 \left\{\begin{array}{ll}
\Delta u_2=2\na_{v_\tau}\cdot \na^2\phi+\D v_\tau\cdot\na\phi+D_\tau(\gamma\xi) \qquad\hbox{on}\quad \Om_t,\\
\pa_{N_t} u_2|_{\Gamma_t}= (\pa_t d_{\G_t}\mu)\circ \Phi^{-1}_{S_t}\cdot D_\tau N_t-\big(D_\tau N_t-\na_{N_t}v_\tau\big)\cdot \na\phi, \\
\pa_{N_b} u_2|_{\Gamma_b}= \big(\na_{N_b}v_\tau-\na_{v_\tau}N_b\big)\cdot \na \phi\big|_{\G_b}
\end{array}\right.
\]
with $D_\tau N_t=-\big((\na v_\tau)^* N_t\big)^\top$.

As a result, one arrives at
\beq\label{D be v}
D_\tau v=\na u_1+\na u_2-\na v_\tau\cdot \na \phi.
\eeq
Substituting this expression together with \eqref{v* expression} back into \eqref{pa v*}, one finally derives
\beq\label{pa v* expression}
\pa_\tau v^*=(D\Phi_{S_t})^{-1}\Big((\na u_1+\na u_2-\na v_\tau\cdot \na \phi)\circ \Phi_{S_t}-\pa_\tau\pa_t d_{\G_t}\mu-D(\pa_\tau d_{\G_t}\mu)\,v^*\Big).
\eeq

\bthm{Lemma}
Let the surface $\G_t\in \Lambda(S_*, s, \delta, {\pi}/{16})$. Then we have
\[
\begin{split}
&\|v^*\|_{H^s(\G_*)} \le C\big(\|d_{\G_t}\|_{H^{s+1}(\G_*)}, \|\pa_t d_{\G_t}\|_{H^s(\G_*)}, \|v\|_{H^{s}(\G_t)}\big),\\
&\|\pa_\tau v^*\|_{H^s(\G_*)}\le C\big(\|d_{\G_t}\|_{H^{s+1}(\G_*)}, \|\pa_t d_{\G_t}\|_{H^s(\G_*)}, \|v\|_{H^{s}(\G_t)}\big)\Big(\|\pa_\tau d_{\G_t}\|_{H^{s+1}(\G_*)}+\|\pa_t \pa_\tau d_{\G_t}\|_{H^s(\G_*)}+\|D_\tau v\|_{H^{s}(\G_t)}\Big).
\end{split}
\]
\ethm
\begin{proof}
Applying Lemma \ref{trace thm PG} and Lemma \ref{est:elliptic} directly,  the proof can be finished.
\end{proof}

\medskip

\subsection{The new quantity $\mN_a$}
We introduce  the new quantity $\mN_a$ related to the mean curvature $\ka$ for some large constant $a>0$:
\beq\label{def:d}
\mN_a=\cN(\ka)\circ\Phi_S+a^3 d_{\G_t},
\eeq
which is used to rewrite the water-waves problem $\mbox{(WW)}$.
The constant $a$ is taken large enough, which is  useful in the following proposition and also in the iteration section. 

Before we derive the equation for $\mN_a$, we need to make sure that we can recover the upper surface $\G_t$ from it.   In order to do this, we need also the  boundary condition at the two end points $p_{l*}, p_{r*}$:
\beq\label{point condition for d}
d_{\G_t}(p_{i*})=d_i,\quad i=l, r,
\eeq
where $d_i=d_i(t)$ is some function to be given later. 

We define the operator $\cK$  by
%from equation \eqref{def:d} and the boundary conditions point condition for $d_i$ above:
\[
\cK: \, d_{\G_t}\rightarrow (\mN_a, d_l, d_r).
\]
When $a$ is chosen to be big enough, we prove that the operator $\cK$ is invertible. 
%We need to prove that the operator $\cK$ is invertible.

%One can prove  immediately  that for any given $(\mN_a, d_l, d_r)$, there exists a unique $d_{\G_t}$ with corresponding  estimate.

\bthm{Proposition}\label{d and  N ka}
{\it
Let  $s \in[1, 5.5]$ and the set $B_{\delta_1}$ is defined by
\beno
B_{\delta_1}\triangleq \{ (\mN_a,d_l, d_r)\,|\, \| \mN_a-  \mN_{a*}\|_{H^{s}(\Gamma_*)}, |d_l|, |d_r|<\delta_1\},
\eeno
where $\mN_{a*}$ is the value of $\mN_a$ taken $\Gamma_*$. Then, there exists $A>0$ such that when $a\ge A>0$,  the operator $\cK$ is a diffeomorphism from $\Lam(S_*, s+3, \del, \pi/16)\subset H^{s+3}(\G_*)$ to $\cK\big(\Lam(S_*, s+3, \del, \pi/16)\big)\subset B_{\delta_1}$ and there  holds
\[
\|d_{\G_t}\|_{H^{s+3}(\G_*)}\le C\big(\|\mN_a-\mN_{a*}\|_{H^{s}(\G_*)}+|d_l|+|d_r|\big),
\]
where the constant $C$ depends on $\delta, \delta_1, a, \G_*$. Moreover, one also has the following estimate involving the constant $a^{-1}$:
\[
\|d_{\G_t}\|_{H^{s_1}(\G_*)}\le a^{s_1-s_2-3} C\big(\|\mN_a-\mN_{a*}\|_{H^{s_2}(\G_*)}+|d_l|+|d_r|\big)
\]
with $8.5 \ge s_1\ge s_2\ge 1$ and $s_1\le s_2+3$.
}
\ethm

\begin{proof}
%The idea for the proof is similar as  in Lemma 2.2 \cite{SZ2}, but there is also a big difference in our case due to the boundary points.

Since the problem is quasilinear, we first linearize  the operator $\cK$ around $\G_*$ (where $d_{\G_*}=0$) as in (2.4) \cite{SZ2} to obtain the linearized operator $\mathcal L(\G_*)$:
\[
\mathcal L(\G_*)d \triangleq \Big(\f{\del\mN_a}{\del d_{\G}}(\G_*)d, \ d(p_{l*}), d(p_{r*})\big),
\]
where the linearization for $\mN_a$ is derived by a variation  with parameter $\tau$:
\beq\label{variation of N ka}
\begin{split}
\f{\del \mN_a}{\del d_{\G_t}}(\G_*)d=&a^3 d+D_\tau\cN(\ka)\big|_{\tau=0}\\
=&a^3 d+\cN_*\big(D_\tau\ka\big|_{\tau=0}\big)
-\na_{N_*}(\mu\,d)_{ex}\cdot\na\cH(\ka_*)-\na_{(\na \cH(\ka_*))^\top}(\mu\,d)_{ex}\cdot N_*\\
&\quad +\na_{N_*}\D^{-1}\big(2\na(\mu\,d)_{ex}\cdot\na^2\cH(\ka_*)+\D(\mu \,d)_{ex}\cdot\na \cH(\ka_*),\, (\na_{N_b}(\mu\,d)_{ex}-\na_{(\mu\,d)_{ex}}N_b)\cdot\na\cH(\ka_*)\big),
\end{split}
\eeq
where \eqref{commutator DN} is applied and $\N_*, \ka_*, \cN_*$ are defined on $\G_*$.  Besides, $v|_{\tau=0}=(\mu \,d)_{ex}$ denotes the extension of $\mu\, d$ from $\G_t$ to $\Om_t$ (for example, use the trace theorem in \cite{MW1}). The variation is defined with the velocity filed $v(\tau): S_\tau\rightarrow \R^2$ where $S_\tau$ is a family of surfaces.

Moreover, one has
\[
D_\tau \ka\big|_{\tau=0}=\f{\del (\ka\circ\Phi_{S})}{\del d_{\G}}(\G_*)d,
\]
where we quote the computation from \cite{SZ2} directly:
\[
\f{\del(\ka\circ \Phi_{S_t})}{\del d_{\G}}(\G_*)d=-(\mu\cdot N_*)\D_{\G_*}d-2(\na_{N_*}\mu\cdot\tau_*)\na_{\tau_*}d -(N_*\cdot\D_{\G_*}\mu)\,d-2(\Pi_*\cdot\tau_*)\tau_*\cdot\na_{\tau_*}(d\mu).
\]
Summing up these expression above, we can conclude that
\[
\f{\del \mN_a}{\del d_{\G_t}}(\G_*)d=a^3 d-\cN_*\big((\mu\cdot N_*)\D_{\G_*}d\big)+r_d
\]
with the remainder term
\[
\begin{split}
r_d=r_d(\pa^2 d)&=-\cN_*\Big(2(\na_{N_*}\mu\cdot\tau_*)\na_{\tau_*}d +(N_*\cdot\D_{\G_*}\mu)\,d+2(\Pi_*\cdot\tau_*)\tau_*\cdot\na_{\tau_*}(d\mu)\Big)\\
&\quad -\na_{N_*}(\mu\,d)_{ex}\cdot\na\cH(\ka_*)-\na_{(\na \cH(\ka_*))^\top}(\mu\,d)_{ex}\cdot N_*\\
&\quad +\na_{N_*}\D^{-1}\big(2\na(\mu\,d)_{ex}\cdot\na^2\cH(\ka_*)+\D(\mu \,d)_{ex}\cdot\na \cH(\ka_*),\, (\na_{N_b}(\mu\,d)_{ex}-\na_{(\mu\,d)_{ex}}N_b)\cdot\na\cH(\ka_*)\big).
\end{split}
\]
As a result,  if one can prove that this linear operator $\mathcal  L(\G_*)$ is invertible and continuous, one immediately has that $\cK$ is a diffeomorphism near $\G_*$ in $\Lam(S_*, s, \del, \pi/16)$ by the inverse function theorem (see for example Theorem 1.2.3\cite{Chang})  and the desired estimate will follow.

So all we need is to prove now is the invertibility of $\mathcal  L(\G_*)$, which  can be done by applying standard elliptic analysis on the following system
\[
\mathcal  L(\G_*)d =(h, d_l, d_r),
\]
for some $h\in H^{s}(\G_*)$ and $ d_l, d_r\in\R$.

 Firstly,  we  consider about the variation solution in the space $\mathcal V=\{d\in H^{2.5}(\G_*)\,\big|\, d(p_{i*})=0, i=l, r\}$.  Here  a standard analysis is used to obtain zero Dirichlet boundary condition at $p_i$ for $\tilde d=d-f$, where
\[
 f(p_{i*})=d_i\ (i=l, r),\qquad\hbox{and}\quad \|f\|_{ H^s(\G_*)}\le C(|d_l|+|d_r|).
\]

In order to apply Lax-Milgram's Theorem, we define
\[
A(d_1, d_2)= \int_{\G_*}\Big(a^3 d_1-\cN_*\big((\mu\cdot N_*)\D_{\G_*}d_1\big)+r_{d_1}\Big) \big(1-a^{-1}(\mu\cdot N_*)\D_{\G_*}\big)d_2\,ds,
\]
and
\[
F(d_2)=\int_{\G_*}h\big(1-a^{-1}(\mu\cdot N_*)\D_{\G_*}\big)d_2\,ds.
\]

One easily shows that the conditions in Lax-Milgram's Theorem are all satisfied when $a$ is large enough, so the linear system admits a unique solution $d\in \mathcal V$ and the estimate in $\mathcal V$ follows. 
Here we remark that the boundary conditions $d_i$ are used for the operator $\Delta_{\G_*}$, while no boundary conditions are needed due to the symmetry property of $\cN$.

Secondly, we  prove the higher-order estimates. For example, setting
\[
u_1=\cN_*\big((\mu\cdot N_*)\D_{\G_*}d)-\cN_*\Big(2\big((\na_{N_*}\mu\cdot \tau_*)+(\Pi_*\cdot\tau_*)\tau_*\cdot\mu\big)\na_{\tau_*}d\Big)
\]
and taking $\D_{\G_*}$ on the equation of $d$, we derive the linear system for $u_1$ with Dirichlet boundary conditions. Therefore, the higher-order estimates are proved.

Moreover, noticing that there is $a^3$ in $\mN_a$ and the linearized equation as well, so one can obtain elliptic estimates with weight $a$ for the linearized equation using interpolations. Consequently, one can prove the desired weighted estimate and the proof is finished.
\end{proof}

\bthm{Remark}\label{inverse Aa}{\it Following the proof above, we know immediately that the elliptic operator 
\[
\cA_a\triangleq a^3+\cA (d_{\G_t})
\]
 is invertible with Dirichlet boundary conditions at $p_{i*}$, and corresponding elliptic estimates follow naturally.
}
\ethm

As a result, as long as we have $d_{\G_t}$, we   retrieve $\Phi_{S_t}$  to define the surface $\G_t$ and the corresponding domain $\Om_t$.

\bigskip

\subsection{Evolution of $\mN_a$ and the  boundary conditions at $p_{i}$}
Now, we are in a position to derive the evolution equation of $\mN_a$ on $\G_*$ from the water-waves problem $\mbox{(WW)}$. The boundary conditions at the corner points are adjusted to the version for the new quantity $\mN_a$. In the end, the evolution equation for $d_i=d_{\G_t}(p_{i*})$ $(i=l, r)$ are derived.

Firstly, since $\mN_a$ is related to $\ka$ defined on $\G_t$, we begin with the equation of $\ka$. In fact, it has been proved and used in \cite{SZ2, MW} that the mean curvature $\kappa$ satisfies the equation
\beq\label{ka equation}
D^2_t\ka=-N_t\cdot \D_{\G_t}D_t v+2\si \Pi(\tau_t)\cdot  \na_{\tau_t}\na\ka_{\cH}+r_1,
\eeq
where $r_1$ is the lower order term
\[
\begin{split}
r_1=& 2\Big[\cD\big((\na v)^*N_t\big)^\top+\Pi\big((\na_{\tau_t} v)^\top\big)\Big]\cdot \na_{\tau_t} v+\D_{\G_t}v\cdot\big((\na v)^*N_t\big)^\top+2\Pi(\tau_t)\cdot (\na v)^2\tau_t\\
&
\  -2(\na_{\tau_t} v\cdot N_t)\Big(\Pi(\tau_t)\cdot\na_{N_t}v
+\Pi(\tau_t)\cdot(\na v)^*N_t\Big)
-2N_t\cdot D^2v\big(\tau_t,\,(\na_{\tau_t}v)^\top\big)\\
&
\ -2(\na_{N_t}v\cdot N_t)(\Pi(\tau_t)\cdot \na_{\tau_t} v)+N_t\cdot \na v\big((\D_{\G_t}v)^\top\big) -\kappa\big|\big((\na v)^*N_t\big)^\top\big|^2+2\Pi(\tau_t)\cdot \na_{\tau_t}\na P_{v,v}.
\end{split}
\]
One can tell that the highest-order terms in $r_1$ are like $\pa^2 v,\,\pa N_t$.

Substituting  Euler equation
\[
D_t v=-\si \na \ka_{\cH}-\na P_{v,v}+\mathbf{g},
\]
into the equation above, a direct computation leads to
\beq\label{ka equation}
D_t^2 \kappa=\sigma \Delta_{\Gamma_t} \cN( \kappa)+\widetilde{R}_1,
\eeq
where
\beno
\widetilde{R}_1&=&\si [N_t,\,\D_{\G_t}]\cdot \na \ka_{\cH}+N_t\cdot \Delta_{\Gamma_t}\na P_{v, v}+2\si \Pi(\tau_t)\cdot  \na_{\tau_t}\na\ka_{\cH}+r_1\\
&=& -\si\D_tN_t \cdot \na \ka_{\cH}+N_t\cdot \Delta_{\Gamma_t}\na P_{v, v}+r_1,
\eeno
and one can tell that the highest-order terms in $\widetilde{R}_1$ are like $\pa\ka,\,\pa^2 v$.

Secondly,  acting $\cN$ on both sides of equation \eqref{ka equation} leads to
\[
D_{t}^2 \cN(\ka)=\sigma \cN \Delta_{\Gamma_t} \cN(\ka) +\bar R_1\circ \Phi^{-1}_{S_t}\qquad\hbox{on}\quad \G_t,
\]
where $\bar R_1$ is defined by
\beno
\bar R_1=\Big(D_t[D_{t}, \cN]\ka+[D_{t}, \cN]D_t\ka+\cN\widetilde{R}_1\Big)\circ\Phi_{S_t}
\eeno
with the commutator expressed by \eqref{commutator DN} and $D_t v$ replaced by $-\si\na\ka_\cH-\na P_{v,v}+{\bf g}$ from Euler equation.

Moreover, substituting $\cN(\ka)=\mN_a\circ\Phi^{-1}_{S_t}-a^3d_{\G_t}\circ \Phi_{S_t}^{-1}$ into both sides of the equation above,
 we obtain  immediately that
\[
D^2_t\big(\mN_a\circ\Phi^{-1}_{S_t}\big)=\sigma \cN \Delta_{\Gamma_t} \big(\mN_a\circ\Phi^{-1}_{S_t}\big)
+R_0\circ \Phi_{S_t}^{-1}
\]
with the remainder term
\beq\label{definition of R_0}
R_0=R_a +\bar{R}_1,
\eeq
and
\[
R_a=\Big(a^3D_{t}^2 (d_{\G_t}\circ \Phi_{S_t}^{-1})-a^3\sigma \cN\Delta_{\Gamma_t}\big( d_{\G_t}\circ \Phi_{S_t}^{-1}\big)\Big)\circ \Phi_{S_t}.
\]

To remove the second-order time derivative of $d_{\Gamma_t}\circ \Phi_{S_t}^{-1}$ in $R_a$, we rewrite
\[
D_{t}^2 (d_{\G_t}\circ \Phi_{S_t}^{-1} )=(D^2_{t*}d_{\G_t})\circ  \Phi_{S_t}^{-1},
\]
where
\beq\label{D 2 t d expression}
D^2_{t*}d_{\G_t}=\pa_t^2 d_{\G_t}+\na_{\pa_tv^*}d_{\G_t}+2\na_{v^*}\pa_td_{\G_t}+\na_{v^*}\na_{v^*}d_{\G_t},
\eeq
and one has from \eqref{v* expression}  that
\beq\label{pa t v*}
\pa_t v^*=(D\Phi_{S_t})^{-1}\big[(D_t v)\circ \Phi_{S_t}-\pa^2_t d_{\G_t}\mu+\big((\pa_t d_{\G_t}\mu\circ \Phi^{-1}_{S_t}-v)\cdot \na v\big)\circ \Phi_{S_t}\big]+\pa_t\big((D\Phi_{S_t})^{-1}\big)(v\circ \Phi_{S_t}-\pa_t d_{\G_t} \mu).
\eeq

Moreover, we need to express $\pa_t^2d_{\G_t}$ in terms of $\na v, \pa^2d_{\G_t}, \pa_t d_{\G_t}$. In fact, recalling from \eqref{normal derivative of v} and \eqref{change of Dt}, one obtains immediately
\[
D_t (v\cdot N_t)=D_t\big((\pa_t d_{\G_t}\mu)\circ\Phi^{-1}_{S_t}\cdot N_t\big),
\]
or equivalently
\beq\label{Dt  v and d eqn}
D_t v\cdot N_t=(\mu\circ\Phi^{-1}_{S_t}\cdot N_t)(\pa^2_t d_{\G_t}+\na_{v^*}\pa_t d_{\G_t})\circ\Phi^{-1}_{S_t}+(\na_{v^*}\mu\circ\Phi^{-1}_{S_t}\cdot N_t)\pa_t d_{\G_t}\circ\Phi^{-1}_{S_t}
+\big((\pa_t d_{\G_t}\mu)\circ\Phi^{-1}_{S_t}-v\big)\cdot D_t N_t.
\eeq
Substituting  Euler equation from $\mbox{(WW)}$ into \eqref{Dt  v and d eqn} and expressing $\cN(\ka)$ with $\mN_a$, one arrives at
\beq\label{pa 2 t d}
\begin{split}
\pa_t^2 d_{\G_t}=-\f{1}{\mu\cdot (N_t \circ\Phi_{S_t})} &\Big(\mu\cdot (N_t \circ\Phi_{S_t})\na_{v^*}\pa_td_{\G_t}+\na_{v^*}\mu\cdot (N_t\circ\Phi_{S_t})\pa_t d_{\G_t}+(\pa_td_{\G_t}\mu-v\circ\Phi_{S_t})\cdot (D_tN_t\circ\Phi_{S_t})\\
&\quad +\big(\si\mN_a-\si a^3 d_{\G_t}+\na P_{v,v}\cdot N_t+{\bf g}\cdot N_t\big)\circ \Phi_{S_t}\Big)
\end{split}
\eeq

As a result, $R_a$ is expressed by
\beq\label{Ra expression}
R_a=a^3 \big(\pa_t^2 d_{\G_t}+\na_{\pa_tv^*}d_{\G_t}+2\na_{v^*}\pa_td_{\G_t}+\na_{v^*}\na_{v^*}d_{\G_t}\big)-a^3\sigma \Delta_{\Gamma_t} \cN(d_{\G_t}\circ \Phi_{S_t}^{-1})\circ \Phi_{S_t},
\eeq
where $\pa^2_t d_{\G_t}$, $\pa_tv^*$ are given by \eqref{pa 2 t d}, \eqref{pa t v*}. One knows immediately that $R_a=R_a(a^3\pa^3 d_{\G_t}, a^3\pa\pa_td_{\G_t}, a^3\pa v)$.

\bigskip
Now summing up all these computations above, we get the following lemma.
\bthm{Lemma}\label{lemma:equ mN_a}
{\it We have that $\mN_a$ satisfies the following equation
\beq\label{ka_a equation}
D_{t*}^2 \mN_a+\sigma \cA (d_{\G_t}) \mN_a=R_0\qquad\hbox{on}\quad \G_*,
\eeq
where
\[
R_0=R_0\big(\pa^4 d_{\G_t}, \pa^3\pa_t d_{\G_t}, \pa^2 v, \pa\mN_a, a^3\pa^3 d_{\G_t}, a^3\pa\pa_t d_{\G_t}, a^3\pa v\big)
\]
 is defined in \eqref{definition of R_0} and the operator $\cA(d_{\G_t})$ is defined by
 \beno
  \cA (d_{\G_t}) f\triangleq-\big( \cN \Delta_{\Gamma_t} (f\circ\Phi^{-1}_{S_t})\big)\circ \Phi_{S_t}\qquad\hbox{for some $f$ on $\G_*$}.
 \eeno

}
\ethm

The boundary conditions at the contact points in \mbox{(WW)} were rewritten in  Lemma 7.1 \cite{MW},  which gives the boundary condition at the left contact point $p_l$ for $J=\na\ka_\cH$. By a similar argument,  we find the corresponding boundary conditions on $p_{i*}$ for $\mN_a$.
\bthm{Lemma}\label{lem:cc}
{\it 
We have the following conditions for $\mN_a$ at the contact points $p_{i*}$ $(i=l, r)$ on $ \G_*$
 \beq\label{equ:cca}
\begin{split}
& D_{t*}\cA(d_{\G_t})\mN_a +\f{\sigma^2 }{\beta_c}(\sin\om_i)^2\circ\Phi_{S_t}\big( \na_{\tau_t}( \cA(d_{\G_t})\mN_a\circ\Phi^{-1}_{S_t})\big)\circ\Phi_{S_t}
=R_{c, i}\qquad \textrm{at}\quad p_{l*},\\
& D_{t*}\cA(d_{\G_t})\mN_a -\f{\sigma^2 }{\beta_c}(\sin\om_i)^2\circ\Phi_{S_t}\big( \na_{\tau_t}( \cA(d_{\G_t})\mN_a\circ\Phi^{-1}_{S_t})\big)\circ\Phi_{S_t}
=R_{c, i}\qquad \textrm{at}\quad p_{r*},
\end{split}
\eeq
where the remainder term $R_{c, i}=R_{c, i}(\pa^5d_{\G_t}, \pa^4\pa_t d_{\G_t}, \pa^3 v, \pa^2\mN_a, \pa\pa_t\mN_a)$ is expressed in details  as below:
\[
\begin{split}
R_{c, i}=&R_{ca}-\f1\si D_{t*}R_0-\f{\si}{\be_c}(\sin\om_i)^2\big(\na_{\tau_t}(R_0\circ\Phi^{-1}_{S_t})\big)\circ \Phi_{S_t}+(D^2_{t}r_c)\circ\Phi_{S_t}\\
&\quad
 -\f{\si}{\be_c}\Big[\big(D_{t*}(\sin\om_i)^2\circ\Phi_{S_t}\big)\big(\na_{\tau_t}(D_{t*}\mN_a\circ\Phi^{-1}_{S_t})\big) +(\sin\om_i)^2\Big((D_t\tau_t-\na_{\tau_t}v)\cdot\na\big( D_{t*}\mN_a\circ\Phi^{-1}_{S_t}\big)\Big)\Big]\circ \Phi_{S_t}\\
&\quad
 -\f{\sigma^2 }{\beta_c}D_t\Big[\big(D_t(\sin\om_i)^2\big)\na_{\tau_t}(\mN_a\circ\Phi^{-1}_{S_t})
-\f{\sigma^2 }{\beta_c}(\sin\om_i)^2\big(D_t\tau_t-\na_{\tau_t}v\big)\cdot\na(\mN_a\circ\Phi^{-1}_{S_t})\Big]\circ \Phi_{S_t},
\end{split}
\]
with
\[
R_{ca}=a^3 D^3_{t*}d_{\Gamma_t}+a^3 D^2_{t*}\Big(\f{\sigma^2 }{\beta_c}(\sin\om_i)^2 \big(\na_{\tau_t}(d_{\Gamma_t}\circ \Phi^{-1}_{S_t} )\big)\circ\Phi_{S_t}\Big)\\
=R_{ca}(a^3\pa^4d_{\G_t}, a^3\pa^2\pa_t d_{\G_t}, a^3\pa^2v)
\]
and
\[
\begin{split}
r_c=&\Big(\sigma  \sin\om_i  \big(\na_{\tau_t}  \na P_{v,v}   -\na_{N_t} v\na_{\tau_t}v\cdot N_t+\na_{\tau_t}v\cdot\na v\big)\cdot N_t+\sigma  \sin\om_i \big(\na_{\tau_t}v\cdot ((\na v)^*N_t)^\top\big)\Big) \,\tau_b\cdot N_t
\\
&\quad -\sigma  \tau_b\cdot ((\na v)^*N_t)^\top\,(\na_{\tau_t}v\cdot N_t)(\tau_b\cdot N_t)
-((\na v)^*N_t)^\top \cdot \na \ka_\cH+\f{\sigma^2 }{\beta_c}(\sin\om_i)^2 \na_{\tau_t} N_t\cdot  \na \ka_\cH\\
&\quad  +\beta_c (\tau_b\cdot N_t)\na_{\tau_b} v\cdot\na P_{v,v}
-\beta_c (\tau_b\cdot N_t)\na_ {\tau_b}\D^{-1}\big(-2tr(\na D_tv\na v)+2tr(\na v\cdot \na v\na v), \,D_t(v\cdot \na_v N_b)\big)\\
&\quad -\beta_c (\tau_b\cdot N_t)\na_ {\tau_b}\D^{-1}\big(2\na v\cdot \na^2P_{v,v} +\D v\cdot \na P_{v,v},
(\na_{N_b}v-\na_vN_b)\cdot\na P_{v,v}\big).
\end{split}
\]
Here $D^2_{t*}d_{\G_t}, \pa_t v^*$ are expressed by \eqref{D 2 t d expression}, \eqref{pa t v*} and \eqref{pa 2 t d}, and $D_tv$ is replaced by $-\si \na\ka_\cH-\na P_{v,v}+{\bf g}$.

}
 \ethm

\begin{proof} The boundary conditions will be proved for the left corner point $p_l$, and the case for the right one is similar. To begin with, we know from Lemma 7.1 in \cite{MW} that the conditions at the corner points can be written under the form of $J=\na \ka_\cH$:
\[
D_t J=\f{\si^2}{\be_c}\sin \om_l(\na_{\tau_t}J\cdot N_t)\tau_b+r_J\qquad\hbox{at}\quad p_l,
\]
where
\[
\begin{split}
r_J=&\Big(\sigma  \sin\om_l \big(\na_{\tau_t}  \na P_{v,v}   -[D_t, \na_{\tau_t} ] v\big)\cdot N_t-\sigma  \sin\om_l(\na_{\tau_t}v\cdot D_t N_t)+\sigma  \tau_b\cdot D_t N_t\,(\na_{\tau_t}v\cdot N_t)
\\
&\quad- \beta_c D_t \na P_{v,v}\cdot \tau_b\Big) \,\tau_b.
\end{split}
\]
Applying the inner product with $N_t$ on both sides of the equation above  and noticing that 
\[
\tau_b\cdot N_t=-\sin\om_l,
\]
 one obtains
\[
D_t\cN(\ka)=-\f{\si^2}{\be_c}(\sin \om_l)^2\na_{\tau_t}\cN(\ka)+r_c,
\]
where
\[
r_c=D_tN_t\cdot\na\ka_\cH+\f{\si^2}{\be_c}(\sin \om_l)^2\na_{\tau_t}N_t\cdot\na\ka_\cH+r_J\cdot N_t.
\]
Besides, substituting $\cN(\ka)=\mN_a\circ\Phi^{-1}_{S_t}-a^3 d_{\G_t}\circ\Phi^{-1}_{S_t}$ into the equation above and noticing that
\[
(\na_{\tau_t}f)\circ\Phi_{S_t}=\na_{(\na\Phi_{S_t})^*\tau_t\circ\Phi_{S_t}}(f\circ\Phi_{S_t}),
\]
one can change the equation for $\cN(\ka)$ into an equation for $\mN_a$:
\[
 D_{t*}\mN_a = -\f{\sigma^2 }{\beta_c}(\sin\om_l)^2\big( \na_{\tau_t}( \mN_a\circ\Phi^{-1}_{S_t})\big)\circ\Phi_{S_t}
 +R_{c0}\qquad \textrm{at}\qquad p_{l*},
\]
where
\[
R_{c0}=r_c\circ\Phi_{S_t}+\Big(a^3 D_{t*}d_{\Gamma_t}+a^3 \f{\sigma^2 }{\beta_c}(\sin\om_l)^2 \big(\na_{\tau_t}(d_{\Gamma_t}\circ \Phi^{-1}_{S_t} )\big)\circ\Phi_{S_t}\Big)
\]
with $r_c$ expressed in the statement of this lemma.
Moreover,  one uses Euler equation from $\mbox{(WW)}$ to replace $D_tv$ by $-\si\na\ka_\cH-\na P_{v,v}+{\bf g}$ and applying \eqref{commutator Dt D-1}  and the expressions for $D_tN_t$, $D_t\tau_t$ in the beginning of Section 4.

Secondly, applying $D_{t*}$ twice on the equation above, we get
\[
D^3_{t*}\mN_a=-D^2_{t*}\Big(\f{\sigma^2 }{\beta_c}(\sin\om_l)^2\big( \na_{\tau_t}( \mN_a\circ\Phi^{-1}_{S_t})\big)\circ\Phi_{S_t}\Big)+D^2_{t*}R_{c0}.
\]
Checking term by term and applying \eqref{ka_a equation}, we finally derive the desired condition.
Again, one needs to notice that the remainder terms involving $D_t v$ should be replaced by $-\si\na\ka_\cH-\na P_{v,v}+{\bf g}$ from Euler equation, and $D^2_{t*}d_{\G_t}, \pa^2_td_{\G_t},  \pa_t v^*$ are expressed by \eqref{D 2 t d expression}, \eqref{pa t v*} and \eqref{pa 2 t d}.

\end{proof}

\bthm{Remark}{\it In the following text,  the choices of signs at the left and the right corner points $p_{i*}$ are the same as in \eqref{equ:cca}. 
}\ethm

\bthm{Remark}\label{rmk: com cond}{\it 
By the compatibility conditions \eqref{eq:com cond}, it is easy to get
\beno
 D_{t*}\cA(d_{\G_t})\mN_a|_{t=0}\pm\bigg(\f{\sigma^2 }{\beta_c}(\sin\om_i)^2\circ\Phi_{S_t}\big( \na_{\tau_t}( \cA(d_{\G_t})\mN_a\circ\Phi^{-1}_{S_t})\big)\circ\Phi_{S_t}\bigg)(0)
=R_{c, i}|_{t=0}\quad \textrm{at}\quad p_{i*}\ ( i=l, r).
\eeno
}
\ethm

In the end, except for the equation of $\mN_a$ and the boundary conditions, we also need the evolution equation for $d_{\G_t}$ at the corner points in the iteration.

Notice that at the corner points $p_l, p_r$, the velocity $v$ is tangential along the bottom $\G_b$ and the unit vector $\mu$ is  defined to be tangential along $\G_{b*}$ as well. Combining \eqref{normal derivative of v},  one has
\[
v=(\pa_t d_{\G_t}\mu)\circ\Phi^{-1}_{S_t}\quad\hbox{at}\quad p_i\ (i=l, r).
\]
Substituting this equality into \eqref{pa 2 t d}, one derives the evolution equation for $d_i(t)=d_{\G_t}(p_{i*})$ $ (i=l, r)$:
\beq\label{d c eqn}
 d_i''(t)=\mathfrak{B}_i,\quad i=l, r,
\eeq
where
\[
\begin{split}
\mathfrak{B}_i\triangleq-\f{1}{\mu\cdot (N_t \circ\Phi_{S_t})} &\Big(\mu\cdot (N_t \circ\Phi_{S_t})\na_{v^*}\pa_td_{\G_t}+\na_{v^*}\mu\cdot (N_t\circ\Phi_{S_t})\pa_t d_{\G_t}+\si\mN_a-\si a^3d_{\G_t}\\
&\quad +\na P_{v,v}\cdot N_t+{\bf g}\cdot N_t\Big)\big|_{p_{i*}}.
\end{split}
\]

Consequently, we specify the boundary conditions for $d_{\Gamma_t}$ in \eqref{point condition for d}, if the right side term $\mathfrak{B}_i$ and $d_i(0), d'_i(0)$ are known in the iteration section:
\ben\label{equ:boundary condition}
d_i=d_i(0)+d_i'(0)t+\int_0^t\int_0^s \mathfrak{B}_i(\tau)d\tau\,ds\quad \textrm{at} \quad p_{i*} \ (i=l, r).
\een

\bthm{Remark}\label{V0}
{\it According to \eqref{Dt  v and d eqn} and the equation above it, the boundary conditions \eqref{equ:boundary condition}  are in fact corresponding to $D_tv\cdot N_t=-\sigma \cN(\kappa)-\na_{N_t}P_{v, v}+{\bf g}\cdot N_t$ at $p_i$ where $v$ is give by \eqref{v potential expression}.  One will find that these boundary conditions  play an important role in the last section, where we go back to the water-waves system from the iteration.
}
\ethm

\medskip

\subsection{Estimates for the remainder terms $R_0$ and $R_{c, i}$}
We consider the estimates for the remainder terms here. In the iteration scheme, the estimates of the remainder terms $R_0$, $R_{c, i}$ depend on the norms of $d_{\G_t}$, $\pa_td_{\G_t}$. Meanwhile, notice that one can only use the velocity $v=\na\phi$ defined in \eqref{v potential expression}  and relate $N_t, \tau_t, \ka$ with $d_{\G_t}$ as well.  Therefore, $D_t v$, $D_t N_t$, $D_t\tau_t$ and $D_t\ka$  should be estimated in terms of $d_{\G_t}$, $\pa_t d_{\G_t}$ in order to go back to $\mN_a,\pa_t\mN_a, d_i, \pa_t d_i$ by using the operator $\cK$.

To begin with,  we consider about the estimates for $v$, $N_t$, and $\tau_t$.
\bthm{Lemma}\label{Dt v estimate}{\it Let $v=\na\phi$ be defined in \eqref{v potential expression} and $N_t$, $\tau_t$ are the unit normal and tangential vectors defined on $\G_t\in \Lam(S_*, 8.5, \del, \pi/16)$. \\
\noindent (1) The following estimates hold for $N_t$, $\tau_t$, $\ka$ when $s\ge 3$:
\[
\begin{split}
\|N_t\|_{H^s(\G_t)}, \|\tau_t\|_{H^s(\G_t)},  \|\kappa\|_{H^{s-1}(\G_t)}&\le C(\|d_{\G_t}\|_{H^{s+1}(\G_*)}),\\
\|D_t N_t\|_{H^{s-1.5}(\G_t)},\,\|D_t\tau_t\|_{H^{s-1.5}(\G_t)}, \,\|D_t\kappa\|_{H^{s-2.5}(\G_t)}&\le C\big(\|d_{\G_t}\|_{H^{s+0.5}(\G_*)}, \|\pa_td_{\G_t}\|_{H^{s-0.5}(\G_*)}, \|v\|_{H^{s-1}(\Om_t)}\big),
\end{split}
\]
and
\[
\|D^2_t N_t\|_{H^{s-3}(\G_t)},\,\|D^2_t\tau_t\|_{H^{s-3}(\G_t)}\le C\big(\|d_{\G_t}\|_{H^{s}(\G_*)}, \|\pa_td_{\G_t}\|_{H^{s-1}(\G_*)}, \|\pa^2_td_{\G_t}\|_{H^{s-2}(\G_*)}\big);
\]

\noindent (2) For some parameter $\tau$,  the following estimates hold when $s\ge 3$:
\[
\begin{split}
&\|v\|_{H^{s}(\Om_t)}\le C\big(\|d_{\G_t}\|_{H^{s+0.5}(\G_*)}, \|\pa_t d_{\G_t}\|_{H^{s-0.5}(\G_*)}\big),\\
&\|D_t v\|_{H^{s-1.5}(\Om_t)}\le C\big(\|d_{\G_t}\|_{H^{s}(\G_*)}, \|\pa_t d_{\G_t}\|_{H^{s-1}(\G_*)}, \|\pa^2_t d_{\G_t}\|_{H^{s-2}(\G_*)}\big).
\end{split}
\]
Moreover, if $s\geq4$, one has 
\[
\begin{split}
&\|D_\tau v\|_{H^{s-1.5}(\Om_t)}\le C\big(\|d_{\G_t}\|_{H^{s}(\G_*)}, \|\pa_t d_{\G_t}\|_{H^{s-1}(\G_*)}\big)\Big(\|\pa_\tau d_{\G_t}\|_{H^{s-1}(\G_*)}+\|\pa_t\pa_\tau d_{\G_t}\|_{H^{s-2}(\G_*)}\Big),
\\
&\|D_tD_\tau v\|_{H^{s-3}(\Om_t)}\le  C\big(\|d_{\G_t}\|_{H^{s-1.5}(\G_*)}, \|\pa_t d_{\G_t}\|_{H^{s-1.5}(\G_*)}\big)\Big(\|\pa_\tau d_{\G_t}\|_{H^{s-1.5}(\G_*)}+ \|\pa_t\pa_\tau d_{\G_t}\|_{H^{s-1.5}(\G_*)}+ \|\pa^2_t\pa_\tau d_{\G_t}\|_{H^{s-3.5}(\G_*)}\Big).
\end{split}
\]
}
\ethm
\begin{proof} {Step 1: Estimates for $N_t$, $\tau_t$ and $\ka$}.
Firstly, we recall from Section 2 that the free surface $\G_t$ is defined by $\Phi_{S_t}(p)=p+d_{\G_t}(p)\mu(p)$ for $p\in \G_*$. Meanwhile, $N_t, \tau_t$ are defined on $\G_t$.  Parameterizing $\G_*$ with the arc length parameter $s$ and denoting $\Phi_{S_t}=\Phi_{S_t}(s)$, one knows immediately  that the unit tangential vector of this parametric curve reads
\[
\tau_t\circ\Phi_{S_t}(s)=\f{\pa_s \Phi_{S_t}(s)}{|\pa_s \Phi_{S_t}(s)|},\quad\hbox{with}\quad \pa_s \Phi_{S_t}(s) =p'(s)+\pa_s d_{\G_t}(s)\mu(s)+d_{\G_t}(s)\mu'(s).
\]
Consequently,  one obtains from \eqref{change of Dt} that
\[
(D_t\tau_t)\circ\Phi_{S_t}=D_{t*}(\tau_t\circ\Phi_{S_t})=h(D_{t^*}d_{\G_t}(s), D_{t*}\pa_s d_{\G_t}(s), d_{\G_t}(s), \pa_sd_{\G_t}(s)),
\]
with $v^*$ defined in \eqref{v* expression} and $h$ some polynomial function.  Moreover, noticing that
$\na_{v^*}f=(v^*\cdot\tau_*)\pa_s f$
for any function $f$ on $\G_*$, where $\tau_*$ is the unit tangential vector of $\G_*$, one rewrites
\[
D_{t*}f=\pa_tf+(v^*\cdot\tau_*)\pa_s f.
\]
One can see that the higher-order terms in $D_t\tau_t$ are $\pa_s\pa_t d_{\G_t}$, $\pa^2_s d_{\G_t}$, $v$, or we simply write $D_t\tau_t=h(\pa_s\pa_t d_{\G_t}, \pa^2_s d_{\G_t}, v)$, which implies immediately the desired estimate for $\tau_t, D_t\tau_t$ with Lemma \ref{trace thm PG}.
Meanwhile, $N_t, D_tN_t$ and $\ka, D_t\ka$ can be handled in a similar way.

\medskip
\noindent{Step 2: Estimates for $v$ and $D^2_t\tau_t$.}
In fact, recalling that $v=\na\phi$ with $\phi$ satisfying system \eqref{phi system}, one obtains immediately by Lemma \ref{est:elliptic} that
\[
\|v\|_{H^{s}(\Om_t)}\le C(\|d_{\G_t}\|_{H^{s}(\G_*)})\big(\|(\pa_t d_{\G_t}\mu)\circ \Phi^{-1}_{S_t}\cdot N_t\|_{H^{s-0.5}(\G_t)}+\|\xi\gamma\|_{H^{s-2}(\Om_t)}\big),
\]
which implies the desired estimate.

On the other hand, we turn to the estimate for $D_t v$. One has similar as in \eqref{D be v} that
\[
D_t v=\na D_t\phi+[D_t,\na]\phi.
\]
 Checking term by term and applying Lemma \ref{est:elliptic}, one can obtain the desired estimate.

Meanwhile,  a similar computation as in Step 1 leads to
\[
D^2_t\tau_t=h(\pa^3_s d_{\G_t}, \pa^2_s\pa_t d_{\G_t}, \pa_s\pa^2_td_{\G_t}, v, D_t v)
\]
for some polynomial function $h$. Therefore, combining the estimate for $D_tv$, one can obtain the desired estimate for $D^2_t\tau_t$. And $D^2_tN_t$ can be handled similarly as well.

In the end,   analogous analysis and applying Lemma \ref{trace thm PG} and  Lemma \ref{est:elliptic}, one can prove the desired estimates for $D_\tau v, D_t D_\tau v$.
\end{proof}

Now we  turn to the remainder term $R_0$ from equation \eqref{ka_a equation}, where the highest-order terms in $R_0$ are like $\pa^4 d_{\G_t}, \pa^3 \pa_t d_{\G_t}, \pa^2 v, \pa\mN_a$.
 \bthm{Lemma}\label{R0 estimate}{\it
 Let $\G_t\in \Lambda(S_*, 8.5, \delta, \pi/16)$. Then the remainder term $R_0$ satisfies
\[
 \|R_0\|_{H^{4}(\G_*)} \leq a^{3/2} C(\|\mN_a\|_{H^{5.5}(\G_*)}, \|\pa_t \mN_a\|_{H^{4}(\G_*)}, |d_i|, |d'_i|\big)
\]
and
\[
\|R_0\|_{H^{2.5}(\G_*)} \leq C(\|\mN_a\|_{H^{5.5}(\G_*)}, \|\pa_t \mN_a\|_{H^{4}(\G_*)}, |d_i|, |d'_i|\big),
\]
where  the index $i$ stands for the summation of $i=l, r$.

Moreover, one also has for a parameter $\tau$ the following estimate:
\[
\|\pa_{\tau}R_0\|_{H^{2.5}(\G_*)} \leq a^{3/2} C\big(\|\mN_a\|_{H^{5.5}(\G_*)}, \|\pa_t \mN_a\|_{H^{4}(\G_*)}, |d_i|, |d'_i|\big)\big( \|\pa_\tau\mN_a\|_{H^{4}(\G_*)}+\|\pa_t \pa_\tau\mN_a\|_{H^{2.5}(\G_*)}
+\|\pa_\tau d_i|+|\pa_t\pa_\tau d_i|\big).
\]
 }
 \ethm

\begin{proof}  { Step 1: Estimate for $R_0$.}
Recalling the definition of $R_0$, we have that
\[
\|R_0\circ\Phi^{-1}_{S_t}\|_{H^{4}(\G_t)} \leq \|R_a\|_{H^4(\G_*)}+ \|[D_{t}, \cN] D_t\ka\|_{H^{4}(\G_t)}+\|D_t[D_{t}, \cN] \ka\|_{H^{4}(\G_t)}+\|\cN\widetilde{R}_1\|_{H^{4}(\G_t)}.
 \]

For  the first term in the inequality above, one has directly from the expression of $R_a$ the following estimate:
\[
\|R_a\|_{H^4(\G_*)}\le C\big(a^3\|d_{\G_t}\|_{H^7(\G_*)}, a^3\|\pa_td_{\G_t}\|_{H^5(\G_*)}, a^3\|v\|_{H^{5.5}(\Om_t)}\big)
\]
where the polynomial $C$ is linear with $a^3$.

Recalling \eqref{commutator DN} and \eqref{Dt k expression},
we can have 
\[
\begin{split}
\|[D_t,\,\cN]D_t\ka\|_{H^{4}(\G_t)}\leq& \|\D^{-1}\Big(2\na v\cdot \na^2(D_t\ka)_{\cH}+\D v\cdot\na (D_t\ka)_{\cH},\,(\na_{N_b}v-\na_{v}N_b)\cdot\na (D_t\ka)_{\cH}\Big)\|_{H^{4}(\Om_t)}\\
&+ \|\na_{N_t}v\cdot \na (D_t\ka)_{\cH}\|_{H^{4}(\G_t)}+\|\na_{\tau_t} (D_t\ka)_{\cH}(\na_{\tau_t} v\cdot N_t)\|_{H^{4}(\G_t)}\\
\leq&C(\|d_{\G_t}\|_{H^{6}(\G_*)})\|v\|_{H^{5.5}(\Om_t)}\,\|D_t\ka\|_{H^{5}(\G_t)}\\
\leq& C\big(\|d_{\G_t}\|_{H^{8}(\G_*)}, \|\pa_td_{\G_t}\|_{H^{7}(\G_*)}\big),
\end{split}
\]
where Lemma \ref{trace thm PG},  Lemma \ref{est:elliptic} and Lemma \ref{Dt v estimate} are applied.

On the other hand,  for the term $\|D_t[D_{t}, \cN] \ka\|_{H^{4}(\G_t)}$,  one denotes
\[
 w=\D^{-1}\big(2\na v\cdot \na^2\ka_{\cH}+\D v\cdot\na \ka_{\cH},\,(\na_{N_b}v-\na_{v}N_b)\cdot\na \ka_{\cH}\big)
\]
to have
\[
D_t[D_{t}, \cN] \ka=D_t\big(\na_{N_t}w\big)-D_t\big(\na_{N_t}v\cdot \na \ka_{\cH})-D_t\big(\na_{(\na \ka_{\cH})^\top} v\cdot N_t\big).
\]
Consequently, remembering that we used $D_tv=-\si\na\ka_\cH-\na P_{v,v}+{\bf g}$,  one can prove
 \[
 \begin{split}
 \| D_t[D_{t}, \cN] f\|_{H^{4}(\G_t)}
 \leq& \big\|[D_t,  \na_{N_t}]w \big\|_{H^{4}(\Gamma_t)}+\big\|  \na_{N_t}D_tw \big\|_{H^{4}(\Gamma_t)}+\|D_t(\na_{N_t}v\cdot \na \ka_{\cH})\|_{H^{4}(\Gamma_t)}\\
 &\quad +\|D_t(\na_{(\na g_{\cH})^\top} v\cdot N_t)\|_{H^{4}(\Gamma_t)}\\
\le & C(\|d_{\G_t}\|_{H^{6}(\G_*)}, \|v\|_{H^{5.5}(\Om_t)}, \|D_t v\|_{H^{5.5}(\Om_t)})
\big(\|\ka\|_{H^{5}(\G_t)}+\|D_t\ka\|_{H^{5}(\G_t)}\big)\\
\le & C\big(\|d_{\G_t}\|_{H^{8}(\G_*)}, \|\pa_t d_{\G_t}\|_{H^{7}(\G_*)}\big),
\end{split}
\]
where one used the following estimate according to \eqref{commutator Dt D-1}:
\[
 \begin{split}
\|D_tw\|_{H^{5.5}(\Om_t)}\le& C(\|d_{\G_t}\|_{H^{5}(\G_*)})\Big(\big\|D_t(2\na v\cdot \na^2\ka_{\cH}+\D v\cdot\na \ka_{\cH})\big\|_{H^{3.5}(\Om_t)}+\big\|D_t\big((\na_{N_b}v-\na_{v}N_b)\cdot\na \ka_{\cH}\big)\big\|_{H^{4}(\G_b)}\\
&\qquad\quad +\big\|2\na v\cdot \na^2w+\D v\cdot\na w\big\|_{H^{4}(\Om_t)}+\big\|(\na_{N_b}v-\na_{v}N_b)\cdot\na w\big\|_{H^{4.5}(\G_b)}\Big)\\
\le &  C(\|d_{\G_t}\|_{H^{6}(\G_*)}, \|v\|_{H^{5.5}(\Om_t)}, \|D_t v\|_{H^{5.5}(\Om_t)})\big(\|\ka\|_{H^{5}(\G_t)}+\|D_t \ka_\cH\|_{H^{5.5}(\Om_t)}\big),
\end{split}
\]
and also
\[
\begin{split}
\|D_t \ka_\cH\|_{H^{5.5}(\Om_t)}&\le \|\cH(D_t \ka)\|_{H^{5.5}(\Om_t)}+\|\ka\|_{H^{5.5}(\Om_t)}\\
&\le C(\|d_{\G_t}\|_{H^{6}(\G_*)})\big(\|D_t \ka\|_{H^{5}(\G_t)}+\|v\|_{H^{5.5}(\Om_t)}\|\ka\|_{H^{5}(\G_t)}\big),
\end{split}
\]
combining  \eqref{commutator Dt H} and Lemma \ref{est:elliptic}.
Therefore,  applying Lemma \ref{Dt v estimate} again, we arrive at
\[
\|[D_{t}, \cN] D_t\ka\|_{H^{4}(\G_t)}+\|D_t[D_{t}, \cN] \ka\|_{H^{4}(\G_t)}\le C\big(\|d_{\G_t}\|_{H^{8}(\G_*)}, \|\pa_td_{\G_t}\|_{H^{7}(\G_*)}\big).
 \]

\medskip

 To finish the estimate for $\|R_0\|_{H^{4}(\G_*)}$, we still  need to deal with $\|\cN\widetilde{R}_1\|_{H^{4}(\G_t)}$. In fact, from the expression of $\widetilde {R}_1$, we know that the higher-order terms are like $\pa_td_{\G_t}$, $\pa^2 d_{\G_t}$,  $\na \kappa_\cH$, $\pa^2 v$, $\D_{\G_t}N_t$ and $\Delta_{\Gamma_t}\na P_{v,v}$. For example, from the definition of $P_{v,v}$, we have
 \beno
 \|\na P_{v, v}\|_{H^{7}(\Gamma_t)} \leq C \big(\|d_{\Gamma_t}\|_{H^{8.5}(\Gamma_*)})(1+\|v\|_{H^{7}(\Gamma_t)})^2\leq C(\|d_{\Gamma_t}\|_{H^{8.5}(\Gamma_*)}, \|\pa_td_{\Gamma_t}\|_{H^{7}(\Gamma_*)}\big).
 \eeno

 Then, checking term by term, it is easy to show that
\beno
\|\cN\widetilde{R}_1\|_{H^{4}(\G_t)}\leq C(\|d_{\G_t}\|_{H^{8}(\G_*)}, \|\pa_t d_{\G_t}\|_{H^{7}(\G_*)} \big).
\eeno
Summing these estimates up, we have the estimate for $R_0$:
\[
\|R_0\|_{H^{4}(\G_*)} \leq C(\|d_{\G_t}\|_{H^{8}}, \|\pa_t d_{\G_t}\|_{H^{7}(\G_*)},  a^3\|d_{\G_t}\|_{H^7(\G_*)}, a^3\|\pa_td_{\G_t}\|_{H^5(\G_*)}, a^3\|v\|_{H^{5.5}(\Om_t)}\big),
\]
and notice here that the right side is linear with respect to $a^3$. The desired estimate for $\|R_0\|_{H^{4}(\G_*)}$ follows from Proposition \ref{d and  N ka},  Proposition \ref{pa t d estimate} and  Lemma \ref{Dt v estimate}. 

Moreover, one  has by an analogous analysis that
\[
\|R_0\|_{H^{2.5}(\G_*)} \leq C(\|d_{\G_t}\|_{H^{6.5}}, \|\pa_t d_{\G_t}\|_{H^{5.5}(\G_*)},  a^3\|d_{\G_t}\|_{H^{5.5}(\G_*)}, a^3\|\pa_td_{\G_t}\|_{H^{3.5}(\G_*)}, a^3\|v\|_{H^{4}(\Om_t)}\big),
\]
which leads to the desired estimate  by using again Proposition \ref{d and  N ka},  Proposition \ref{pa t d estimate} and  Lemma \ref{Dt v estimate}.

\medskip

\noindent{Step 2: Estimate for $\pa_{\tau} R_0$.} The analysis is similar as before, so we omit the details. In fact, since
\[
R_0=R_0\big(\pa^4 d_{\G_t}, \pa^3\pa_t d_{\G_t}, \pa^2 v, \pa\mN_a, a^3\pa^3 d_{\G_t}, a^3\pa\pa_t d_{\G_t}, a^3\pa v\big)
\]
and  $\pa_\tau$ acting on $R_0$ results in extra $\pa_\tau$ (or $D_\tau$ on $\G_t$),  similar analysis as before leads to  the desired estimates.
\end{proof}

\medskip

On the other  hand,  we prove the estimate of the reminder term $R_{c, i}$.
 \bthm{Lemma}\label{lem:R_c}
{\it Let $\G_t\in \Lambda(S_*, 8.5, \delta, \pi/16)$ and $\tau$ be a parameter. Then  the remainder term $R_{c, i}$  ($i=l, r$) from \eqref{equ:cca} satisfies
 \[
\big\|D_{t_*}R_{c, i}\big\|_{H^1(\G_*)} \leq  a\, C\big(\|\mN_a\|_{H^{5.5}(\G_*)}, \|\pa_t\mN_a\|_{H^{4}(\G_*)}, \|\pa^2_t\mN_a\|_{H^{2.5}(\G_*)}, |d_i|, |d'_i|, |d''_i|\big),
 \]
and
 \[
 \begin{split}
 \big\|\pa_\tau R_{c, i}\big\|_{H^1(\G_*)} &\leq  a \,C\big(\|\mN_a\|_{H^{5.5}(\G_*)}, \|\pa_t\mN_a\|_{H^{4}(\G_*)},  |d_i|, |\pa_td_i|\big)\times\\
 &\qquad\quad\Big(\|\pa_\tau\mN_a\|_{H^{4}(\G_*)}+\|\pa_t\pa_\tau\mN_a\|_{H^{2.5}(\G_*)}+|\pa_\tau d_i|+ |\pa_t\pa_\tau d_i|\Big).
\end{split} \]
}
 \ethm
\begin{proof}
The proof is analogous to the proof of the previous lemma, so we omit most of the details. Firstly, since one knows from Lemma \ref{lem:cc} that
\[
R_{c, i}=R_{c, i}\big(\pa^5d_{\G_t}, \pa^4\pa_td_{\G_t}, \pa^3 v, \pa^2\mN_a, \pa\pa_t\mN_a, a^3\pa^4d_{\G_t}, a^3\pa^2\pa_td_{\G_t}, a^3\pa^2 v\big),
\]
one has
\beq\label{estimate middle for R ci}
\begin{split}
\big\|R_{c, i}\big\|_{H^1(\G_*)}&\le  C\big(\|d_{\G_t}\|_{H^6(\G_*)}, \|\pa_t d_{\G_t}\|_{H^5(\G_*)}, \|v\|_{H^{4.5}(\Om_t)}, \|\mN_a\|_{H^3(\G_*)}, \|\pa_t\mN_a\|_{H^2(\G_*)},\\
&\qquad a^3\|d_{\G_t}\|_{H^5(\G_*)}, a^3\|\pa_t d_{\G_t}\|_{H^3(\G_*)},   a^3\|v\|_{H^{3.5}(\Om_t)}\big)
\end{split}
\eeq
where the right side is linear with respect to $a^3$.

On the other hand, we know that $D_{t*}$ acting on $R_{c, i}$ results in extra `$\pa_t$' and `$\na_{v^*}$' terms in each term of $R_{c, i}$, which implies
\[
\begin{split}
\big\|D_{t*} R_{c, i}\big\|_{H^1(\G_*)}
&\le C\big(\|d_{\G_t}\|_{H^7(\G_*)}, \|\pa_t d_{\G_t}\|_{H^6(\G_*)}, \|\pa^2_t d_{\G_t}\|_{H^5(\G_*)},
 \|v\|_{H^{5.5}(\Om_t)}, \|D_tv\|_{H^{4.5}(\Om_t)},
\|\mN_a\|_{H^4(\G_*)}, \\
&\qquad \quad\|\pa_t\mN_a\|_{H^3(\G_*)},   \|\pa^2_t\mN_a\|_{H^2(\G_*)}, a^3\|d_{\G_t}\|_{H^6(\G_*)}, a^3\|\pa_t d_{\G_t}\|_{H^5(\G_*)}, a^3\|\pa^2_t d_{\G_t}\|_{H^3(\G_*)}, \\
&\qquad\quad a^3\|v\|_{H^{4.5}(\Om_t)}, a^3\|D_tv\|_{H^{3.5}(\Om_t)}\big)
\end{split}
\]
and the right side is linear with $a^3$ as well. 
Consequently, applying Proposition \ref{d and  N ka},  Proposition \ref{pa t d estimate} and Lemma \ref{Dt v estimate}, we finish the proof for the first estimate.

On the other hand, the estimate for $\pa_\tau R_{c, i}$ follows in a similar way while notice that $\pa_\tau$ only adds extra $\pa_\tau$ on each term in \eqref{estimate middle for R ci} (except $D_\tau$ on $v$).
\end{proof}

\medskip

\subsection{Estimates for time derivatives of $\Gamma_t$}

Since the norms involving time derivatives of $\mN_a, d_i$ will be used in the energy estimates,  we  need  to consider  the estimates for $\pa_td_{\G_t}$, $\pa^2_td_{\G_t}$ in terms of times derivatives of  $\mN_a$ and $d_i$.

\bthm{Propsition}\label{pa t d estimate}{\it  Let $\G_t\in \Lam(S_*, 8.5, \del, \pi/16)$ for $t\in [0,T]$ and $s_1, s_2$ satisfy $8.5 \ge s_1\ge s_2\ge 1$ and $s_1\le s_2+3$. 

 Then the following estimates hold:
\[
\begin{split}
&\|\pa_td_{\G_t}\|_{H^{s_1}(\G_*)}\le a^{s_1-s_2-3} \,C\big(\big\|\mN_a-\mN_{a*}\big\|_{H^{s_2}(\G_*)}, |d_i|\big)\big(\big\|\pa_t\mN_a\big\|_{H^{s_2}(\G_*)}+\big|d'_i\big|\big),\\
&
\|\pa^2_td_{\G_t}\|_{H^{s_1}(\G_*)}\le a^{s_1-s_2-3} \,C\big(\big\|\mN_a-\mN_{a*}\big\|_{H^{s_2}(\G_*)}, \big\|\pa_t\mN_a\big\|_{H^{s_2}(\G_*)}, |d_i|, \big|d'_i\big|\big)\big(\big\|\pa^2_t\mN_a\big\|_{H^{s_2}(\G_*)}+\big|d''_i\big|\big),
\end{split}
\]
where we use the index $i$ for  the summation on $i=l, r$ for the sake of convenience.
}
\ethm
\begin{proof} The proof follows from the proof of Proposition \ref{d and  N ka}. In fact, denoting
\[
v_1=(\pa_t d_{\G_t}\mu)\circ\Phi^{-1}_{S_t}\quad \hbox{and}\quad D_{t1}=\pa_t+\na_{v_1},
\]
one has from the definition of $\mN_a$ that
\[
\begin{split}
\pa_t\mN_a&=a^3\pa_t d_{\G_t}+D_{t1}\cN(\ka)\circ\Phi_{S_t}\\
&=a^3\pa_t d_{\G_t}+\cN\big(D_{t1}\ka\big)\circ\Phi_{S_t}+\big([D_{t1}, \cN]\ka\circ \Phi^{-1}_{S_t}\big)\circ\Phi_{S_t}.
\end{split}
\]

On the other hand,  expressing $D_{t1}\ka$ similarly as in Proposition \ref{d and  N ka}  and substituting it into the expression of $\pa_t\mN_a$, one obtains a linear system for $\pa_td_{\G_t}$ combining \eqref{point condition for d}:
\beq\label{eqn for pa t d}
\left\{
\begin{array}{ll}
a^3\pa_td_{\G_t}-\cN\Big(\big((\mu\circ \Phi^{-1}_{S_t}\big)\cdot N_t\big)\D_{\G_t}(\pa_td_{\G_t}\circ\Phi^{-1}_{S_t}) \Big)\circ\Phi_{S_t}+r_{\pa_t d}=\pa_t\mN_a,\\
\pa_t d_{\G_t}(p_{i*})=\pa_t d_i,\quad i=l, r,
\end{array}\right.
\eeq
where $r_{\pa_t d}$ contains remainder terms of $\pa_t d_{\G_t}$ and $d_{\G_t}$, which  can be written explicitly.

Noticing that $(\mu\circ \Phi^{-1}_{S_t}\big)\cdot N_t\ge c_1>0$ for $t\in [0,T]$ and some constant  $c_1>0$,  so when $a$ is large enough, standard elliptic analysis as before leads to  the desired estimates for $\pa_t d_{\G_t}$.  In the end, one can also prove similar estimates for $\pa^2_t d_{\G_t}$ and the proof can be finished.
\end{proof}

In the following sections, one also needs to consider the estimates for $\pa_\tau d_{\G_t}$, $\pa_t\pa_\tau d_{\G_t}$ and $\pa^2_t\pa_\tau d_{\G_t}$.
\bthm{Corollary}\label{pa tau d estimate}{\it  Let $\G_t\in \Lam(S_*, 8.5, \del, \f{\pi}{16})$ for $t\in [0,T]$ and $s_1, s_2$ satisfy $8.5 \ge s_1\ge s_2\ge 1$ and $s_1\le s_2+3$. 

Then the following inequalities hold:
\[
\begin{split}
&\|\pa_\tau d_{\G_t}\|_{H^{s_1}(\G_*)}\le a^{s_1-s_2-3}C\big(\big\|\mN_a-\mN_{a*}\big\|_{H^{s_2}(\G_*)}, |d_i|\big)\Big(\big\|\pa_\tau\mN_a\big\|_{H^{s_2}(\G_*)}+\big|\pa_\tau d_i\big|\Big),\\
&\|\pa_t\pa_\tau d_{\G_t}\|_{H^{s_1}(\G_*)}\le a^{s_1-s_2-3}C\big(\big\|\mN_a-\mN_{a*}\big\|_{H^{s_2}(\G_*)}, \big\|\pa_t\mN_a\big\|_{H^{s_2}(\G_*)} , |d_i|, |\pa_t d_i|\big)\times\\
&\qquad\qquad\qquad  \qquad\quad \Big(\big\|\pa_\tau\mN_a\big\|_{H^{s_2}(\G_*)}+\big\|\pa_t\pa_\tau\mN_a\big\|_{H^{s_2}(\G_*)}+\big|\pa_\tau d_i\big|+\big|\pa_t\pa_\tau d_i\big|\Big);
\end{split}
\]
Moreover, one has
\[
\begin{split}
\|\pa^2_t\pa_\tau d_{\G_t}\|_{H^{s_1}(\G_*)}
&\le a^{s_1-s_2-3} C\big(\big\|\mN_a-\mN_{a*}\big\|_{H^{s_2}(\G_*)}, \big\|\pa_t\mN_a\big\|_{H^{s_2}(\G_*)}, \big\|\pa^2_t\mN_a\big\|_{H^{s_2}(\G_*)},  |d_i|, |\pa_t d_i|, |\pa^2_t d_i|\big)\times\\
&\quad
 \Big( \big\|\pa_\tau\mN_a\big\|_{H^{s_2}(\G_*)}+\big\|\pa_t\pa_\tau\mN_a\big\|_{H^{s_2}(\G_*)}+\big\|\pa^2_t\pa_\tau\mN_a\big\|_{H^{s_2}(\G_*)}+ \big|\pa_\tau d_i\big|+\big|\pa_t\pa_\tau d_i\big|+\big|\pa^2_t\pa_\tau d_i\big|\Big).
\end{split}
\]
}
\ethm
\begin{proof} The proof is analogous to the proof of Proposition \ref{pa t d estimate} and is  left to the readers.

\end{proof}

\bigskip

 \section{The linear problem on the free surface}

Assuming that the free surface $\G_t\in \Lambda(S_*, 8.5, \delta, \f{\pi}{16})$ is known already, we consider the following linear system of $(f(t,p), d_l(t), d_r(t))$ with $p\in \G_*$:
 \beq\label{equ:linear}
\left\{\begin{array}{ll}
D_{t*}^2 f+\sigma \cA (d_{\G_t})f=g_1\qquad\hbox{on}\quad \G_*\\
D_{t*}\cA(d_{\G_t}) f\pm\f{\sigma^2}{\beta_c}( \sin \om_i)^2\circ\Phi_{S_t}  (\na_{\tau_t}(\cA(d_{\G_t})f\circ \Phi_{S_t}^{-1}))\circ\Phi_{S_t}= g_{2,i}\quad \textrm{at} \quad p_{i*},\\
d''_i(t)=B_i,\qquad i=l, r
\end{array}
\right.
\eeq
with initial data
\[
f|_{t=0}=f_0, \  D_{t*} f|_{t=0}=f_1 ,\ d_i(0)=d_{i, 0},\ d'_i(0)=d_{i, 1}.
\]
Here the right sides  $g_1,\, g_2, B_i$ and the initial data $f_0, f_1, d_{i, 0}, d_{i, 1}$ are given $(i=l, r)$.

To solve the  linear system above, we introduce 
\[
F\triangleq (a^3+\cA (d_{\G_t}))f\triangleq\cA_a f.
\]
This makes sense thanks to Remark \ref{inverse Aa}, so we could retrieve $f$ from $F$ as long as we have Dirichlet boundary conditions for $f$.

 By direct calculations, we get
 \beq\label{equ:linear1}
\left\{\begin{array}{ll}
D_{t*}^2 F+\sigma \cA (d_{\G_t})F= \cA_a g_1+[\cA_a,D_{t*}^2 ]f\qquad\hbox{on}\quad \G_*,\\
D_{t*}  F\pm\f{\sigma^2}{\beta_c}( \sin \om_i)^2\circ\Phi_{S_t}  (\na_{\tau_t}( F\circ \Phi_{S_t}^{-1}))\circ\Phi_{S_t}= g_{2,i}\quad \textrm{at} \quad p_{i*}.
\end{array}
\right.
\eeq
Consequently, system \eqref{equ:linear} is equivalent to the following system
 \beq\label{equ:linear2}
\left\{\begin{array}{ll}
D_{t*}^2 F+\sigma \cA (d_{\G_t})F= \cA_a g_1+[ \cA_a,D_{t*}^2 ]f\qquad\hbox{on}\quad \G_*,\\
D_{t*}  F\pm\f{\sigma^2}{\beta_c}( \sin \om_i)^2\circ\Phi_{S_t}  (\na_{\tau_t}( F\circ \Phi_{S_t}^{-1}))\circ\Phi_{S_t}= g_{2,i}\quad \textrm{at} \quad p_{i*},\\
(a^3+\cA (d_{\G_t}))f=F,\\
D_{t*}^2 f\pm\sigma \cA (d_{\G_t})f=   g_1\quad \textrm{at} \quad p_{i*},\\
d''_i(t)=B_i,\qquad i=l, r.
\end{array}
\right.
\eeq

\bthm{Remark}
Here, the reason of using  higher-order boundary conditions in \eqref{equ:linear} is that we need high-order energy estimates of $f$ later,  while second or higher-order time derivatives on the right side should be avoided. In fact, we  derived the high-order boundary conditions containing $D_tv$ and $\pa_t^2 d_{\Gamma_t}$ in the right-side part $g_2$ in the previous section, where Euler equation is applied to replace them. In contrast, we do not have Euler equation in the iteration scheme, which means we have to avoid second or higher order of time derivatives on the right side from now on. As a result,  we give the high-order boundary conditions in the first place. \ethm

We firstly prove the existence of the solution to our linear problem.
\bthm{Proposition}\label{existence of linear system}{\it
Let $(f_0, f_1)\in H^{5.5}(\Gamma_*)\times H^{4}(\Gamma_*)$,  $g_1\in C\big([0,T];H^{4}(\Gamma_{*})\big)$, $g_{2,i}\in C^1\big([0,T]; H^1(\G_*)\big)$  and $d_{i, 0}, d_{i, 1} \in\R$ be given and $\om_*\in(0, \f{\pi}{16})$. If the compatibility conditions at $t=0$ hold, namely
\beno
D_{t*}  \cA (d_{\G_t})f\big|_{t=0}\pm\f{\sigma^2}{\beta_c}( \sin \om_i)^2\circ\Phi_{S_t}  (\na_{\tau_t}( F\circ \Phi_{S_t}^{-1}))\circ\Phi_{S_t}\big|_{t=0}= g_{2, i}\big|_{t=0}\quad \textrm{at} \quad p_{i*},
\eeno
then there exists some $T>0$ such that the system \eqref{equ:linear} has a unique solution on $[0, T]$.
}
\ethm
\begin{proof}
To solve system \eqref{equ:linear}, we need to solve system \eqref{equ:linear2} equivalently. 

First of all, thanks to Remark \ref{inverse Aa},  we notice that the commutator $[\cA_a,D_{t*}^2 ]f$ is a lower-order term of $F, D_{t*}F$, and we can write it as $\cR_1F+\cR_0D_{t*}F$, where $\cR_1$ means at most first-order derivative of $F$ and $\cR_0$ means at most zero-order derivative of $D_{t*}F$.  The operator $\cA_a$ is invertible since Dirichlet boundary conditions for $f$ are given in \eqref{equ:linear2} by $D_{t*}^2 f\pm\sigma \cA (d_{\G_t})f=   g_1$. 

Consequently,  to finish the proof, we only need to prove the existence of the solution $F$ to the following system
\beq\label{equ:linear3}
\left\{\begin{array}{ll}
D_{t*}^2 F+\sigma \cA_aF= \cR_1F+\cR_0D_{t*}F+\cA_a g_1\qquad\hbox{on}\quad \G_*\\
D_{t*}  F\pm\f{\sigma^2}{\beta_c}( \sin \om_i)^2\circ\Phi_{S_t}  (\na_{\tau_t}( F\circ \Phi_{S_t}^{-1}))\circ\Phi_{S_t}= g_{2, i}\quad \textrm{at} \quad p_{i*},
\end{array}
\right.
\eeq
with initial data (which can be derived from \eqref{equ:linear} and are omitted here)
\[
F|_{t=0}=F_0, \  D_{t*} F|_{t=0}=F_1.
\]
The equation here is a linear hyperbolic equation of third order with mixed boundary conditions.
To solve this system, we follow the proof of \cite{MI}, while we keep in mind that the difference between our paper and \cite{MI} lies in the third-order elliptic operator $\cA_a=\cA(d_{\G_t})+a^3$ studied in Remark \ref{inverse Aa}. As a result, we mainly focus on the estimates related to $\cA_a$, and we only outline the proof for the remainder parts.

%In fact, $\cA(d_{\G_t})$ satisfies  $\lambda+\cA(d_{\G_t})\geq0$ when $\lambda\geq \lambda_0$ for some $
%\lambda_0$. Therefore, following the steps in \cite{MI}, we can show that the system \eqref{equ:linear} is locally %well-posed in $H^{s+1.5}(\Gamma_*)\times H^{s}(\Gamma_*)$ and the details are omitted here.
\medskip
Introducing 
\[
U(t)=\big(F(t, p), D_{t*}F(t, p)\big)^*,
\]
we rewrite system \eqref{equ:linear3} into an equivalent system again:
\beq\label{semigroup system}
\begin{split}
&\f d{dt}U(t)=\tilde \cA(t)U(t)+G(t)\quad \hbox{on}\ \G_*\\
&\cB(t)U(t)=D_{t*}g_{2,i},\quad\hbox{at}\ p_{i*},\quad U(0)=U_0,
\end{split}
\eeq
where 
\[
G(t)=\left(\begin{array}{ll}
0\\
\cA_ag_1
\end{array}\right),\quad 
\tilde \cA(t)=\left(\begin{matrix}
-\na_{v^*} & 1\\
-\si \cA_a+\cR_1 &-\na_{v^*}+\cR_0
\end{matrix}\right)
\]
and 
\[
\cB(t)=\left(\begin{matrix}
-\si \cA_a+\cR_1, 
& c_0\big(\na_{\tau_t}(\cdot\circ \Phi^{-1}_{S_t})\big)\circ \Phi_{S_t}
\end{matrix}\right)
\]
with the notation
\[
c_0=c_0(t)=\f{\sigma^2}{\beta_c}( \sin \om_i)^2\circ\Phi_{S_t}> 0.
\]
Moreover, we recall that  $\cR_1=\cR_1(\pa)$ is some first-order operator, $\cR_0=\cR_0(\pa)$ is some zero-order operator on $\G_*$
(the details here can be computed in a direct way and are omitted since they are lower-order terms as mentioned in the beginning of the proof). The boundary condition $\cB(t)U(t)=D_{t*}g_{2,i}$ is derived from the boundary condition of $F$ by taking $D_{t*}$ and applying the equation of $F$, since we use $H^{2.5}(\G_*)$ norm for $F$ as a starting point for proving the existence.

We plan to prove firstly the existence of the solution to system \eqref{semigroup system} when the time $t=t_0$ is fixed in $\tilde \cA(t)$ and the boundary satisfies $\cB(t)U=0$, where we can  use  standard semigroup method as in Lemma 2.4 of \cite{MI}. As the second step, we can build energy inequality for this system and prove the existence of the solution to our system \eqref {equ:linear3} of $F$  when the time is fixed in the operators  as well. At last, we can use the method of Cauchy's polygonal line to recover the existence of the solution to system \eqref {equ:linear3} as in section 4 of \cite{MI}.

To finish the first step and use semigroup method for proving the existence, the key point lies in  the estimates for $\tilde \cA(t)$ under a proper norm. In fact,  for $U=(u, v)^*\in H^{2.5}(\G_*)\times H^1(\G_*)$, we define the following norm
\[
\|U\|_{H}=\big(U, U\big)_{H},
\]
where
\[
\begin{split}
\big(U_1, U_2\big)_{H}\triangleq&\si\big((\cN\Delta_{\G_t}\bar u_1)\circ\Phi_{S_t}, \, (\Delta_{\G_t}\bar u_2)\circ\Phi_{S_t}\big)
+(u_1, u_2)+\big((\na_{\tau_t}\bar v_1)\circ\Phi_{S_t},\,(\na_{\tau_t}\bar v_1)\circ\Phi_{S_t}\big)\\
&+(v_1, v_2)
\end{split}
\]
with $U_j=(u_j, v_j)^*$, $\bar U_j=U_j\circ \Phi^{-1}_{S_t}$ and $(\cdot, \cdot)$  the $L^2$ inner product on $\G_*$. 

Similarly as in Lemma 2.2 \cite{MI}, we want to show that the following estimate holds for some constant $C$ depending on $\G_*, \G_t$:
\beq\label{A semigroup estimate}
\big(\tilde \cA(t)U, U\big)_H+\big(U, \tilde \cA(t)U\big)_H\le C\Big((U, U)_H+\sum_i|\cB(t)U|^2\big|_{p_{i*}}\Big).
\eeq
As long as we have this estimate, we can prove as in Lemma 2.3 \cite{MI} that
\[
\|(\lam I-\tilde \cA(t))^{-1}\|_H\le \f 1{\lam-\lam_0}
\]
for any $\lam>\lam_0$ with some constant $\lam_0>0$. As a result, when $\cB(t)U=0$ at $p_{i*}$, we can prove the existence of the solution to system \eqref{semigroup system} when the time $t=t_0$ is fixed in $\tilde \cA(t)$ by Hille-Yosida's theorem.

Now it remains to focus on the proof for \eqref{A semigroup estimate}. 
In fact, one has by a direct calculation that
\[
\begin{split}
&\big(\tilde \cA(t)U, U\big)_H+\big(U, \tilde \cA(t)U\big)_H\\
&=2\si\int_{\G_t}\cN\Delta_{\G_t}\bar v\ \Delta_{\G_t}\bar vds-2\si
\int_{\G_t}\cN\Delta_{\G_t}\na_{v^*}\bar u\,\Delta_{\G_t}\bar uds+2\int_{\G_*}u\ vds-2\int_{\G_*}\na_{v^*}\bar u\,\bar uds\\
&\quad+ 2\int_{\G_t}\na_{\tau_t}\bar v\,\na_{\tau_t}\big(-\si\cA_au+\cR_1u-2\na_{v^*}v+\cR_0v\big)\circ\Phi^{-1}_{S_t}ds
+2\int_{\G_*}\bar v\,\big(-\si\cA_au+\cR_1u-2\na_{v^*}v+\cR_0v\big)ds\\
&:=I_1+I_2+\cdots+I_6
\end{split}
\]
where we need to pay attention to the higher-order terms $I_1, I_2, I_5$ , since the boundary condition will be used there. The other terms can be handled directly and hence are omitted here.

For $I_5$, the main part is the $\cA_a$ part, so an integration by parts leads to
\[
\begin{split}
I_5&=2\si\int_{\G_t}\na_{\tau_t}\bar v\ \na_{\tau_t}\cN\Delta_{\G_t}\bar u ds+l.o.t.\\
&=2\si\na_{\tau_t}\bar v\ \cN\Delta_{\G_t}\bar u\big|^{p_{l}}_{p_{r}}-I_1+l.o.t..
\end{split}
\]
Due to the boundary condition $\cB(t)U=D_{t*}g_{2,i}$, we know immediately that 
\[
\na_{\tau_t}\bar v=\pm\f 1{c_0}\big(\si\cA_a u-\cR_1u+\cB(t)U \big)\circ \Phi^{-1}_{S_t}\qquad\hbox{at}\ p_i (i=l, r),
\]
which implies
\[
2\si\na_{\tau_t}\bar v\ \cN\Delta_{\G_t}\bar u\big|^{p_{l}}_{p_{r}}=-\f{\si^2}{c_0}\sum_i\big(\cN\Delta_{\G_t}\bar u\big)^2\big|_{p_i}-\f1{c_0}\sum_i\cR_1u\,\cN\Delta_{\G_t}\bar u\big|_{p_i}+\f1{c_0}\sum_i\cB(t)U\circ\Phi^{-1}_{S_t}\,\cN\Delta_{\G_t}\bar u\big|_{p_i}
\]
Substituting this into $I_5$, we arrive at
\[
I_5\le -I_1-\f12\f{\si^2}{c_0}\sum_i\big(\cN\Delta_{\G_t}\bar u\big)^2\big|_{p_i}+C\big(\|u\|^2_{H^{2.5}(\G_*)}+\sum_i|\cB(t)U|^2\big|_{p_{i*}}\big)
\]
for some constant $C$ depending on $\G_*$ and $\G_t$, so it means 
\[
I_1+I_5\le -\f12\f{\si^2}{c_0}\sum_i\big(\cN\Delta_{\G_t}\bar u\big)^2\big|_{p_i}+C\big(\|u\|^2_{H^{2.5}(\G_*)}
+\sum_i|\cB(t)U|^2\big|_{p_{i*}}\big).
\]
Moreover, $I_2$ can also be handled by using integration by parts, the trace theorem and combining the inequality above.

As a result, we can conclude the estimate \eqref{A semigroup estimate}, where we also use  elliptic properties for $\cA_a$ from Remark \ref{inverse Aa}.

In this way, we can continue to follow the proof in \cite{MI}, and our proof can be finished.

\end{proof}

%The case without boundary conditions has in fact been proved in \cite{SZ2} (there is a difference in a lower-order term).

\medskip

Next, we give the energy estimates for system \eqref{equ:linear}.

Denoting by
\[
\bar f=f\circ \Phi^{-1}_{S_t},
\]
one has $\bar f$ defined on $\G_t$ and
\[
D_{t*}f=(D_t \bar f)\circ \Phi_{S_t},\quad D^2_{t*}f=(D^2_t\bar f)\circ\Phi_{S_t},\quad\hbox{and}\quad  \cA(d_{\G_t})f=-(\cN\D_{\G_t}\bar f)\circ \Phi_{S_t}.
\]
Consequently, \eqref{equ:linear} is equivalent to the following  linear problem of $(\bar f, d_l,  d_r)$:
 \beq\label{equ:linear on Gt}
\left\{\begin{array}{ll}
D_{t}^2 \bar f-\sigma \cN\D_{\G_t}\bar f=\bar g_1\qquad\hbox{on}\quad \G_t,\\
D_{t}\cN\D_{\G_t} \bar f\pm\f{\sigma^2}{\beta_c}( \sin \om_i)^2\na_{\tau_t}\cN\D_{\G_t}\bar f= \bar g_{2,i}\quad \textrm{at} \quad p_{i},\\
d''_i(t)=B_i,\quad i=l, r,
\end{array}
\right.
\eeq
with initial data
\[
\bar f|_{t=0}=f_0\circ \Phi^{-1}_{S_0},\  D_t \bar f|_{t=0}=f_1\circ\Phi_{S_0},\  d_i(0)=d_{i, 0},\ d'_i(0)=d_{i, 1}.
\]
Here we used the notation
\[
\bar g_1=g_1\circ \Phi^{-1}_{S_t},\quad \bar g_{2,i}=g_{2,i}\circ \Phi^{-1}_{S_t}.
\]
\medskip

The  higher-order energy $E_h(t,\bar f,\pa_t\bar f)$ and the higher-order dissipation $F_{h, i}(t, \bar f, \pa_t \bar f)$  for  system \eqref{equ:linear on Gt} are defined by
\[
\begin{split}
E_h(t, \bar f, \pa_t \bar f)&= \big\|\na_{\tau_t}\cN\D_{\G_t}D_t\bar f\big\|^2_{L^2(\G_t)}+\si\big\|\na \cH(\D_{\G_t}\cN\D_{\G_t}\bar f) \big\|^2_{L^2(\Om_t)}+\|D_t\bar f\|^2_{L^2(\Gamma_t)} +\|\bar f\|^2_{L^2(\Gamma_t)},\\
F_{h, i}(t, \bar f, \pa_t \bar f)&= ( \sin \om_i)^2\big|\na_{\tau_t}\cN\D_{\G_t}D_t\bar f\big|^2 ,\qquad \textrm{at}\quad p_i\ (i=l, r).
\end{split}
\]

Meanwhile, we also define the lower-order energy and dissipation
\beq\label{lower energy}
\begin{split}
E_l(t, \bar f, \pa_t \bar f)&=\big\|\na \cH(D_t\D_{\G_t}\bar f) \big\|^2_{L^2(\Om_t)}
+\si  \big\|\na_{\tau_t}\cN\D_{\G_t}\bar f\big\|^2_{L^2(\G_t)}+\|D_t\bar f\|^2_{L^2(\Gamma_t)}+\|\bar f\|^2_{L^2(\Gamma_t)},\\
F_{l, i}(t, \bar f, \pa_t \bar f)&= ( \sin \om_i)^2\big|\na_{\tau_t}\cN\D_{\G_t}\bar f\big|^2 \qquad \textrm{at}\quad p_i\ (i=l, r).
\end{split}
\eeq
%which will be used later in the iteration section.

 \bthm{Proposition}\label{prop:energy est}(Higher-order energy estimates)
{\it Let $\mN_a\in C^0([0, T]; H^{5.5}(\Gamma_{*}))\cap C^1([0, T]; H^{4}(\Gamma_{*}))$ and $d_i, d'_i\in \R$ ($i=l,r$) be given, and $g_1\in C\big([0,T];H^{4}(\Gamma_{*})\big)$, $g_{2,i}\in C^1\big([0,T]; H^1(\G_*)\big)$, $(f_0, f_1)\in H^{5.5}(\Gamma_{*})\times H^{4}(\Gamma_{*})$.
Assume moreover that the contact angles $\om_i$ satisfy
\[
\min_i\sin \om_i\ge c_0\qquad \hbox{for some }\ c_0>0.
\]

Then  we have the following energy estimates for the linear problem \eqref{equ:linear}:
\[
\begin{split}
\| f\|^2_{H^{5.5}(\G_*)}+\|\pa_t f\|^2_{H^{4}(\G_*)}\le& e^{Q_1 t}Q_2(0)\big(\|f_0\|^2_{H^{5.5}(\G_*)}+\|f_1\|^2_{H^{4}(\G_*)}\big)\\
&\qquad +e^{Q_1 t}\int^t_0Q_1\big(\|g_1\|^{2}_{H^{4}(\G_{*})}+\|D_{t*} g_{2,i}\|^2_{H^1(\G_{*})}\big)dt',
\end{split}
\]
and
\[
\big|d_i(t)\big|+\big|d'_i(t)\big|\le \big|d_i(0)\big|+(1+t)\,\big|d'_i(0)\big|+\int^t_0\int^{t'}_0\big|B_i(\tau)\big|d\tau\,dt'+\int^t_0\big|B_i(t')\big|dt',
\]
where  $Q_1$ is a polynomial of the norms $\|d_{\G_t}\|_{H^{5.5}(\G_*)}$, $\|\pa_td_{\G_t}\|_{H^{5.5}(\G_*)}$, $\|\pa^2_t d_{\G_t}\|_{H^{4}(\G_*)}$,  $\|v\|_{H^6(\Om_t)}$, $\|D_t v\|_{H^{4.5}(\Om_t)}$, $Q_2(t)$ is a polynomial of $\|d_{\G_t}(t)\|_{H^{5.5}(\G_*)}$, $\|\pa_td_{\G_t}(t)\|_{H^{4}(\G_*)}$, $\|v(t)\|_{H^{4.5}(\Om_t)}$.
}
\ethm

\medskip

To prove this energy estimate, we firstly prove the estimate in forms of $E_h(t,\bar f,\pa_t \bar f)$ and $F_{h, i}(t,\bar f,\pa_t \bar f)$.

\bthm{Lemma}\label{energy}{\it Under the assumptions of Proposition \ref{prop:energy est}, we have for system \eqref{equ:linear on Gt} that
\[
\begin{split}
\pa_t E_h(t,\bar f,\pa_t \bar f) +\sum_iF_{h,i}(t,\bar f,\pa_t \bar f) \leq
 Q_1\,\Big(E_h(t,\bar f,\pa_t \bar f) +\|\bar f\|^2_{H^{5.5}(\G_t)}+\|D_t\bar f\|^2_{H^{4}(\G_t)}+\|\bar g_1\|^2_{H^{4}(\G_t)}+\|D_t\bar g_{2,i}\|^2_{H^1(\G_t)}\Big).
\end{split}
\]
%where $Q_1$ is is a polynomial of the norms $\|d_{\G_t}\|_{H^{5.5}(\G_*)}$, $\|\pa_td_{\G_t}\|_{H^{5.5}(\G_*)}$, $\|\pa^2_t d_{\G_t}\|_{H^{4}(\G_*)}$, $\|v\|_{H^6(\Om_t)}$, $\|D_t v\|_{H^{4.5}(\Om_t)}$.
}
\ethm

\begin{proof} (Proof for  Lemma \ref{energy})
Applying $\mathcal N \Delta_{\Gamma_t}$ on  both sides of the equation of $\bar f$ in \eqref{equ:linear on Gt} to obtain
\[
\cN\D_{\G_t}D^2_t\bar f-\si (\cN\D_{\G_t})^2\bar f=\cN\D_{\G_t}\bar g_1\qquad \hbox{on}\quad \G_t.
\]
Multiplying $-\Delta_{\Gamma_t}\cN\D_{\G_t} D_t\bar f$ on  both sides of this equation  and integrating on $\Gamma_t$, one obtains
\beq\label{integral eqn}
-\int_{\Gamma_t} \mathcal{N} \Delta_{\Gamma_t } D^2_t \bar f\cdot \Delta_{\Gamma_t}\cN\D_{\G_t} D_t\bar fds
 +\sigma\int_{\Gamma_t} (\mathcal{N} \Delta_{\Gamma_t})^2\bar f
\cdot \Delta_{\Gamma_t}\cN\D_{\G_t} D_t\bar fds
=-\int_{\Gamma_t} \mathcal{N} \Delta_{\Gamma_t} \bar g_1\cdot  \Delta_{\Gamma_t}\cN\D_{\G_t} D_t\bar fds.
\eeq
Next, we deal with the above integrals one by one.

\noindent{\bf Step 1. The first term on the left side.}
Firstly, one rewrites the first term in \eqref{integral eqn} as
\[
\begin{split}
-\int_{\Gamma_t} \mathcal{N} \Delta_{\Gamma_t } D^2_t \bar f\cdot \Delta_{\Gamma_t}\cN\D_{\G_t} D_t\bar fds
&=\int_{\Gamma_t} \na_{\tau_t}\cN\D_{\G_t} D^2_t \bar f\cdot \na_{\tau_t}\mathcal{N} \Delta_{\Gamma_{t}}D_t\bar fds-\sum_i\cN\D_{\G_t}D^2_t\bar f\cdot \na_{\tau_t}\mathcal{N} \Delta_{\Gamma_{t}}D_t\bar f\big|_{p_i} \\
&\triangleq  \f12 \pa_t\big\|\na_{\tau_t}\cN\D_{\G_t}D_t\bar f\big\|^2_{L^2(\G_t)}+A_1+A_2,
\end{split}
\]
where
\[
\begin{split}
&A_1=\int_{\Gamma_t} [\na_{\tau_t}\cN\D_{\G_t}, D_t] D_t \bar f\cdot \na_{\tau_t}\mathcal{N} \Delta_{\Gamma_{t}}D_t\bar fds-\f12\int_{\G_t}(\cD\cdot v)\big|\na_{\tau_t}\cN\D_{\G_t}D_t\bar f\big|^2ds,\\
&A_2= -\cN\D_{\G_t}D^2_t\bar f\cdot \na_{\tau_t}\mathcal{N} \Delta_{\Gamma_{t}}D_t\bar f\big|^{p_l}_{p_r}.
\end{split}
\]

%Secondly, we will deal with  $A_1, A_2$ separately.

\noindent{\bf Term $A_1$}.  In fact,  a direct computation leads to
\[
[\na_{\tau_t}\cN\D_{\G_t}, D_t] D_t \bar f=[\na_{\tau_t}, D_t]\cN\D_{\G_t}D_t \bar f
+\na_{\tau_t}[\cN, D_t]\D_{\G_t}D_t \bar f+\na_{\tau_t}\cN[\D_{\G_t}, D_t]D_t \bar f.
\]
Applying \eqref{commutator DN},  \eqref{commutator surface delta}, Lemma \ref{est:elliptic} and Lemma \ref{Dt v estimate} (1),  we have
\beq\label{A1}
|A_1|\le C\big(\|d_{\G_t}\|_{H^{4}(\G_*)}, \|\pa_t d_{\G_t}\|_{H^2(\G_*)},
\|v\|_{H^{4.5}(\Om_t)}\big)\big(\|D_t\bar f\|^2_{H^4(\G_t)} +E_h(t,\bar f, \pa_t \bar f)\big).
\eeq

\noindent{\bf Term $A_2$}. Recalling the condition at the corner points from system \eqref{equ:linear on Gt} and applying  $D_t$ on the both sides again, one derives
\[
\cN\D_{\G_t}D^2_t\bar f\pm\f{\si^2}{\be_c}(\sin\om_i)^2\na_{\tau_t}\cN\D_{\G_t}D_t \bar f=R_{c2, i}\qquad\hbox{at}\quad p_i\ (i=l,r),
\]
with
\[
R_{c2, i}=D_t \bar g_{2,i}-[D^2_t, \cN\D_{\G_t}]\bar f-\f{\si^2}{\be_c}D_t(\sin\om_i)^2\na_{\tau_t}\cN\D_{\G_t} \bar f
-\f{\si^2}{\be_c}(\sin\om_i)^2[D_t, \na_{\tau_t}\cN\D_{\G_t}] \bar f.
\]

Substituting this equality into $A_2$, one obtains
\beq\label{A2}
A_2=\f{\si^2}{\be_c}\sum_iF_{h, i}(t, \bar f, \pa_t\bar f)+A_{2R},
\eeq
with
\[
A_{2R}=-R_{c2,i}\cdot \na_{\tau_t}\cN\D_{\G_t}D_t\bar f\big|^{p_l}_{p_r}.
\]
Moreover, analogous analysis as before shows that the remainder term $A_{2R}$ satisfies the estimate
\[
\begin{split}
|A_{2R}|&\le \del F_h(t, \bar f, \pa_t\bar f)+(\sin\om_i)^{-2} C_\del \|R_{c2}\|^2_{L^\infty(\G_t)}\\
&\le \del F_h(t, \bar f, \pa_t\bar f)+(\sin\om_i)^{-2} C_\del \,C\big(\|d_{\G_t}\|_{H^5(\G_*)}, \|\pa_t d_{\G_t}\|_{H^3(\G_*)}, \|v\|_{H^{5.5}(\Om_t)}, \|D_t v\|_{H^{4.5}(\Om_t)}\big)\times\\
&\qquad\quad\big(\|\bar f\|^2_{H^5(\G_t)}+\|D_t\bar f\|^2_{H^4(\G_t)}+\|D_t\bar g_{2,i}\|^2_{H^1(\G_t)}\big),
\end{split}
\]
with $C_\del$ a constant depending on $\del^{-1}, \be_c, \si$.

\medskip

\noindent{\bf Step 2. The second term on the left side.} We firstly rewrite this integral by Green's formula as follows:
\[
\begin{split}
\sigma\int_{\Gamma_t} (\mathcal{N} \Delta_{\Gamma_t})^2\bar f
\cdot \Delta_{\Gamma_t}\cN\D_{\G_t} D_t\bar f\,ds
&=\si\int_{\Om_{t}}\na\cH(\D_{\G_t}\cN\D_{\G_t}\bar f)\cdot \na\cH(\Delta_{\Gamma_t}\cN\D_{\G_t} D_t\bar f)dX\\
&\triangleq \f\si 2  \pa_t \big\|\na\cH(\D_{\G_t}\cN\D_{\G_t}\bar f)\big\|^2_{L^2(\Om_t)}+A_3+A_4,
\end{split}
\]
where the remainder terms
\[
\begin{split}
A_3&=\si\int_{\Gamma_{t}}\na\cH(\D_{\G_t}\cN\D_{\G_t}\bar f)\cdot [\na\cH\Delta_{\Gamma_t}\cN\D_{\G_t}, D_t]\bar f dX\\
&=\si \int_{\Gamma_{t}}\na\cH(\D_{\G_t}\cN\D_{\G_t}\bar f)\cdot [\na\cH, D_t]\Delta_{\Gamma_t}\cN\D_{\G_t} \bar fdX+\si\int_{\Gamma_{t}}\na\cH(\D_{\G_t}\cN\D_{\G_t}\bar f)\cdot \na\cH[\Delta_{\Gamma_t}\cN\D_{\G_t}, D_t]\bar f dX\\
&\triangleq A_{31}+A_{32}
\end{split}
\]
and 
\[
A_4=-\f\si 2\int_{\Om_t}(\dive v)\big|\na\cH(\D_{\G_t}\cN\D_{\G_t}\bar f)\big|^2dX.
\]

\noindent{\bf Term $A_{31}$.} To handle $A_{31}$,  denoting by
\[
g=\Delta_{\Gamma_{t}}\cN\D_{\G_t} \bar f,
\]
 a direct computation leads to
\[
A_{31}=\int_{\Om_t}\na g_{\cH}\cdot [ \na, D_t] g_{\cH}dX+ \int_{\Om_t}\na g_{\cH}\cdot \na[\cH, D_t]  g dX.
\]
Noticing that
\[
[\na, D_t]g_{\cH}=-\na v\cdot \na g_{\cH},
\]
and due to the definition of $g_\cH$ and
\beno
[\mathcal{H} , D_t ] g=\Delta^{-1} \big(2\na v \cdot \na^2 g_{\cH}+\Delta v  \cdot \na g_{\cH}, \,(\na_{N_b }v-\na_{v }N_b )\cdot \na g_{\cH}\big),
\eeno
we know
\[
\begin{split}
\int_{\Om_t}\na g_{\cH}\cdot \na[\cH, D_t]  g dX&=- \int_{\Om_t} \Delta g_{\cH} \cdot[\cH, D_t]  gdX+\int_{\Gamma_t} \na_{N_t} g_{\cH}\cdot[\cH, D_t]  gds+\int_{\Gamma_b}\na_{N_b} g_{\cH}\cdot[\cH, D_t]  gds\\
 &=0.
\end{split}\]

Consequently,  one can show that
\beq\label{A31}
|A_{31}|\le C\big(\|d_{\G_t}\|_{H^{2.5}(\G_*)}, \|v\|_{H^{2.5}(\Om_t)}\big)E_h(t, \bar f, \pa_t\bar f).
\eeq

\medskip

\noindent{\bf Term $A_{32}$ and $A_4$.} Applying Lemma \ref{Harmonic extension H1 estimate}, Lemma \ref{est:elliptic} and recalling  the commutators \eqref{commutator DN}, \eqref{commutator surface delta},  one obtains
\beq\label{A32}
\begin{split}
|A_{32}|+|A_4|&\le C(\|d_{\G_t}\|_{H^{2.5}(\G_t)})\|[\Delta_{\Gamma_t}\cN\D_{\G_t}, D_t]\bar f\|_{H^{0.5}(\G_t)}E_h(t,\bar f, \pa_t\bar f)^{1/2}\\
&\le C\big(\|d_{\G_t}|_{H^{5.5}(\G_t)}, \|v\|_{H^6(\Om_t)}\big)\big(\|\bar f\|^2_{H^{5.5}(\G_t)}+E_h(t,\bar f, \pa_t\bar f)\big).
\end{split}
\eeq

\medskip

As a result,  combining \eqref{A1}, \eqref{A2}, \eqref{A31} and \eqref{A32}, we finally conclude that the left side of equation \eqref{integral eqn} becomes
\[
\mbox{Left side}=\f12\pa_t E_h(t,\bar f,\pa_t \bar f)+\f{\sigma^2}{\beta_c}F_h(t, \bar f,\pa_t \bar f)+R_L,
\]
where the remainder term
\[
R_L=A_1+A_{2R}+A_{31}+A_{32}+A_4,
\]
satisfying the estimate
\[
\begin{split}
|R_L|
&\le \del F_h(t, \bar f, \pa_t\bar f)+(\min_i \{\sin\om_i\})^{-2} C_\del \,C\big(\|d_{\G_t}\|_{H^{5.5}(\G_*)}, \|\pa_t d_{\G_t}\|_{H^3(\G_*)}, \|v\|_{H^{6}(\Om_t)}, \|D_t v\|_{H^{4.5}(\Om_t)}\big)\times\\
&\qquad\quad\big(E_h(t,\bar f,\pa_t \bar f)+\|\bar f\|^2_{H^{5.5}(\G_t)}+\|D_t\bar f\|^2_{H^4(\G_t)}+\|D_t\bar g_{2,i}\|^2_{H^1(\G_t)}\big).
\end{split}
\]

\medskip

\noindent{\bf Step 3. The right side.} Integrating by parts on the right side leads to
\[
\mbox{Right side}
=\int_{\G_t} \na_{\tau_t}\cN\D_{\G_t}\bar g_1\cdot \na_{\tau_t}\cN\D_{\G_t}D_t\bar f\,ds
-\cN\D_{\G_t}\bar g_1\cdot \na_{\tau_t}\cN\D_{\G_t}D_t\bar f\big|^{p_l}_{p_r}.
\]
Therefore, we have
\[
\begin{split}
\mbox{Right side}&\le \del F_h(t, \bar f, \pa_t\bar f)+(\min_i \{\sin\om_i\})^{-2}C_\del \|\cN\D_{\G_t}\bar g_1\|^2_{H^1(\G_t)}+\|\na_{\tau_t}\cN\D_{\G_t}\bar g_1\|_{L^2(\G_t)}E_h(t,\bar f,\pa_t \bar f)^{1/2}\\
&\le \del F_h(t, \bar f, \pa_t\bar f)+(\min_i \{\sin\om_i\})^{-2}C_\del \,C(\|d_{\G_t}\|_{H^4(\G_*)})\Big(E_h(t,\bar f,\pa_t \bar f)+\|\bar g_1\|^2_{H^4(\G_t)}\Big).
\end{split}
\]

\noindent{\bf Step 4. $L^2$ energy.}
Multiplying $D_t \bar f$ on the both sides of \eqref{equ:linear on Gt} and integrating on $\G_t$ to get
\[
\begin{split}
\pa_t\|D_{t} \bar f\|^2_{L^2(\Gamma_t)}&\leq C(\|\cN\D_{\G_t}\bar f\|_{L^2(\Gamma_t)}+\|\bar g_1\|_{L^2(\Gamma_t)})\|D_{t} \bar f\|_{L^2(\Gamma_t)}+C(\|d_{\G_t}\|_{H^{2.5}(\G_*)}, \|v\|_{H^{2.5}(\Om_t)})\|D_t\bar f\|^2_{L^2(\G_t)}\\
&\le C(\|d_{\G_t}\|_{H^3(\G_*)}, \|v\|_{H^{2.5}(\Om_t)})\Big(\|\bar f\|^2_{H^3(\G_t)}+\|\bar g_1\|^2_{L^2(\G_t)}+E_h(t,\bar f,\pa_t \bar f)\Big).
\end{split}
\]
Moreover,  using 
\[
\pa_t\|  \bar f\|^2_{L^2(\Gamma_t)}=2 \int_{\Gamma_t} D_t\bar f\cdot \bar f ds+\int_{\G_t}(\cD\cdot v)|\bar f|^2ds,
\]
 we have immediately that
\beno
 \pa_t\|  \bar f\|^2_{L^2(\Gamma_t)}
 \le C(\|d_{\G_t}\|_{H^{2.5}(\G_*)}, \|v\|_{H^{2.5}(\Om_t)})\,E_h(t, \bar f,\pa_t \bar f).
\eeno

\medskip

Combing all the  estimates from Step 1 to Step 4, we get the desired estimate.
\end{proof}

Moreover, we replace the energy $E_h(t,\bar f,\pa_t \bar f)$ with $\|\bar f\|_{H^{5.5}(\G_t)}$ and $\|D_t\bar f\|_{H^{4}(\G_t)}$, which would be more convenient to use.

\bthm{Lemma}\label{E(t) equivalence}{\it
\noindent (1) There exists  a polynomial $Q_l$ of $\|d_{\G_t}\|_{H^4(\G_*)}$, $\|v\|_{H^3(\Om_t)}$ such that
\[
E_l(t,\bar f,\pa_t \bar f)\le Q_l(\|d_{\G_t}\|_{H^4(\G_*)},\|v\|_{H^3(\Om_t)}) \big(\|D_t\bar f\|^2_{H^{2.5}(\G_t)}+\|\bar f\|^2_{H^4(\G_t)}\big),
\] and
\[
\|D_t\bar f\|^2_{H^{2.5}(\G_t)}+\|\bar f\|^2_{H^4(\G_t)}\le Q_l(\|d_{\G_t}\|_{H^4},\|v\|_{H^3(\Om_t)})\,E_l(t,\bar f,\pa_t \bar f);
\]
\noindent (2) There exists   a polynomial $Q_h$  of $\|d_{\G_t}\|_{H^5(\G_*)}$, $\|v\|_{H^{4.5}(\Om_t)}$ such that
\[
E_h(t,\bar f,\pa_t \bar f)\le Q_h(\|d_{\G_t}\|_{H^5(\G_*)},\|v\|_{H^{4.5}(\Om_t)}) \big(\|D_t\bar f\|^2_{H^{4}(\G_t)}+\|\bar f\|^2_{H^{5.5}(\G_t)}\big),
\] and
\[
\|D_t\bar f\|^2_{H^{4}(\G_t)}+\|\bar f\|^2_{H^{5.5}(\G_t)}\le Q_h(\|d_{\G_t}\|_{H^5(\G_*)},\|v\|_{H^{4.5}(\Om_t)})\,E_h(t,\bar f,\pa_t \bar f).
\]
}\ethm
\begin{proof}
In fact,  the first inequality in (1) can be  proved directly by applying  Lemma \ref{Harmonic extension H1 estimate}, \eqref{commutator surface delta}, Lemma \ref{trace thm PG} and Lemma \ref{est:elliptic}.

We focus on  the second inequality in (1).  Due to the definition of the norm $H^{2.5}(\G_t)$, one has
\[
\|D_t\bar f\|_{H^{2.5}(\G_t)}\le C(\|d_{\G_t}\|_{H^{2.5}(\G_*)})\big(\|\D_{\G_t}D_t\bar f\|_{H^{0.5}(\G_t)}+\|D_t\bar f\|_{L^2(\G_t)}\big).
\]
Besides, one has  by Lemma \ref{trace thm PG} that
\[
\begin{split}
\|\D_{\G_t}D_t\bar f\|_{H^{0.5}(\G_t)}&\le \|D_t\D_{\G_t}\bar f\|_{H^{0.5}(\G_t)}+\|[\D_{\G_t}, D_t]\bar f\|_{H^{0.5}(\G_t)}\\
&\le C(\|d_{\G_t}\|_{H^{2.5}(\G_*)})\big(\|\cH(D_t\D_{\G_t}\bar f)\|_{H^1(\Om_t)}+\|v\|_{H^3(\Om_t)}\|\bar f\|_{H^{2.5}(\G_t)}\big).
\end{split}
\]

Noticing that
\[
\|\cH(D_t\D_{\G_t}\bar f)\|_{L^2(\Om_t)}\le C\big(\|\na\cH(D_t\D_{\G_t}\bar f)\|_{L^2(\Om_t)}+\|D_t\D_{\G_t}\bar f\|_{L^2(\G_t)}\big),
\]
 and combining all the inequalities  above one can derive the following estimate
\beq\label{Dt f H2.5 estimate}
\|D_t\bar f\|_{H^{2.5}(\G_t)}\le C(\|d_{\G_t}\|_{H^{2.5}(\G_*)}, \|v\|_{H^3(\Om_t)})\big(E(t,\bar f,\pa_t \bar f)^{1/2}+\eps\,\|\bar f\|_{H^4(\G_t)}\big),
\eeq
where one uses an interpolation for $\|\bar f\|_{H^{2.5}(\G_t)}$ and $\eps$ is a small constant.

Secondly, one writes  that
\[
\|\bar f\|_{H^4(\G_t)}\le C(\|d_{\G_t}\|_{H^{4}(\G_*)})\big(\|\bar f\|_{L^2(\G_t)}+\|\D_{\G_t}\bar f\|_{H^2(\G_t)}\big),
\]
where the term $\|\D_{\G_t}\bar f\|_{H^2(\G_t)}$ needs to be handled.
Applying Lemma \ref{DN-1 operator} (2),  we obtain immediately
\[
\|\D_{\G_t}\bar f\|_{H^2(\G_t)}\le C(\|d_{\G_t}\|_{H^{2.5}(\G_*)})\big(\|\cN\D_{\G_t}\bar f\|_{H^1(\G_t)}+\|\D_{\G_t}\bar f\|_{L^2(\G_t)}\big),
\]
which leads to the following inequality immediately:
\beq\label{f H4 estimate}
\begin{split}
\|\bar f\|_{H^4(\G_t)}&\le C(\|d_{\G_t}\|_{H^{4}(\G_*)})\big(\|\bar f\|_{L^2(\G_t)}+\|\cN\D_{\G_t}\bar f\|_{H^1(\G_t)}\big)\\
&\le  C(\|d_{\G_t}\|_{H^{4}(\G_*)}) E_l(t, \bar f,\pa_t \bar f)^{1/2}.
\end{split}
\eeq
In the end, one only need to combine \eqref{Dt f H2.5 estimate} with \eqref{f H4 estimate} and choose $\eps$ small enough to finish the proof for (1). Moreover, the proof for (2) follows in the same way.
\end{proof}

\medskip

Now, we are in the position to prove the Proposition \ref{prop:energy est}.
\begin{proof}(Proof for Proposition \ref{prop:energy est})
In fact,  multiplying $Q_2$ defined in the Proposition \ref{prop:energy est} on both sides of the inequality in Lemma \ref{energy} and using Lemma \ref{E(t) equivalence} leads to
\[
Q_2\pa_t E_h(t,\bar f,\pa_t \bar f) +Q_2\sum_i F_{h,i}(t,\bar f,\pa_t \bar f)\leq
 Q_1Q_2\,\Big(E_h(t,\bar f,\pa_t \bar f) +\|\bar g_1\|^2_{H^{4}(\G_t)}+\|D_t\bar g_{2,i}\|^2_{H^1(\G_t)}\Big).
\]
On the other hand, acting $\pa_t$ on $Q_2$,  one has
\[
(\pa_t Q_2)\,E_h(t,\bar f,\pa_t \bar f)\le Q_1 E_h(t,\bar f,\pa_t \bar f)\le Q_1 Q_2E_h(t,\bar f,\pa_t \bar f),
\]
where $Q_1$ is defined in the Proposition \ref{prop:energy est}. Summing up these two inequalities and applying Lemma \ref{E(t) equivalence} again, one one derives
\[
\begin{split}
&Q_2 E_h(t, \bar f, \pa_t\bar f)+e^{Q_1t}\int^t_0 e^{-Q_1t'}Q_2\sum_i F_{h,i}(t',\bar f,\pa_{t} \bar f)dt'\\
&\le
e^{Q_1 t}Q_2\big|_{t=0}E_h(t, \bar f, \pa_t\bar f)\big|_{t=0}+e^{Q_1t}\int^t_0 e^{-Q_1t'}Q_1 \big(\|\bar g_1\|^2_{H^{4}(\G_{t'})}+\|D_t\bar g_{2,i}\|^2_{H^1(\G_{t'})}\big)dt'.
\end{split}
\]

On the other hand,   one knows that
\[
\begin{split}
\|f\|^2_{H^{5.5}(\G_*)}+\|\pa_t f\|^2_{H^4(\G_*)}&\le C\big(\|d_{\G_t}\|_{H^{5.5}(\G_*)}, \|\pa_td_{\G_t}\|_{H^4(\G_*)}, \|v\|_{H^{4.5}(\Om_t)}\big)\big(\|f\|^2_{H^{5.5}(\G_*)}+\|D_{t*} f\|^2_{H^4(\G_*)}\big)\\
&\le Q_2\big(\|\bar f\|^2_{H^{5.5}(\G_t)}+\|D_{t} \bar f\|^2_{H^4(\G_t)}\big).
\end{split}
\]
Combining this with Lemma \ref{E(t) equivalence}, we derive immediately the desired energy estimate for $f$. Moreover, the estimate for $d(t)$ follows immediately by integrating with respect to time twice on its equation.

\end{proof}

\medskip

In the end, the estimate for the second-order time derivatives of $f$ is  considered here, which is based on  the linear system \eqref{equ:linear}.
\bthm{Corollary}\label{pa 2 t f estimate}{\it Under the assumptions of Proposition \ref{prop:energy est}, one has the following estimate for $f$ in system \eqref{equ:linear} when $s\ge 3$:
\[
\begin{split}
\|\pa^2_tf\|_{H^s(\G_*)}\le &\|g_1\|_{H^s(\G_*)}+C\big(\|\mN_a\|_{H^{s-1}(\G_*)}, \|\pa_t \mN_a\|_{H^{s-2}(\G_*)},  a^{-1}\|\pa^2_t \mN_a\|_{H^{s-2}(\G_*)}, |d_i|, |d'_i|, a^{-1}|d''_i|\big)\times\\
&\qquad\big(\|f\|_{H^{s+3}(\G_*)}+\|\pa_t f\|_{H^{s+1}(\G_*)}\big),
\end{split}
\]
where the index $i$ stands for the summation on $i=l,r$.
When $1\le s\le 3$, the norms for $\mN_a$, $\pa_t\mN_a, \pa^2_t\mN_a$ on the right side should be higher than $H^1(\G_*)$.
}
\ethm
\begin{proof} To begin with, one knows from  the equation for $f$ in system \eqref{equ:linear} that
\[
\pa^2_t f=-\na_{\pa_tv^*}f-2\na_{v^*}\pa_t f-\na^2_{v^*}f+\si\cA(d_{\G_t})f+g_1,
\]
where $v^*$ is defined by \eqref{v* expression} and $\pa_t v^*$ is expressed by \eqref{pa t v*}.

One can have immediately 
\[
\|\pa^2_tf\|_{H^s(\G_*)}\le C\big(\|d_{\G_t}\|_{H^{s+2}(\G_*)}, \|\pa_tv^*\|_{H^s(\G_*)}, \|v^*\|_{H^{s+1}(\G_*)}\big)\big(\|f\|_{H^{s+3}(\G_*)}+\|\pa_t f\|_{H^{s+1}(\G_*)}\big)+\|g_1\|_{H^s(\G_*)}.
\]
Moreover,  \eqref{pa t v*} leads to 
\[
\|\pa_tv^*\|_{H^s(\G_*)}\le C\big(\|d_{\G_t}\|_{H^{s+1}(\G_*)}, \|\pa_t d_{\G_t}\|_{H^{s+1}(\G_*)}, \|\pa^2_t d_{\G_t}\|_{H^s(\G_*)}, \|v\|_{H^{s+1.5}(\Om_t)}, \|D_t v\|_{H^{s+0.5}(\Om_t)}\big)
\]
and
\[
\|v^*\|_{H^{s+1}(\G_*)}\le C\big(\|d_{\G_t}\|_{H^{s+2}(\G_*)}, \|\pa_t d_{\G_t}\|_{H^{s+1}(\G_*)}, \|v\|_{H^{s+1.5}(\Om_t)}\big),
\]
so one obtains
\[
\begin{split}
\|\pa^2_tf\|_{H^s(\G_*)}\le &\|g_1\|_{H^s(\G_*)}+C\big(\|d_{\G_t}\|_{H^{s+2}(\G_*)}, \|\pa_t d_{\G_t}\|_{H^{s+1}(\G_*)},  \|\pa^2_t d_{\G_t}\|_{H^s(\G_*)}, \|v\|_{H^{s+1.5}(\Om_t)}\big)\times\\
&\qquad \big(\|f\|_{H^{s+3}(\G_*)}+\|\pa_t f\|_{H^{s+1}(\G_*)}\big).
\end{split}
\]
As a result,  the desired estimate follows from Proposition \ref{d and  N ka} and Proposition \ref{pa t d estimate}.
\end{proof}

\section{Local well-posedness}
In this section, we use a standard iteration method and a fixed point theorem to prove that the system of \eqref{ka_a equation}, \eqref{d c eqn} admits a unique solution $(\mN_a, d_l, d_r)$ satisfying
\beno
\mN_a \in C^0\big([0, T]; H^{5.5}(\Gamma_*)\big)\cap C^1\big([0,T]; H^{4}(\G_*)\big)
\quad\hbox{and}\quad  \big|d_i\big|\le L_0,\ \big|d'_i\big|\le L_1\ (i=l, r),
\eeno
where the constants $L_0, L_1$ are specified later.
Consequently, we can construct the free surface $\Phi_{S_t}$ by $d_{\G_t}\in C([0,T]; H^{8.5}(\G_*))\cap C^1([0,T]; H^{7}(\G_*))$ from Proposition \ref{d and  N ka}, Proposition \ref{pa t d estimate}, and we define the velocity $v\in H^{7.5}(\Om_t)$ by \eqref{v potential expression}. In the end,  we prove that the system of $(\mN_a, d_l, d_r)$ is equal to the  system  $\mbox{(WW)}$ to finish the proof for the local well-posedness of the water-waves problem.

\subsection{Settings}
Before we start the iteration scheme, we need to clarify some settings.
First of all, we will use the set of surfaces $\Lam_*= \Lam(S_*,8.5,\del, \pi/16)$ from Definition \ref{Lambda neighborhood}.
The following  set $\Sigma$ is defined for the bounds related to $\mN_a$ and $d_l, d_r$.
\bthm{Definition}{\it
The set $\Sigma$ is defined as the collection of $(\mN_a, d_l, d_r)$ which satisfies the conditions for the initial data:
\[
\|\mN_a(0,\cdot)-\mN_{a*}\|_{H^{5.5}(\Gamma_{*})}, |d_i(0)|\leq \delta_1,
\quad\|\pa_t\mN_{a}(0,\cdot)\|_{H^{4}(\Gamma_*)}, |d'_i(0)|\leq L\ (i=l, r),
\]
as well as the higher-order bounds  for $\mN_a$ and  $d_i$:
\[
\begin{split}
&\|\mN_a\|_{C([0, T]; H^{5.5}(\Gamma_*))}, |d_i|_{C([0,T])}\leq L_0,\quad  \|\pa_t\mN_a\|_{C([0, T]; H^{4}(\Gamma_*))}, |d'_i|_{C([0,T])}\leq L_1,\\
& \|\pa^2_t\mN_a\|_{C([0, T]; H^{2.5}(\Gamma_*))}, |d''_i|_{C([0,T])}\leq L_2,
\end{split}
\]
where the constants $T, \del_1, L, L_0, L_1, L_2$ are given with small $\del_1$.
}
\ethm

Meanwhile, we introduce the set of initial data
\[
\begin{split}
\mathcal{I}(\eps, A_1)\triangleq &\big\{ \,\big((\mN_a)_I, (\pa_t\mN_a)_I, (d_i)_I, (d'_i)_I\big)\,\big|  \,
\big((\mN_a)_I, (\pa_t\mN_a)_I\big)\in H^{5.5}(\Gamma_*)\times H^{4}(\Gamma_*), (d_i), (d'_i)_I\in \R, \\
&\qquad \|(\mN_a)_I-\mN_{a*}\|_{H^{5.5}(\Gamma_*)},  |(d_i) _I|\le \eps, \big\|(\pa_t\mN_a)_{I}\big\|_{H^{4}(\Gamma_*)}, \big|(d'_i)_I\big|\le A_1,  \, (\om_i)_{I}\in (0, \pi/16), i=l, r\big\}
\end{split}
\]
where $0<\eps<<\delta_1$ and $A_1>0$ some large constant, and $(\om_i)_{I}$ are the corresponding  contact angles. When $\eps$ is taken sufficiently small,  $(\om_i)_{I}\in (0, \pi/16)$ will be satisfied naturally.

\subsection{The iteration scheme} Fixing the initial data $\big((\mN_a)_I,(\pa_t\mN_a)_I, (d_i)_I, (d'_i)_I \big)\in \mathcal{I}(\eps, A_1)$, we plan to use the linear system \eqref{equ:linear}  to generate an iteration sequence of $(\mN_a, d_l, d_r)$.  Then, we will show that this sequence converges to $(\mN_a, d_l, d_r)$ by a fixed point theorem.

To get started, for any  given $(\mN_a, d_l, d_r)\in \Sigma$, using Proposition \ref{d and  N ka}, we define the free surface $\G_t$ (i.e. $\Phi_{S_t}$) by the function $d_{\G_t}\in C([0,T];H^{8.5}(\G_*))$ with the condition at  the corner point $p_{i*}$
\[
d_{\G_t}(p_{i*})=d_i, \ i=l, r
\]
and the velocity $v$ is defined by \eqref{v potential expression}.

Let  $(\widetilde{\mN_a}, \widetilde d_l, \widetilde d_r)$ be the solution of the following linear system
 \ben\label{equ:map1}
\left\{\begin{array}{ll}
D^2_{t*} \widetilde{\mN_a}+\sigma \cA (d_{\G_t})\widetilde{\mN_a}=R_0 \qquad\hbox{on}\quad \G_*,\\
D_{t*}\cA(d_{\G_t})\widetilde{\mN_a}\pm\f{\sigma^2 }{\beta_c}(\sin \om_i)^2(\na_{\tau_t}(\cA(d_{\G_t})\widetilde{\mN_a}\circ \Phi_{S_t}^{-1}))\circ\Phi_{S_t}= R_{c, i}\qquad \textrm{at} \quad p_{i*},\\
 \widetilde d''_i(t)=B_i \quad (i=l, r)
\end{array}
\right.
\een
with initial data
\[
\widetilde{\mN_a}|_{t=0}=(\mN_a)_I, \  \pa_t  \widetilde{\mN_a}|_{t=0}=(\pa_t\mN_a)_I,\  \widetilde d_i(0)=(d_i)_I, \ \widetilde d'_i(0)=(d'_i)_I.
\]
Here the right side $R_0, R_{c, i}$ and $B_i$ are defined in \eqref{ka_a equation}, \eqref{equ:cca} and \eqref{d c eqn} , which depend on the known functions $(\mN_a, d_l, d_r)\in \Sigma$.

By Remark \ref{rmk: com cond}, it is easy to get that the initial data satisfy the following  compatibility condition
\beno
D_{t*}\cA(d_{\G_t})\widetilde{\mN_a}|_{t=0}\pm\f{\sigma^2 }{\beta_c}(\sin \om_i)^2(\na_{\tau_t}(\cA(d_{\G_t})\widetilde{\mN_a}\circ \Phi_{S_t}^{-1}))\circ\Phi_{S_t}|_{t=0}= R_{c, i}|_{t=0}\quad \textrm{at} \quad p_{i*}.
 \eeno

Applying Proposition \ref{existence of linear system} and Proposition \ref{prop:energy est}, we can prove immediately that  \eqref{equ:map1} admits a unique solution $(\widetilde{\mN_a}, \widetilde d_l, \widetilde d_r)$ satisfying the initial data  above.

Moreover, the following higher-order energy estimates hold for $\forall t\in [0,T]$:
\beq\label{energy estimate for N ka}
\begin{split}
\big\|\widetilde{\mN_a}\big\|^2_{H^{5.5}(\Gamma_*)}+\big\|  \pa_t \widetilde{\mN_a}\big\|^2_{H^{4}(\Gamma_*)}
\le& e^{Q_1 t}Q_2(0)\big(\big\|(\mN_a)_I\big\|^2_{H^{5.5}(\G_*)}+\big\|(\pa_t\mN_a)_I\big\|^2_{H^{4}(\G_*)}\big)\\
&\qquad +e^{Q_1 t}\int^t_0Q_1\big(\|R_0\|^{2}_{H^{4}(\G_{t'})}+|D_{t*} R_{c, i}|^2_{H^1(\G_{t'})}\big)dt',
\end{split}
\eeq
and
\[
\big|\widetilde d_i\big|+\big|\widetilde d_i'\big|\le \big|(d_i)_I\big|+(1+t)\,\big|(d'_i)_I\big|+\int^t_0\int^{t'}_0\big|B_i(\tau)\big|d\tau\,dt'+\int^t_0\big|B_i(t')\big|dt',\quad(i=l, r),
\]
where $Q_1$ is a polynomial of the norms $\|d_{\G_t}\|_{H^{5.5}(\G_*)}$, $\|\pa_td_{\G_t}\|_{H^{5.5}(\G_*)}$, $\|\pa^2_t d_{\G_t}\|_{H^{4}(\G_*)}$,  $\|v\|_{H^6(\Om_t)}$, $\|D_t v\|_{H^{4.5}(\Om_t)}$, $Q_2(0)$ is a polynomial of $\|d_{\G_t}(0)\|_{H^{5.5}(\G_*)}$, $\|\pa_td_{\G_t}(0)\|_{H^{4}(\G_*)}$, $\|v(0)\|_{H^{4.5}(\Om_t)}$.

Applying Lemma \ref{Dt v estimate} (2), Proposition \ref{d and  N ka} and Proposition \ref{pa t d estimate}, one has immediately that
\[
\begin{split}
Q_1&=Q_1\big(\|d_{\G_t}\|_{H^{5.5}(\G_*)}, \|\pa_td_{\G_t}\|_{H^{5.5}(\G_*)}, \|\pa^2_t d_{\G_t}\|_{H^{4}(\G_*)}, \|v\|_{H^6(\Om_t)}, \|D_t v\|_{H^{4.5}(\Om_t)}\big)\\
&\le C\big(\|d_{\G_t}\|_{H^7(\G_*)}, \|\pa_t d_{\G_t}\|_{H^{5.5}(\G_*)}, \|\pa^2_td_{\G_t}\|_{H^{4}(\G_*)} \big)\\
&\le C\big(\big\|\mN_a-\mN_{a*}\big\|_{H^{4}(\G_*)}, \big\|\pa_t\mN_a\big\|_{H^{2.5}(\G_*)}, \big\|\pa^2_t\mN_a\big\|_{H^{1}(\G_*)}, |d_i|, \big|d'_i\big|, \big|d''_i\big|\big)\\
&
\le C\big(L_0, L_1, L_2).
\end{split}\]

On the other hand,  the estimates for $\|R_0\|_{H^{4}(\G_*)}$ and $\|D_{t*} R_{c, i}\|_{H^1(\G_*)}$ are handled  in Lemma \ref{R0 estimate} and Lemma \ref{lem:R_c}:
\[
\begin{split}
&\big\|R_0\big\|_{H^{4}(\G_*)} +\big\|D_{t*}R_{c, i}\big\|_{H^1(\G_*)} \\
&\le a^{3/2} C\big(\big\|\mN_a\big\|_{H^{5.5}(\G_*)}, \big\|\pa_t\mN_a\big\|_{H^{4}(\G_*)}, \big\|\pa^2_t\mN_a\big\|_{H^{2.5}(\G_*)}, |d_i|, \big|d'_i\big|, \big|d''_i\big|\big)
\le a^{3/2} C\big(L_0, L_1, L_2).
\end{split}
\]
Moreover, we have from Lemma \ref{Dt v estimate} (2), Propositon \ref{d and  N ka} and  Proposition \ref{pa t d estimate} that
\[
Q_2(0)\le C(\del_1, L).
\]

Therefore, combining these estimates above  and going back to \eqref{energy estimate for N ka}, one obtains that
\[
\big\|\widetilde{\mN_a}\big\|^2_{H^{5}(\Gamma_*)}+\big\|  \pa_t \widetilde{\mN_a}\big\|^2_{H^{4}(\Gamma_*)}\le e^{C(L_0, L_1, L_2)T}\big(C(\del_1, L, A_1)+a^{3}T\,C(L_0, L_1,L_2) \big).
\]
Meanwhile, one also has from \eqref{d c eqn} the following estimate for $i=l ,r$:
\[
\begin{split}
\big|\widetilde d_i\big|+\big|\widetilde d_i'\big|&\le \big|(d_i)_I\big|+(1+T)\big|(d'_i)_I\big|+T\,C(a^3\|d_{\G_t}\|_{H^2(\G_*)}, \|\pa_td_{\G_t}\|_{H^2(\G_*)}, \|\ka_a\|_{H^2(\G_*)}, \|v\|_{H^1(\G_*)})\\
&\le C(\del_1, L, A_1)+T\,C(L_0, L_1).
\end{split}
\]

Choosing $T$ small enough compared to $a, L_0, L_1, L_2$,  and $L_0, L_1$ large enough compared to $L, A_1$ in the inequalities above, we have
\[
\big\|\widetilde{\mN_a}\big\|_{H^{5.5}(\Gamma_*)}, \big|\widetilde d_l\big|, \big|\widetilde d_r\big|\le L_0,\quad
\big\|  \pa_t \widetilde{\mN_a}\big\|_{H^{4}(\Gamma_*)}, \big|\widetilde d_l'\big|, \big|\widetilde d'_r\big|\le L_1.
\]
On the other hand, using  system \eqref{equ:map1}, Lemma \ref{R0 estimate} and Corollary \ref{pa 2 t f estimate} with $s=2.5$, we arrive at
\[
\begin{split}
\|\pa^2_t\widetilde{\mN_a}\|_{H^{2.5}(\G_*)}\le &\|R_0\|_{H^{2.5}(\G_*)}+C\big(\|\mN_a\|_{H^{1.5}(\G_*)}, \|\pa_t \mN_a\|_{H^{1}(\G_*)},  a^{-1}\|\pa^2_t \mN_a\|_{H^{1}(\G_*)}, |d_i|, |d'_i|, a^{-1}|d''_i|\big)\times\\
&\qquad\big(\|\widetilde{\mN_a}\|_{H^{5.5}(\G_*)}+\|\pa_t \widetilde{\mN_a}\|_{H^{3.5}(\G_*)}\big)\\
\le & C(L_0, L_1,  a^{-1} L_2).
\end{split}
\] 
Besides, one has directly the estimate
\[
\big|\widetilde d''_i\big|\le C(L_0, L_1).
\] 

Therefore, when $a$ taken sufficiently large, we  have
\[
\big\|  \pa^2_t \widetilde{\mN_a}\big\|_{H^{2.5}(\Gamma_*)}, \big|\widetilde d_l''\big|,  \big|\widetilde d_r''\big|\le L_2,
\]
with  large $L_2$ compared to $L_0, L_1$.

Consequently, checking the initial data for $(\widetilde{\mN_a}, \widetilde d_l, \widetilde d_r)$, one  concludes that when $L\ge A_1$ is sufficiently large,  one will have $(\widetilde{\mN_a}, \widetilde d_l, \widetilde d_r)\in \Sigma$.

\medskip
Based on the analysis above, we can define the iteration map $\cF$ on $\Sigma$ by solving the linear system \eqref{equ:map1}:
\beno
\cF\big( \mN_a, d_l, d_r\big) \triangleq (\widetilde{\mN_a}, \widetilde d_l, \widetilde d_r),
\eeno
with initial data  $\big ((\mN_a)_I,(\pa_t\mN_a)_I), (d_i)_I, (d'_i)_I\big) \in \mathcal{I}(\eps, A_1)$.  When $(\mN_a, d_l, d_r) \in \Sigma$, we have proved that
\[
(\widetilde{\mN_a}, \widetilde d_l, \widetilde d_r)=\cF\big( \mN_a, d_l, d_r\big) \in \Sigma.
\]

\bthm{Remark}\label{angles}
{\it
Thanks to the assumptions on the reference domain $\Om_*$ such that $\om_{i*}\in (0, \pi/16)$,  we know that when $T$ and the constant $\del_1$ in $\Sigma$ are chosen small enough, the contact angles generated by both $( \mN_a, d_l, d_r\big)$ and $(\widetilde{\mN_a}, \widetilde d_l, \widetilde d_r)$  lie in the same interval $(0, \pi/16)$ as well.  
}
\ethm

\subsection{The contraction mapping}
Firstly, we introduce a lower-order norm $\|\cdot\|_{\Sigma}$ for $\Sigma$:
\[
\begin{split}
\big\|\big(\mN_a,  d_l, d_r\big)\big\|_{\Sigma}&\triangleq\|\mN_a\|_{C([0,T];H^{4}(\Gamma_*))}+\big\|  \pa_t \mN_a \big\|_{C(0,T]; H^{2.5}(\Gamma_*))}+a^{-3/2}\big\|  \pa^2_t \mN_a \big\|_{C([0,T]; H^{1}(\Gamma_*))}\\
&\quad +\sum_i\Big(| d_i|_{C([0,T])}+\big| d'_i\big|_{C([0,T])}+a^{-3/2}\big|  d''_i\big|_{C([0,T])}\Big).
\end{split}
\]

To prove that  $\cF$ is a contraction mapping, we use the Banach fixed point theorem in a variational way.  In fact, let's consider a one-parameter family $(\mN_a(\tau), d_l(\tau), d_r(\tau))\in \Sigma$ with initial data
\[
\big((\mN_a)_I(\tau), (\pa_t\mN_a)_I(\tau),  (d_i)_I(\tau), (d'_i)_I(\tau)\big).
\]

Correspondingly, we  write
\[
\cF( \mN_a(\tau),d_l(\tau), d_r(\tau)) = \big(\widetilde{\mN_a}(\tau), \widetilde d_l(\tau), \widetilde d_r(\tau)\big).
\]
Taking the variation with respect to $\tau$ on both sides of system \eqref{equ:map1} leads to a system  for $\big(\pa_\tau \widetilde{\mN_a}, \pa_\tau \widetilde d_l, \pa_\tau \widetilde d_r\big)$ immediately:
 \beq\label{equ:map}
 \left\{\begin{array}{ll}
D^2_{t*} \pa_\tau \widetilde{\mN_a}+\sigma \cA (d_{\G_t})\pa_\tau\widetilde{\mN_a} =G_1
\qquad\hbox{on}\quad \G_*,\\
D_{t*}\cA(d_{\G_t})\pa_\tau\widetilde{\mN_a}\pm\f{\si^2}{\be_c}(\sin\om_i)^2\circ\Phi_{S_t}\big(\na_{\tau_t}(\cA(d_{\G_t})\pa_\tau\widetilde{\mN_a}\circ\Phi^{-1}_{S_t})\big)\circ\Phi_{S_t}=G_{2, i}\qquad\hbox{at}\quad p_{i*},\\
 \pa^2_t\pa_\tau\widetilde d_i=\pa_\tau B_i\quad(i=l, r)
\end{array}
\right.
\eeq
with
\[
\begin{split}
&G_1\circ \Phi^{-1}_{S_t}=(\pa_\tau R_0)\circ \Phi^{-1}_{S_t}
-[D_\tau, D^2_t]\widetilde{\mN_a}\circ \Phi^{-1}_{S_t}
+\sigma [D_\tau, \cN\D_{\G_t}]\widetilde{\mN_a}\circ \Phi^{-1}_{S_t},\\
&G_{2, i}\circ \Phi^{-1}_{S_t}=(\pa_\tau  R_{c, i})\circ \Phi^{-1}_{S_t}+[D_{t}\cN\D_{\G_t},  D_\tau]\widetilde{\mN_a}\circ \Phi^{-1}_{S_t}
-\f{\si^2}{\be_c}D_\tau(\sin\om_i)^2\na_{\tau_t}\cN\D_{\G_t}\big(\widetilde{\mN_a}\circ\Phi^{-1}_{S_t}\big)\\
&\qquad \qquad\qquad -\f{\si^2}{\be_c}(\sin\om_i)^2[D_\tau, \na_{\tau_t}\cN\D_{\G_t}]\widetilde{\mN_a}\circ\Phi^{-1}_{S_t},
\end{split}
\]
where $R_{0}, R_{c, i}$ are defined before
and  moreover
\[
\pa_\tau B=\pa_\tau B\big(a^3\pa_\tau \pa d_{\G_t}, \na_{v^*}\pa_t\pa_\tau d_{\G_t}, \pa_\tau v^*,  \pa_\tau \mN_a, \pa_\tau \na P_{v,v}\big).
\]

In fact, the condition at the corner point in \eqref{equ:map} is rewritten from the condition in system \eqref{equ:map1}:
\[
D_{t*}\cA(d_{\G_t})\widetilde{\mN_a}\pm\f{\sigma^2 }{\beta_c}(\sin \om_i)^2(\na_{\tau_t}(\cA(d_{\G_t})\widetilde{\mN_a}\circ \Phi_{S_t}^{-1}))\circ\Phi_{S_t}= R_{c, i}\qquad \textrm{at} \quad p_{i*},
\]
which is equivalent to
\[
D_t\cN\D_{\G_t}\big(\widetilde{\mN_a}\circ \Phi_{S_t}^{-1}\big)\pm\f{\sigma^2 }{\beta_c}(\sin \om_i)^2\na_{\tau_t}\cN\D_{\G_t}\big(\widetilde{\mN_a}\circ \Phi_{S_t}^{-1}\big)=R_{c, i}\circ \Phi_{S_t}^{-1}\qquad \textrm{at} \quad p_{i}\ (i=l, r).
\]
Taking $D_\tau$ on both sides of the equation above, one can derive the desired condition in \eqref{equ:map}.

To simplify the notations, we denote
\[
F(\tau)=\big(\pa_\tau\widetilde{\mN_a}\big)\circ \Phi^{-1}_{S_t}(\tau).
\]
 and rewrite \eqref{equ:map} as
\beq\label{F linear system}
 \left\{\begin{array}{ll}
D^2_{t}F-\sigma \mathcal{N} \Delta_{\Gamma_{t}}F =G_1\circ \Phi^{-1}_{S_t}\qquad\hbox{on}\quad \G_t,\\
D_{t}\cN\D_{\G_t}F\pm\f{\si^2}{\be_c}(\sin\om_i)^2\na_{\tau_t}\cN\D_{\G_t}F=G_{2, i}\circ \Phi^{-1}_{S_t}\qquad\hbox{at}\quad p_i,\\
\pa^2_t\pa_\tau\widetilde d_i=\pa_\tau B_i\qquad (i=l, r).
\end{array}
\right.
\eeq

In order to prove the contraction mapping, we need a lower-order energy estimate for $(F, \pa_\tau \widetilde d_l, \pa_\tau \widetilde d_r)$ under the norm $\|\cdot\|_\Sigma$, which is similar to the energy estimates obtained in Section 5.
 \bthm{Proposition}\label{prop:energy est lower order}(Lower-order energy estimates)
{\it Let $(\mN_a,\,d_l, \,d_r)\in \Sigma$  be given.
Assume that the contact angles $\om_i$ satisfy
\[
\min_i\sin \om_i\ge c_0\qquad \hbox{for some }\ c_0>0.
\]

Then  we have the following energy estimates for the linear problem \eqref{F linear system}:
\[
\begin{split}
&
\| F\|^2_{H^4(\G_t)}+\|D_{t} F\|^2_{H^{2.5}(\G_t)}\\
&\le e^{C(L_0, L_1) t}C(\del_1, L)\big(\|F(0)\|^2_{H^4(\G_0)}+\|D_t F(0)\|^2_{H^{2.5}(\G_0)}\big)+a^{3/2}e^{C(L_0, L_1) t}C(L_0, L_1)\int^t_0\Big( \|\pa_\tau\mN_a\|^2_{H^{4}(\G_*)}\\
&\qquad\qquad +\|\pa_{t}\pa_\tau\mN_a\|^2_{H^{2.5}(\G_*)}
+a^{-3}\|\pa^2_{t}\pa_\tau\mN_a\|^2_{H^{1}(\G_*)}+|\pa_\tau d_i|^2+|\pa_t\pa_\tau d_i|^2+a^{-3}|\pa^2_t\pa_\tau d_i|^2\Big)dt',
\end{split}
\]
and
\[
\begin{split}
\big|\pa_\tau d_i(t)\big|+\big|\pa_t\pa_\tau d_i(t)\big|&\le \big|\pa_\tau d(0)\big|+(1+t)\,\big|\pa_t\pa_\tau d_i(0)\big|\\
&\qquad +C(L_0, L_1)\int^t_0\big( \|\pa_\tau\mN_a\|_{H^{1}(\G_*)}+\|\pa_{t}\pa_\tau\mN_a\|_{H^{1}(\G_*)} +|\pa_\tau d_i|+|\pa_t\pa_\tau d_i|\big) dt',
\end{split}
\]
where the constant $C$ depends on $L_0, L_1$ and the index $i$ in the estimates stands for the summation on $i=l, r$.
}
\ethm
\begin{proof} The proof is much similar as the proof for the higher-order energy estimates, and the lower-order energy and dissipation \eqref{lower energy} are used here.

To begin with, we  apply $\mathcal N \Delta_{\Gamma_t}$ on   both sides of the equation of $F$ from \eqref{F linear system} to obtain
\[
\cN\D_{\G_t}D^2_tF-\si (\cN\D_{\G_t})^2F=\cN\D_{\G_t}(G_1\circ \Phi^{-1}_{S_t})\qquad \hbox{on}\quad \G_t.
\]
Multiplying $D_t \Delta_{\Gamma_t} F$ on  both sides the equation above and integrating on $\Gamma_t$, one has
\beq\label{integral equation lower order}
\int_{\Gamma_t} \mathcal{N} \Delta_{\Gamma_t } D^2_t F\cdot D_t\Delta_{\Gamma_{t}} Fds
 -\sigma\int_{\Gamma_t} \mathcal{N} \Delta_{\Gamma_t} \mathcal{N} \Delta_{\Gamma_{t}}F
\cdot D_t   \Delta_{\Gamma_t}Fds
=\int_{\Gamma_t} \mathcal{N} \Delta_{\Gamma_t}( G_1\circ \Phi^{-1}_{S_t})\cdot  D_t   \Delta_{\Gamma_t}Fds.
\eeq

Firstly, for the first term on the left side, one derives
\[
\int_{\Gamma_t} \mathcal{N} \Delta_{\Gamma_t } D^2_t F\cdot D_t\Delta_{\Gamma_{t}} Fds
=\f12\pa_t\big\| \na \cH( D_t \Delta_{\Gamma_{t } } F)\big\|^2_{L^2(\Om_t)}+A_{l1},
\]
where
\[
A_{l1}=\int_{\Om_t}[ \na \cH, D_t]D_t \Delta_{\Gamma_{t } } F\cdot \na \cH (D_t   \Delta_{\Gamma_t}F)dX
-\f12\int_{\Om_t}(\dive v)\big| \na \cH( D_t \Delta_{\Gamma_{t } } F)\big|^2dX
+\int_{\Gamma_t} [\Delta_{\Gamma_{t } }, D^2_t ]F \cdot  \cN D_t   \Delta_{\Gamma_t}F ds.
\]
Applying \eqref{commutator Dt H}, \eqref{commutator surface delta},  Proposition \ref{d and  N ka}, Corollary \ref{pa tau d estimate}, Lemma \ref{Dt v estimate} and the bounds in $\Sigma$ implies
\[
\begin{split}
|A_{l1}|&\le C\big(\|d_{\G_t}\|_{H^3(\G_*)}, \|v\|_{H^{3.5}(\Om_t)}, \|D_t v\|_{H^{2.5}(\Om_t)}\big)\big(E_l(t, F, \pa_t F)+\|F\|^2_{H^{2.5}(\G_t)}+\|D_tF\|^2_{H^{2.5}(\G_t)}\big)\\
&\le C(L_0, L_1)\big(E_l(t, F, \pa_t F)+\|F\|^2_{H^{2.5}(\G_t)}+\|D_tF\|^2_{H^{2.5}(\G_t)}\big).
\end{split}
\]

Secondly, for the second term on the left side, one can show that
\[
-\si\int_{\Gamma_{t}} \mathcal{N} \Delta_{\Gamma_t} \mathcal{N} \Delta_{\Gamma_t }F
\cdot D_t   \Delta_{\Gamma_t}F\,ds
=\f\si2\pa_t\big\|\na_{\tau_t} \mathcal{N} \Delta_{\Gamma_t }F\big\|^2_{L^2(\G_t)}
+ A_{l2}+A_{l3},
\]
with the remainder term
\[
\begin{split}
A_{l2}&=\si \int_{\Gamma_{t}}\na_{\tau_t} \mathcal{N} \Delta_{\Gamma_t }F
\cdot [ \na_{\tau_t}\cN, D_t]  \Delta_{\Gamma_t}F ds
-\f\si 2\int_{\G_t}(\cD\cdot v)\big|\na_{\tau_t} \mathcal{N} \Delta_{\Gamma_t }F\big|^2 ds\\
&\le C(\|d_{\G_t}\|_{H^{4}(\G_*)}, \|v\|_{H^{2.5}(\Om_t)}\big)\big(E_l(t, F, \pa_t F)+\|F\|^2_{H^{3.5}(\G_t)}\big)\\
&\le C(L_0, L_1)\big(E_l(t, F, \pa_t F)+\|F\|^2_{H^{3.5}(\G_t)}\big),
\end{split}
\]
and the term at the corner points
\[
A_{l3}=-\si\na_{\tau_t} \mathcal{N} \Delta_{\Gamma_t }F\cdot   \mathcal{N}  D_t \Delta_{\Gamma_t}F\,\big|^{p_l}_{p_r}.
\]

To handle $A_{l3}$, we use the condition for the corner points in \eqref{F linear system}:
\[
D_{t}\cN\D_{\G_t}F\pm\f{\si^2}{\be_c}(\sin\om_i)^2\na_{\tau_t}\cN\D_{\G_t}F=G_{2, i}\circ \Phi^{-1}_{S_t}\qquad\hbox{at}\quad p_i\ (i=l, r),
\]
which can be rewritten as
\[
\mathcal{N}  D_t \Delta_{\Gamma_t}F=\mp\f{\si^2}{\be_c}(\sin\om_i)^2\na_{\tau_t}\cN\D_{\G_t}F+r_{4, i}\qquad\hbox{at}\quad p_i,
\]
where
\[
r_{4, i}=-[D_t, \cN]\D_{\G_t}F+G_{2, i}\circ \Phi^{-1}_{S_t}.
\]
Consequently, we have
\[
A_{l3}=\f{\si^3}{\be_c}F_l(t, F, \pa_tF)-\si \na_{\tau_t} \mathcal{N} \Delta_{\Gamma_t }F\cdot r_{4, i},
\]
where the remainder term satisfies
\[
\big|\si \na_{\tau_t} \mathcal{N} \Delta_{\Gamma_t }\,F\cdot r_{4, i}\big|\le \del F_l(t, F, \pa_tF)
+(\sin\om_i)^{-2}C_\del\|r_{4, i}\|^2_{H^1(\G_t)}.
\]
Moreover, checking he terms in $G_2$ carefully and apply Lemma \ref{lem:R_c} to show that
\[
\begin{split}
\|r_{4, i}\|_{H^1(\G_t)}
&\le C\big(\|d_{\G_t}\|_{H^5(\G_*)}, \|\pa_td_{\G_t}\|_{H^4(\G_*)}, \|v\|_{H^{5.5}(\Om_t)}, \|\mN_a\|_{H^3(\G_*)}, \|\pa_t\mN_a\|_{H^2(\G_*)},\|\widetilde{\mN_a}\|_{H^5(\G_*)},
\\
&\qquad\quad   \|\pa_t\widetilde{\mN_a}\|_{H^4(\G_*)}\big)
\Big(\|\pa_\tau d_{\G_t}\|_{H^5(\G_*)}+\|\pa_t\pa_\tau d_{\G_t}\|_{H^4(\G_*)}+\|D_\tau v\|_{H^{5.5}(\Om_t)}+a\|\pa_\tau\mN_a\|_{H^4(\G_*)}\\
&\qquad\quad  +a\|\pa_t\pa_\tau\mN_a\|_{H^{2.5}(\G_*)}+a|\pa_\tau d_i|+a|\pa_t\pa_\tau d_i|
\Big),
\end{split}
\]
and applying Proposition \ref{d and  N ka}, Corollary \ref{pa tau d estimate}, Lemma \ref{Dt v estimate} and the bounds in $\Sigma$ leads to
\[
\|r_{4, i}\|_{H^1(\G_t)}\le a\, C\big(L_0, L_1)\Big(\|\pa_\tau\mN_a\|_{H^4(\G_*)}+\|\pa_t\pa_\tau\mN_a\|_{H^{2.5}(\G_*)}+|\pa_\tau d_i|+|\pa_t\pa_\tau d_i|
\Big).
\]

On the other hand, for the right side of \eqref{integral equation lower order}, we have
\[
\begin{split}
\int_{\Gamma_t} \mathcal{N} \Delta_{\Gamma_t} (G_1\circ \Phi^{-1}_{S_t})\cdot  D_t   \Delta_{\Gamma_t}Fds
&= \int_{\Om_t}\na\cH\big(\Delta_{\Gamma_t} (G_1\circ \Phi^{-1}_{S_t})\big)\cdot \na\cH(D_t   \Delta_{\Gamma_t}F)dX\\
&\le C(\|d_{\G_t}\|_{H^{2.5}(\G_*)})\big\|G_1\circ \Phi^{-1}_{S_t}\big\|_{H^{2.5}(\G_t)}E_l(t, F, \pa_t F)^{1/2}.
\end{split}
\]
Checking the expression of $G_1$, one obtains by Lemma \ref{Dt v estimate} (1) and lemma \ref{R0 estimate} that
\[
\begin{split}
\big\|G_1\circ \Phi^{-1}_{S_t}\big\|_{H^{2.5}(\G_t)}&\le C\big(\|d_{\G_t}\|_{H^{6.5}(\G_*)}, \|\pa_td_{\G_t}\|_{H^{5.5}(\G_*)}, \|v\|_{H^{5}(\Om_t)}, \|\mN_a\|_{H^{3.5}(\G_*)}, \|\widetilde{\mN_a}\|_{H^{5.5}(\G_*)},
\\
&\qquad\quad   \|\pa_t\widetilde{\mN_a}\|_{H^{3.5}(\G_*)}\big)
\Big(\|\pa_\tau d_{\G_t}\|_{H^{6.5}(\G_*)}+\|\pa_t\pa_\tau d_{\G_t}\|_{H^{5.5}(\G_*)}
+\|\pa^2_t\pa_\tau d_{\G_t}\|_{H^{2.5}(\G_*)}\\
&\qquad\quad+\|D_\tau v\|_{H^{5}(\Om_t)}+a^{3/2}\big(\|\pa_\tau\mN_a\|_{H^{4}(\G_*)}+\|\pa_t\pa_\tau\mN_a\|_{H^{2.5}(\G_*)}+|\pa_\tau d_i|+|\pa_t\pa_\tau d_i|\big)
\Big),
\end{split}
\]
where $\|\pa^2_t\pa_\tau d_{\G_t}\|_{H^{2.5}(\G_*)}$ comes from the term $[D_\tau, D^2_t]\widetilde{\mN_a}\circ \Phi^{-1}_{S_t}$ in $G_1$.
Applying Proposition \ref{d and  N ka}, Corollary \ref{pa tau d estimate}, Lemma \ref{Dt v estimate}  (2)  and the bounds in $\Sigma$ again leads to
\[
\begin{split}
\big\|G_1\circ \Phi^{-1}_{S_t}\big\|_{H^{2.5}(\G_t)}&\le a^{3/2}C\big(\|\mN_a\|_{H^{3.5}(\G_*)}, \|\pa_t\mN_a\|_{H^{2.5}(\G_*)}, \|\widetilde{\mN_a}\|_{H^{5.5}(\G_*)}, \|\pa_t\widetilde{\mN_a}\|_{H^{3.5}(\G_*)}, |d_i|, |\pa_t d_i|\big)\times\\
&\qquad\quad \Big(\|\pa_\tau\mN_a\|_{H^{4}(\G_*)}+\|\pa_t\pa_\tau\mN_a\|_{H^{2.5}(\G_*)}
+a^{-3/2}\|\pa^2_t\pa_\tau\mN_a\|_{H^{1}(\G_*)}+|\pa_\tau d_i|\\
&\qquad\quad+|\pa_t\pa_\tau  d_i|+a^{-3/2}|\pa^2_t\pa_\tau d_c|\Big)\\
&\le a^{3/2}C(L_0, L_1)\Big(\|\pa_\tau\mN_a\|_{H^{4}(\G_*)}+\|\pa_t\pa_\tau\mN_a\|_{H^{2.5}(\G_*)}
+a^{-3/2}\|\pa^2_t\pa_\tau\mN_a\|_{H^{1}(\G_*)}+|\pa_\tau d_i|\\
&\qquad\quad+|\pa_t\pa_\tau  d_i|+a^{-3/2}|\pa^2_t\pa_\tau d_i|\Big).
\end{split}
\]

Consequently, since the energy estimates for $\|D_tF\|_{L^2(\G_t)}, \|F\|_{L^2(\G_t)}$ are similar as in the higher-order case,  summing up all the estimates above leads to
\[
\begin{split}
&\pa_t E_l(t, F, \pa_t F)+\sum_i F_{l, i}(t, F, \pa_t F)\\
&\le C(L_0, L_1)\Big( E_l(t, F, \pa_t F)+\|F\|^2_{H^{3.5}(\G_t)}+\|D_tF\|^2_{H^{2.5}(\G_t)}+a^{3/2}\big(\|\pa_\tau\mN_a\|_{H^{3.5}(\G_*)}+\|\pa_t\pa_\tau\mN_a\|_{H^{2.5}(\G_*)}\\
&\qquad \quad
+a^{-3/2}\|\pa^2_t\pa_\tau\mN_a\|_{H^{1}(\G_*)}+|\pa_\tau d_i|+|\pa_t\pa_\tau  d_i|+a^{-3/2}|\pa^2_t\pa_\tau d_i|\big)\Big).
\end{split}
\]
Applying Lemma \ref{E(t) equivalence} with the bounds in $\Sigma$ for the coefficient,
the desired energy estimate for $F$ can be finished. Moreover, the energy estimate for $\pa_\tau d_i$ $(i=l, r)$ can be done similarly and more easily.
\end{proof}

\bigskip

In the end, we finish the iteration.
\bthm{Proposition}{\it
For any $0<\eps<<\delta_1$ and $A_1>0$ and initial data $\big((\mN_a)_I, (\pa_t\mN_a)_I, (d_i)_I, (d'_i)_I\big)\in \mathcal I(\eps, A_1)$, there exists $L, L_0, L_1$ and small constants $a^{-1}, T$ such that $\cF$ defined on $\Sigma$ has a fixed point.
}
\ethm
\begin{proof} Thanks to Proposition \ref{prop:energy est lower order}, when $T$ is small enough we have the following estimate
\beq\label{pa t 0 and 1 estimate}
\begin{split}
&\|\pa_\tau \widetilde{\mN_a}\|^2_{H^4(\G_*)}+\|\pa_t\pa_\tau\widetilde{\mN_a}\|^2_{H^{2.5}(\G_*)}+\sum_i\big(|\pa_\tau d_i|^2+|\pa_t\pa_\tau d_i|^2\big)\\
&\le  C(L_0, \del_1, L)\Big(\|\pa_\tau (\mN_a)_I\|^2_{H^4(\G_*)}+\big\|\pa_\tau\big(\pa_t \mN_a\big)_I\big\|^2_{H^{2.5}(\G_*)}+\sum_i\big(|\pa_\tau (d_i)_I|^2+|\pa_\tau (\pa_t d_i)_I|^2\big)\Big)\\
&\qquad +T\, a^{3}C( L_0, L_1,  \del_1, L)\big\|(\pa_\tau\mN_a, \pa_\tau d_l, \pa_\tau d_r)\big\|^2_\Sigma,
\end{split}
\eeq
where one can see that the estimates for the second-order time derivatives are still missing.

In fact, applying Corollary \ref{pa 2 t f estimate} with $s=1$, $f$ replaced by $\pa_\tau \widetilde{\mN_a}$ and $g_1$ replaced by $G_1$, one obtains immediately that
\beq\label{pa t 2 estimate}
\begin{split}
\|\pa^2_t\pa_\tau \widetilde{\mN_a}\|_{H^1(\G_*)}\le &\|G_1\|_{H^1(\G_*)}+C\big(\|\mN_a\|_{H^{1}(\G_*)}, \|\pa_t \mN_a\|_{H^{1}(\G_*)},  a^{-1}\|\pa^2_t \mN_a\|_{H^{1}(\G_*)}, |d_i|, |d'_i|, a^{-1}|d''_i|\big)\times\\
&\qquad\big(\|\pa_\tau \widetilde{\mN_a}\|_{H^{4}(\G_*)}+\|\pa_t \pa_\tau \widetilde{\mN_a}\|_{H^{2}(\G_*)}\big).
\end{split}
\eeq
Moreover, one has similarly as before the estimate
\[
\begin{split}
\|G_1\|_{H^1(\G_*)}&\le C\big(L_0, L_1)\Big(\|\pa_\tau \mN_a\|_{H^4(\G_*)}+\|\pa_t\pa_\tau \mN_a\|_{H^{2.5}(\G_*)}+a^{-3/2}\|\pa^2_t\pa_\tau \mN_a\|_{H^1(\G_*)}+|\pa_\tau d_i|+|\pa_t\pa_\tau d_i|\\
&\qquad \quad+a^{-3/2}|\pa^2_t\pa_\tau d_i|\Big).
\end{split}
\]
As a result, substituting this inequality and  \eqref{pa t 0 and 1 estimate} into \eqref{pa t 2 estimate},  one arrives at
\[
\begin{split}
\|\pa^2_t\pa_\tau \widetilde{\mN_a}\|^2_{H^1(\G_*)}&\le C(L_0, L_1, \del_1, L)\big\|(\pa_\tau\mN_a, \pa_\tau d_c)\big\|^2_\Sigma\\
&\quad +C(L_0, L_1, \del_1, L)\Big(\|\pa_\tau (\mN_a)_I\|^2_{H^4(\G_*)}+\big\|\pa_\tau\big(\pa_t \mN_a\big)_I\big\|^2_{H^{2.5}(\G_*)}+\sum_i\big(|\pa_\tau (d_i)_I|^2+|\pa_\tau (\pa_t d_i)_I|^2\big)\Big).
\end{split}
\]

On the other hand, using the equation for $\pa_\tau\widetilde d_i$ directly leads to the estimate
\[
\big|\pa^2_t\pa_\tau\widetilde d_i\big|^2\le C(L_0, L_1)\big( \|\pa_\tau\mN_a\|^2_{H^{1}(\G_*)}+\|\pa_{t}\pa_\tau\mN_a\|^2_{H^{1}(\G_*)} +|\pa_\tau d_i|^2+|\pa_t\pa_\tau d_i|^2\big).
\]

Summing up all the estimates above, we finally conclude that
\[
\begin{split}
&\big\|(\pa_\tau\widetilde{\mN_a}, \pa_\tau\widetilde d_l, \pa_\tau\widetilde d_r)\big\|^2_\Sigma\\
&\le  (Ta^{3}+a^{-3})C(L_0, L_1, \del_1, L)\big\|(\pa_\tau\mN_a, \pa_\tau d_l, \pa_\tau d_r)\big\|^2_\Sigma + C(L_0, L_1, \del_1, L)\Big(\|\pa_\tau (\mN_a)_I\|^2_{H^4(\G_*)}+\big\|\pa_\tau\big(\pa_t \mN_a\big)_I\big\|^2_{H^{2.5}(\G_*)} \\
&\qquad +\sum_i\big(|\pa_\tau (d_i)_I|^2+|\pa_\tau (\pa_t d_i)_I|^2\big)\Big).
\end{split}
\]
When we  take $a^{-1}$  and then  $T$ sufficiently small, we have
\[
\begin{split}
&\big\|(\pa_\tau\widetilde{\mN_a}, \pa_\tau\widetilde d_l, \pa_\tau\widetilde d_r)\big\|_\Sigma \\
&\leq \f12  \|\big( \pa_\tau\mN_a, \pa_\tau d_l, \pa_\tau d_r \big)\|_{\Sigma} +C\Big(\|\pa_\tau (\mN_a)_I\|_{H^4(\G_*)}+\big\|\pa_\tau\big(\pa_t \mN_a\big)_I\big\|_{H^{2.5}(\G_*)}+\sum_i\big(|\pa_\tau (d_c)_I|+|\pa_\tau (\pa_t d_c)_I|\big)\Big).
\end{split}
\]
Therefore, if we fix  initial data in $\mathcal{I}(\eps, A_1)$, we construct  a Cauchy sequence $\{ \big((\mN_a)^n, d_l^n, d^n_r \big)\}$ in $\Sigma$. Consequently,  there exist a unique  $(\mN_a, d_l, d_r)\in \Sigma$ such that
\beno
\big(  (\mN_a)^n, d_l^n, d_r^n \big)\to  \big(\mN_a, d_l, d_r\big)\qquad\hbox{in}\quad \Sigma,
\eeno
which implies that $\cF$ defined on  $\Sigma$ admits a fixed point.
\end{proof}

\medskip

\subsection{Back to the Euler equation}
In last section, we have proved the unique existence of the solution $(\mN_a, d_l, d_r)$ of \eqref{ka_a equation} and \eqref{d c eqn} for given initial data. As a result,  we construct the moving domain $\Om_t$ and the velocity field $v$ by $d_{\Gamma_t}$. 

For the moment, we are able to show  that the water-waves system $\mbox{(WW)}$ is satisfied by this velocity $v$ and the  pressure
\[
P=\sigma \kappa_{\cH}+ P_{v, v}.
\]

To begin with, we firstly recall the definition of $v=\na \phi$ by \eqref{v potential expression}. So one has 
\beno
D_t v=D_t \na \phi=\na(D_t \phi)-\na v\cdot \na \phi=\na (D_t\phi-\f12|v|^2).
\eeno
Now, we define the other  pressure
\[
Q=-D_t\phi+\f12|v|^2,
\]
which infers
\beno
D_t v+\na Q=0. 
\eeno

Moreover, we recall that
\[
\dive v=\gamma \xi
\]
with
\[
\gamma=\big(\int_{\Om_t}dX\big)^{-1},\quad \xi=\int_{\Gamma_t} v\cdot N_t\,ds.
\]

We are going to prove
\[
P+gz=Q\quad\hbox{ and} \quad \dive v=0.
\]
In order to do this, we firstly define
\[
V_0\triangleq\na (Q-P-gz)=-D_t v-\na P+{\bf g},
\]
so we have a slightly new equation for $v$ compared to Euler equation:
\beq\label{new v eqn}
D_t v=-\si\na\ka_\cH-\na P_{v,v}+{\bf g}-V_0.
\eeq
We go through the computations for deriving \eqref{ka_a equation} again to check the terms involving $V_0$, which means we need to trace the substitution of $D_t v$ and Euler equation. To begin with, \eqref{ka equation} is rewritten here:
\[
D^2_t\ka=-N_t\cdot \D_{\G_t}D_t v+2\si \Pi(\tau_t)\cdot  \na_{\tau_t}\na\ka_{\cH}+r_1,
\]
where  $r_1$ is the remainder term defined before. Substituting the new equation \eqref{new v eqn} of $v$ into this equation leads to
\beno
D_t^2 \kappa=\sigma \Delta_{\Gamma_t} \cN( \kappa)+\widetilde{\mathfrak{R}}_1,
\eeno
where
\beno
\widetilde{\mathfrak{R}}_1 = \widetilde{R}_1+N_t\cdot \Delta_{\Gamma_t} V_0.
\eeno
Moreover, noticing that the term $R_a$ in \eqref{definition of R_0} contains $D_t v$ by \eqref{pa 2 t d} and \eqref{pa t v*}, one finds these extra terms in $R_a$ as below
\[
-a^3\f{1}{(\mu\circ\Phi^{-1}_{S_t})\cdot N_t}V_0\cdot N_t
+a^3\f{(\mu\circ \Phi_{S_t}^{-1})\cdot\na (d_{\G_t}\circ \Phi_{S_t}^{-1})}{(\mu\circ\Phi^{-1}_{S_t})\cdot N_t}V_0\cdot N_t
-a^3\,V_0\cdot \na (d_{\G_t}\circ \Phi_{S_t}^{-1}).
\]
Similarly, the term $D_t[D_{t}, \cN]\ka$ in $\widetilde R_1$ from \eqref{definition of R_0}  also involves $D_t v$. As a result,  we  finally show that
\beno
D_{t*}^2 \mN_a+\sigma \cA (\kappa_a) \mN_a=\mathfrak{R}_0\qquad\hbox{on}\quad \G_*.
\eeno
where
\[
\begin{split}
\mathfrak{R}_0\circ \Phi^{-1}_{S_t}&=R_0\circ \Phi^{-1}_{S_t}+\cN\big(N_t\cdot \Delta_{\Gamma_t} V_0\big)-a^3\f{1}{(\mu\circ\Phi^{-1}_{S_t})\cdot N_t}V_0\cdot N_t
+a^3\f{(\mu\circ \Phi_{S_t}^{-1})\cdot\na (d_{\G_t}\circ \Phi_{S_t}^{-1})}{(\mu\circ\Phi^{-1}_{S_t})\cdot N_t}V_0\cdot N_t\\
&\qquad
-a^3\,V_0\cdot \na (d_{\G_t}\circ \Phi_{S_t}^{-1})-\na_{N_t} \mathfrak{w}+\na_{N_t}V_0\cdot \na g_{\cH}+\na_{(\na g_{\cH})^\top} V_0 \cdot N_t,
\end{split}
\]
with $R_0$ defined in \eqref{definition of R_0} and
\[
\mathfrak{w}=\Delta^{-1}\big(2\na V_0\cdot \na^2 \ka_{\cH}+\Delta V_0\cdot \na \ka_{\cH}, (\na_{N_b} V_0-\na_{V_0}N_b)\cdot \na \ka_{\cH}\big).
\]
Since we have proved that \eqref{ka_a equation} holds, we know immediately that $V_0$ satisfies
\[
\begin{split}
&\cN\big(N_t\cdot \Delta_{\Gamma_t} V_0\big)-a^3\f{1}{(\mu\circ\Phi^{-1}_{S_t})\cdot N_t}V_0\cdot N_t+a^3\f{(\mu\circ \Phi_{S_t}^{-1})\cdot\na (d_{\G_t}\circ \Phi_{S_t}^{-1})}{(\mu\circ\Phi^{-1}_{S_t})\cdot N_t}V_0\cdot N_t-a^3\,V_0\cdot \na (d_{\G_t}\circ \Phi_{S_t}^{-1})\\
&\qquad-\na_{N_t} \mathfrak{w}+\na_{N_t}V_0\cdot \na g_{\cH}+\na_{(\na g_{\cH})^\top} V_0 \cdot N_t=0.
\end{split}
\]
Moreover, denoting
\[
\begin{split}
R_{V_0}&=\cN\big(2\na_{\tau_t}N_t\cdot \na_{\tau_t}V_0+\D_{\G_t}N_t\cdot V_0\big)-a^3\f{(\mu\circ \Phi_{S_t}^{-1})\cdot\na (d_{\G_t}\circ \Phi_{S_t}^{-1})}{(\mu\circ\Phi^{-1}_{S_t})\cdot N_t}V_0\cdot N_t+a^3\,V_0\cdot \na (d_{\G_t}\circ \Phi_{S_t}^{-1})\\
&\qquad+\na_{N_t} \mathfrak{w}-\na_{N_t}V_0\cdot \na g_{\cH}-\na_{(\na g_{\cH})^\top} V_0 \cdot N_t,
\end{split}\]
the above equation for $V_0$ becomes
\beq\label{equation for V0}
a^3\f{1}{(\mu\circ\Phi^{-1}_{S_t})\cdot N_t}V_0\cdot N_t-\cN\D_{\G_t}(V_0\cdot N_t)+R_{V_0}=0,
\eeq
which is analogous to the linearized equation of $\mN_a$ with respect to $d_{\G_t}$.

Secondly, recalling \eqref{Dt  v and d eqn} (which leads to the equation \eqref{d c eqn} of $d_i$) and using \eqref{new v eqn}, one obtains immediately with \eqref{d c eqn} and Remark \ref{V0} that
\beno
V_0\cdot N_t=0\quad \textrm{at}\quad p_i\ (i=l, r).
\eeno

On the other hand, a direct computation leads to the system for $P-Q+gz$:
\[
\left\{\begin{array}{ll}
\Delta (P-Q+gz)=D_t (\gamma\xi)=-div V_0\qquad \hbox{on}\quad \Om_t,\\
\na_{N_t}(P-Q+gz)\big|_{\G_t}=-V_0\cdot N_t,\quad \na_{N_b}(P-Q+gz)\big|_{\G_b}=0,
\end{array}
\right.
\]
which admits the following elliptic estimate by Theorem 5.3 \cite{MW} and Lemma \ref{trace thm PG} :
\beq\label{P-Q estimate}
\|V_0\|_{H^3(\Om_t)}\le C\|P-Q+gz\|_{H^{4}(\Om_t)}\leq C(L_0)\big(\|V_0\cdot N_t\|_{H^{2.5}(\Gamma_t)}+|\pa_t\xi|+|\xi|\big).
\eeq

\bigskip

Multiplying $\big(1-a^{-1}\D_{\G_t}\big)V_0\cdot N_t$ on  both sides of \eqref{equation for V0} and integrating by parts while using the boundary condition $V_0\cdot N_t\big|_{p_i}=0$,  one can have
\[
\begin{split}
&a^3\|V_0\cdot N_t\|^2_{L^2(\G_t)}+a^2\|\na_{\tau_t}(V_0\cdot N_t)\|^2_{L^2(\G_t)}+a^{-1}\|\D_{\G_t}(V_0\cdot N_t)\|^2_{H^{1/2}(\G_t)}\\
&\le C(L_0)\big(\|V_0\|_{H^{2.5}(\G_t)}+\|V_0\|_{H^2(\Om_t)}+\|V_0\|_{H^{1.5}(\G_b)}\big)\big(\|V_0\cdot N_t\|_{L^2(\G_t)}+a^{-1}\|\D_{\G_t}(V_0\cdot N_t)\|_{L^2(\G_t)}\big).
\end{split}
\]

Moreover, combining this estimate with \eqref{P-Q estimate}, we can conclude that
\[
\begin{split}
&a^3\|V_0\cdot N_t\|^2_{L^2(\G_t)}+a^2\|\na_{\tau_t}(V_0\cdot N_t)\|^2_{L^2(\G_t)}+a^{-1}\|\D_{\G_t}(V_0\cdot N_t)\|^2_{H^{1/2}(\G_t)}\\
&\le C(L_0)\big(\|V_0\cdot N_t\|_{H^{2.5}(\Gamma_t)}+|\pa_t\xi|+|\xi|\big)\big(\|V_0\cdot N_t\|_{L^2(\G_t)}+a^{-1}\|\D_{\G_t}(V_0\cdot N_t)\|_{L^2(\G_t)}\big).
\end{split}
\]
Consequently,  we obtain 
\beno
\|V_0\cdot N_t\|_{H^{2.5}(\Gamma_t)}\leq a^{-1}C(L_0)(|\pa_t\xi|+|\xi|).
\eeno

\medskip

For the moment, it remains to deal with $\xi, \pa_t\xi$. In fact, one has  by a direct calculation that
\beno
\pa_t \xi=\int_{\Om_t}D_t(\gamma\xi)dX+\int_{\Om_t}(\dive v)\gamma\xi dX=-\int_{\Om_t} div V_0 dX+\gamma\xi^2=-\int_{\Gamma_t}V_0\cdot N_tds+\gamma\xi^2.
\eeno
and 
\[
\gamma|\xi|^2\le a^{-1}C(L_0, L_1)|\xi|.
\]
Moreover, one knows
\[
\int_{\Gamma_t}V_0\cdot N_tds \le C(L_0)\|V_0\cdot N_t\|_{H^{2.5}(\G_t)},
\]
due to the fact that the domain is bounded.

Combining these inequalities with the above estimate for $\|V_0\cdot N_t\|_{H^{2.5}(\G_t)}$ and taking $a^{-1}$ sufficiently small, we derive
\beno
|\pa_t \xi |\leq a^{-1}C(L_0, L_1) |\xi|,
\eeno
which implies
\beno
\xi=0\qquad\hbox{since}\quad  \xi|_{t=0}=\int_{\Gamma_t}V_0\cdot N_tds\big|_{t=0} =0.
\eeno

As a result, we know immediately
\[
\dive v=0\quad\hbox{and}\quad V_0=0,
\]
which infers that Euler equation is satisfied by $v$:
\[
D_t v=-\na P+{\bf g}\qquad\hbox{on}\quad \Om_t.
\]

In the end, we will show that the condition for the corner points in \mbox{(WW)}
\beq\label{corner condition in ww}
\beta_c v_i=\sigma(\cos{\om_s}-\cos{\om_i})\qquad\hbox{at}\quad p_i\ (i=l,r)
\eeq
can be derived from
\[
D_{t*}\cA(\ka_a)\mN_a\pm\f{\sigma^2 }{\beta_c}(\sin \om_i)^2\big(\na_{\tau_t}(\cA(\ka_a)\mN_a\circ \Phi_{S_t}^{-1})\big)\circ\Phi_{S_t}= R_{c, i}\qquad\hbox{at}\quad p_{i*}.
\]
In fact, going back to the proofs for Lemma \ref{lem:cc} and Lemma 7.1 \cite{MW}, we know that this equation above is obtained by taking $D_t$ three times  on \eqref{corner condition in ww}, while Euler equation and the equation of $\mN_a$ is applied as well.  As a result,  integrating with respect to time variable three times and remembering that $v\cdot N_b\big|_{\G_b}=0$, we retrieve
 conditions for the contact points. 
 
 Moreover, one  knows from Remark \ref{angles} that the contact angles $\om_i$ stays in $(0, \pi/16)$ when $T$ is sufficiently small.

\bigskip

\section*{Acknowledgements}

The authors would like to thank Prof. Chongchun Zeng for fruitful discussions, and they also want to thank Nikolay Tzvetkov and Fr\'ederic Rousset for the comment on the boundary conditions. The author Mei Ming is supported by NSFC no.11401598.  The author Chao Wang is supported by NSFC no.11701016.


\begin{thebibliography}{99}




\bibitem{ABZ1} T. Alazard, N. Burq, C. Zuily, On the water-wave equations with surface tension. {\it Duke Math. J.}, (3){\bf 158}(2011), 413--499.


\bibitem{ABZ2} T. Alazard, N. Burq, C. Zuily, On the Cauchy problem for water gravity waves. {\it  Invent. Math.}, {\bf 198}(2014), 71--163.

\bibitem{ABZ3}T. Alazard, N. Burq, and C. Zuily, Cauchy theory for the gravity water waves system with non localized initial data, arXiv:1305.0457.

\bibitem{AD} T. Alazard and J.M. Delort, Global solutions and asymptotic behavior for two dimensional gravity water waves,  {\it Ann. Sci. Ec. Norm. Super.}, (5){\bf 48}(2015), 1149--1238..

\bibitem{AL} B. Alvarez-Samaniego and D. Lannes, Large time existence
for 3D water-waves and asymptotics, {\it Invent. Math.},
 {\bf 171}(2008), 485--541.

\bibitem{AM}D.M. Ambrose, Well-posedness of vortex sheets with surface tension. {\it SIAM J. Math. Anal.},  (1){\bf 35}(2003),  211--244.

\bibitem{AM1} D.M. Ambrose and N.  Masmoudi, The zero surface tension
limit of two-dimensional water waves, {\it Comm. Pure Appl.
Math.}, {\bf 58}(2005),  1287--1315.

\bibitem{AM2} D.M. Ambrose and N.  Masmoudi, The zero surface tension
limit of three-dimensional water waves,  {\it Indiana Univ. Math.
J.}, {\bf 58}(2009),  479--521.


\bibitem{BG1} K. Beyer, M. G\"unther, On the Cauchy problem for a capillary drop. I.  Irrotational motion,  {\it Math. Methods
Appl. Sci.}, (12){\bf 21}(1998), 1149--1183.


\bibitem{BEIMR} D. Bonn, J. Eggers, J. Indekeu, J. Meunier and E. Rolley, Wetting and spreading, {Reviews of modern Physics}, {\bf 81}(2009), 739--805.


\bibitem{CDA} A. Carlson, M. Do-Quang,  G. Amberg,  Modeling of dynamic wetting far from equilibrium, {\it Physics of Fluids}, (12){\bf 21}(2009), p121701.


\bibitem{CCFG} A. Castro, D. C\'ordoba, C. Fefferman, F. Gancedo and M. L\'opez-Fern\'andez, Rayleigh-Taylor breakdown for the Muskat problem with applications to water waves. {\it Ann. of Math.}, {\bf 175} (2012), 909--948.

\bibitem{CCFGG} A. Castro, D. C\'ordoba, C. Fefferman, F. Gancedo and J. G\'omez-Serrano, Finite time singularities for the free boundary incompressible Euler equations. {\it Ann. of Math.}, {\bf 178} (2013), 1061--1134.

\bibitem{BG2} K. Beyer, M. G\"unther, The Jacobi equation for irrotational free boundary flows, {\it Analysis (Munich)}, (3){\bf 20}(2000), 237--254.




\bibitem{Chang}  K.-C. Chang, {\it Methods in nolinear analysis}, Springer Monographs in Mathematics, 2014.



\bibitem{CL} D. Christodoulou, H. Lindblad,  On the motion of the free surface of a liquid, {\it Comm. Pure Appl. Math.}, (12){\bf 53}(2000), 1536--1602.

\bibitem{CS} D. Coutand and S. Shkoller, Well-posedness of the
free-surface incompressible Euler equations with or without
surface tension, {\it J. Amer. Math. Soc.}, {\bf 20} (2007), 829--930.

\bibitem{CS2} D. Coutand and S. Shkoller, On the finite-time splash and splat singularities for the 3-D free-surface Euler  equations, {\it Comm. Math. Phys.} {\bf 325} (2014), 143--183.

\bibitem{Craig} W. Craig, An existence theory for water waves and the
Boussinesq and Korteweg-de Vries scaling limits, {\it Comm.
Partial Differential Equations}, {\bf 10}(1985),  787--1003.


%\bibitem{Deng} Y. Deng,  A. D. Ionescu, B. Pausader  and F. Pusateri
%PUSATERI, Global solutions of the gravity-capillary water wave system in 3 dimensions, arXiv:1601.05685v1.

\bibitem{Deng}  Y. Deng, A. D. Ionescu, B. Pausader, and F. Pusateri. Global solutions for the 3D gravity-capillary water
waves system,{\it Acta Math.} 219 (2017), no. 2, 213--402.

\bibitem{GL}  J.-F. Gerbeau, T. Lelivre,  Generalized Navier boundary condition and geometric conservation law for surface tension, {\it Comput. Methods Appl. Mech. Engrg.}, {\bf 198 }(2009), 644--656.

\bibitem{GMS} P. Germain, N. Masmoudi and J. Shatah,
Global solutions for the gravity surface water waves equation in
dimension 3, {\it Ann. of Math.}, {\bf 175}(2012), 691--754.



\bibitem{PG1} P. Grisvard, Elliptic problems in non smooth domains, {\it Pitman Advanced Publishing Program, Boston-London-Melbourne}, 1985.


\bibitem{GT} Y. Guo and I. Tice, Stability of contact lines in fluids: 2D Stokes Flow,   {\it Arch. Ration. Mech. Anal.}, {\bf 227} (2018), no. 2, 767--854.

\bibitem{HIT1} J. Hunter, M. Ifrim and D. Tataru, Two dimensional water waves in holomorphic coordinates,  arXiv:1401.1252.

\bibitem{HIT2} M. Ifrim and D. Tataru, Two dimensional water waves in holomorphic coordinates II: global solutions,  arXiv:1404.7583.

\bibitem{HIT3} M. Ifrim and D. Tataru, The lifespan of small data solutions in two dimensional capillary water waves,  arXiv:1406.5471.

\bibitem{Iguchi} T. Iguchi,  Well-posedness of the initial value problem for capillary-gravity waves, {\it Funkcial. Ekvac.}, (2){\bf 44}(2001), 219--241.


\bibitem{Ig-Ta} T. Iguchi, N. Tanaka and A. Tani,  On a free boundary
problem for an incompressible ideal fluid in two space dimensions,
{\it Adv. Math. Sci. Appl.}, {\bf 9}(1999),  415--472.

\bibitem{MI} M. Ikawa, {\it A mixed problem for hyperbolic equations of second order with non-homogeneous Neumann type boundary condition}, Osaka J. Math. 6 1969 339--374.


\bibitem{IP} A. D. Ionescu and F. Pusateri, Global solutions for the gravity water waves system in 2D, {\it Invent. Math.}, (3){\bf 199}(2015), 653--804.

\bibitem{WuK} R.H. Kinsey and S. Wu, A Priori Estimates for Two-Dimensional Water Waves with Angled Crests, arXiv:1406.7573.


\bibitem{Lannes} D. Lannes,  Well-posedness of the water-wave equations, {\it Journal of the American Math. Society}, (3){\bf 18}(2005), 605--654.

\bibitem{LannesBook} D. Lannes, The water waves problem. Mathematical analysis and asymptotics, {\it Mathematical Surveys and Monographs}, Vol. 188. American Mathematical Society, Providence, RI, 2013. xx+321 pp.

\bibitem{Lannes1} D. Lannes, On the dynamics of floating structures, {\it Ann. PDE}, (1){\bf 3} (2017), Art. 11, 81 pp.

\bibitem{LI} D. Lannes and T. Iguchi, Hyperbolic free boundary problems and applications to wave-structure iterations, arXiv: 1806.07704v1.

\bibitem{LM} D. Lannes and G. M\'etivier, The shoreline problem for the one-dimensional shallow water and Green-Naghdi equations, {\it J. \'Ec. polytech. Math.}, {\bf 5}(2018), 455--518.

\bibitem{Lin} H. Lindblad, Well-posedness for the motion of an incompressible
liquid with free surface boundary, {\it Ann. of Math.},  {\bf
162}(2005), 109--194.


 \bibitem{MW1} M. Ming and C. Wang,  Elliptic estimates for D-N operator on corner domains, {\it Asymptotic analysis}, {\bf 104}(2017), 103--166.


\bibitem{MW}  M. Ming, C. Wang. {\it Water waves problem with surface tension in a corner domain I: A priori estimates with constrained contact angle}, , arXiv:1709.00180.

\bibitem{MZ} M. Ming and  Z. Zhang, Well-posedness of the water-wave problem with surface tension, {\it J. Math. Pures Appl.}, {\bf92}(2009), 429--455.

\bibitem{MZZ} M. Ming, P. Zhang and Z. Zhang, Large time well-posedness to the 3-D capillary-gravity waves in the long-wave regime, {\it Archive for Rational Mechanics and Analysis}, (2){\bf 204}(2012), 387--444.

\bibitem{Na} V.I. Nalimov, The Cauchy-Poisson problem (in
Russian), {\it Dynamika Splosh. Sredy,} {\bf 18}(1974), 104--210.

\bibitem{OT1} M.  Ogawa and A. Tani,  Free boundary problem for an
incompressible ideal fluid
 with surface tension, {\it Math. Models Methods Appl. Sci.}, {\bf 12}(2002), 1725--1740.

\bibitem{OT2} M. Ogawa and A. Tani,  Incompressible perfect fluid motion with
free boundary of finite depth, {\it Adv. Math. Sci. Appl.}, {\bf 13}(2003),  201--223.

\bibitem{Poyferre} T. de Poyferr\'e, A priori estimates for water waves with emerging bottom, arxiv: 1612.04103v1.

\bibitem{RE} W. Ren and W. E, Boundary conditions for the moving contact line problem, {\it Physics of Fluids}, {\bf 19}, 022101(2007), 1--15.

\bibitem{Sch}B.  Schweizer, On the three-dimensional Euler equations with a free boundary subject to surface tension, {\it Ann.
Inst. H. Poincar Anal. Non Linaire}, (6){\bf  22}(2005), 753--781.

\bibitem{SZ} J. Shatah and C. Zeng,  Geometry and a priori estimates for free boundary problems of the Euler equation, {\it Commun. Pure Appl. Math.}, {\bf 61}(2008), 698--744.

\bibitem{SZ2} J. Shatah and C. Zeng,  Local well-posedness for the fluid interface problems. {\it Arch. Ration. Mech. Anal.}, (2){\bf 199}(2011), 653--705.

\bibitem{SA} J.H. Snoeijer and B. Andreotti, Moving Contact Lines: Scales, Regimes, and Dynamical Transitions, {\it Annu. Rev. Fluid Mech.},  {\bf 45}(2013), 269Ð-292.

\bibitem{TZ} I. Tice and Y. Zheng,   Local well-posedness of the contact line problem in 2D Stokes flow. {\it  SIAM J. Math. Anal.},{\bf 49} (2017), no. 2, 899--953.

 \bibitem{Yo1} H. Yosihara, Gravity waves on the free surface
 of an incompressible perfect fluid of finite depth, {\it Publ. Res. Inst. Math. Sci.},
 {\bf 18}(1982),  49--96.

 \bibitem{Yo2} H. Yosihara, Capillary-gravity waves for an incompressible ideal fluid, {\it J. Math, Kyoto Univ.}, (4){\bf 23}(1983),  649--694.


\bibitem{Young} T. Young, An essay on the cohesion of fluids, {\it Philos. Trans. R. Soc. London}, (65){\bf 95}(1805).

\bibitem{WZZZ} C. Wang, Z. Zhang, W. Zhao and Y. Zheng, Local well-posedness and break-down criterion of the incompressible Euler equations with free boundary, to appear in {\it Memoris of the American Mathematical Society}.

%\bibitem{Wang2} X. Wang. On 3D water waves system above a flat bottom. {\it Anal. PDE 10 (2017)}, no. 4, 893--928.    
    
\bibitem{Wang1} X. Wang. Global infinite energy solutions for the 2D gravity water waves system.{\it Comm. Pure Appl. Math.}
71 (2018), no. 1, 90--162.

%\bibitem{Wang3} X. Wang. Global solution for the 3D gravity water waves system above a flat bottom. {\it arXiv:1508.06227}.

\bibitem{Wu1} S. Wu, Well-posedness in Sobolev spaces of the full water wave problem in 2-D, {\it Invent. Math.}, (130) {\bf 1}(1997), 39--72.

\bibitem{Wu2} S. Wu, Well-posedness in Sobolev spaces of the full
water wave problem in 3-D, {\it J. ~Amer. ~Math. ~Soc.,} {\bf 12}(1999), 445--495.

\bibitem{Wu3} S. Wu, A blow-up criteria and the existence of 2d gravity water waves with angled crests,  arXiv:1502.05342.

\bibitem{Wu4}S. Wu, Almost global well-posedness of the 2-D full water wave problem,  {\it Invent. Math.}, {\bf 177}(2009), 45--135.

\bibitem{Wu5} S. Wu, Global well-posedness of the 3-D full water wave problem,  {\it  Invent. Math.}, {\bf 184}(2011), 125--220.


\bibitem{ZZ} P. Zhang and Z. Zhang, On the free boundary problem of  three-dimensional incompressible Euler equations,
{\it Comm. Pure Appl. Math.,} {\bf 61}(2008), 877--940.


\end{thebibliography}
\end{document}